%% file: hex_patterns.tex
\documentclass[11pt]{article}
\usepackage{graphicx}
\usepackage{epsfig,amsmath,eepic}
\input amssym.def
\input amssym.tex

\newcommand{\qed}{\rule{3mm}{3mm}}
\newcommand{\itbf}{\itshape\bfseries}

\newcommand{\zz}{{\frak z}}
\newcommand{\ee}{{\frak e}}

\newcommand{\cA}{{\cal A}}
\newcommand{\cG}{{\cal G}}
\newcommand{\cH}{{\cal H}}
\newcommand{\cL}{{\cal L}}
\newcommand{\cM}{{\cal M}}
\newcommand{\cT}{{\cal T}}

\newbox\meibox
\def\placeunder#1#2#3#4{\setbox\meibox%
\vbox{\hbox{\hskip#4$\hphantom{#2}$}\hbox{$\hphantom{#1}$}}%
\vtop{\baselineskip=0pt\lineskiplimit=\baselineskip%
\lineskip=#3\hbox to \wd\meibox{\hfil\hskip#4$#2$\hfil}%
\hbox to \wd\meibox{\hfil$#1$\hfil}}}
\def\undertilde#1{\mathchoice{%
\placeunder{\vbox to 1.4pt{\hbox{$\displaystyle\widetilde{\,\,\,
}$}\vss}}{\displaystyle#1}{1.5pt}{1.5pt}}%
{\placeunder{\vbox to 1.4pt{\hbox{$\textstyle\widetilde{\,\,
}$}\vss}}{\textstyle#1}{1.5pt}{1.5pt}}%
{\placeunder{\vbox to 1.4pt{\hbox{$\scriptstyle\tilde{
}$}\vss}}{\scriptstyle#1}{1pt}{1pt}}%
{\placeunder{\vbox to 1.4pt{\hbox{$\scriptscriptstyle\tilde{
}$}\vss}}{\scriptscriptstyle#1}{1pt}{1pt}}%
}
\def\underhat#1{\mathchoice{%
\placeunder{\vbox to 1.4pt{\hbox{$\displaystyle\widehat{\,\,\,
}$}\vss}}{\displaystyle#1}{1.5pt}{1.5pt}}%
{\placeunder{\vbox to 1.4pt{\hbox{$\textstyle\widehat{\,\,
}$}\vss}}{\textstyle#1}{1.5pt}{1.5pt}}%
{\placeunder{\vbox to 1.4pt{\hbox{$\scriptstyle\hat{
}$}\vss}}{\scriptstyle#1}{1pt}{1pt}}%
{\placeunder{\vbox to 1.4pt{\hbox{$\scriptscriptstyle\hat{
}$}\vss}}{\scriptscriptstyle#1}{1pt}{1pt}}%
}
\def\underbar#1{\mathchoice{%
\placeunder{\vbox to 1.4pt{\hbox{$\displaystyle\bar{\,\,\,
}$}\vss}}{\displaystyle#1}{1.5pt}{1.5pt}}%
{\placeunder{\vbox to 1.4pt{\hbox{$\textstyle\bar{\,\,
}$}\vss}}{\textstyle#1}{1.5pt}{1.5pt}}%
{\placeunder{\vbox to 1.4pt{\hbox{$\scriptstyle\bar{
}$}\vss}}{\scriptstyle#1}{1pt}{1pt}}%
{\placeunder{\vbox to 1.4pt{\hbox{$\scriptscriptstyle\bar{
}$}\vss}}{\scriptscriptstyle#1}{1pt}{1pt}}%
}

\newtheorem{theorem}{Theorem}
\newtheorem{proposition}[theorem]{Proposition}
\newtheorem{lemma}[theorem]{Lemma}
\newtheorem{corollary}[theorem]{Corollary}
\newtheorem{definition}[theorem]{Definition}
\newtheorem{conjecture}[theorem]{Conjecture}

\begin{document}

\title{Hexagonal circle patterns and integrable systems:\\
Patterns with the multi-ratio property \\ and Lax equations on the
regular triangular lattice}
\author{A.I.\,Bobenko\footnote{E--mail: {\tt bobenko@math.tu-berlin.de}}
\and T.\,Hoffmann\footnote{E--mail: {\tt
timh@sfb288.math.tu-berlin.de}} \and
Yu.B.\,Suris\footnote{E--mail: {\tt
suris@sfb288.math.tu-berlin.de}}}
\date{Fachbereich Mathematik,
Technische Universit\"at Berlin, Str. 17 Juni 136, 10623 Berlin,
Germany} \maketitle

\begin{abstract}
Hexagonal circle patterns are introduced, and a subclass thereof is
studied in detail. It is characterized by the following property: For
every circle the multi-ratio of its six intersection points with
neighboring circles is equal to $-1$. The relation of such patterns
with an integrable system on the regular triangular lattice is
established. A kind of a B\"acklund transformation for circle patterns
is studied. Further, a class of isomonodromic solutions of the
aforementioned integrable system is introduced, including circle
patterns analogons to the analytic functions $z^\alpha$ and $\log z$. 
\end{abstract}

\newpage

\section{Introduction}
\label{Sect introd}

The theory of circle packings and, more generally, of circle
patterns enjoys in recent years a fast development and a growing
interest of specialists in complex analysis. The origin of this
interest was connected with the Thurston's idea about
approximating the Riemann mapping by circle packings, see
\cite{T1}, \cite{RS}. Since then the theory bifurcated to several
subareas. One of them concentrates around the uniformization
theorem of Koebe--Andreev--Thurston, and is dealing with circle
packing realizations of cell complexes of a prescribed
combinatorics, rigidity properties, constructing hyperbolic 3-manifolds, 
etc \cite{T2}, \cite{MR}, \cite{BS}, \cite{H}.

Another one is mainly dealing with approximation problems, and in
this context it is advantageous to stick from the beginning with
fixed regular combinatorics. The most popular are hexagonal
packings, for which the $C^{\infty}$ convergence to the Riemann
mapping was established by He and Schramm \cite{HS}. Similar
results are available also for circle patterns with the
combinatorics of the square grid introduced by Schramm \cite{S}.
It is also the context of regular patterns (more precisely, the
two just mentioned classes thereof) where some progress was
achieved in constructing discrete analogs of analytic functions
(Doyle's spiralling hexagon packings \cite{BDS} and their
generalizations including the discrete analog of a quotient of
Airy functions \cite{BH}, discrete analogs of ${\rm exp}(z)$ and
${\rm erf}(z)$ for the square grid circle patterns \cite{S},
discrete versions of $z^{\alpha}$ and $\log z$ for the same class
of circle patterns \cite{BP}, \cite{AB}). And it is again the
context of regular patterns where the theory comes into interplay
with the theory of integrable systems. Strictly speaking, only one
instance of such an interplay is well--established up to now:
namely, Schramm's equation describing the square grid circle
packings in terms of M\"obius invariants turns out to coincide
with the stationary Hirota's equation, known to be integrable, see
\cite{BP}, \cite{Z}. It should be said that, generally, the
subject of discrete integrable systems on lattices different from
${\Bbb Z}^n$ is underdeveloped at present. The list of relevant
publications is almost exhausted by \cite{ND}, \cite{NS},
\cite{KN}, \cite{A}, \cite{OP}.

The present paper contributes to several of the above mentioned
issues: we introduce a new interesting class of circle patterns,
and relate them to integrable systems. Besides, for this class we
construct, in parallel to \cite{BP}, \cite{AB}, the analogs of the
analytic functions $z^{\alpha}$, $\log z$.

This class is constituted by {\it hexagonal circle patterns}, or,
in other words, by circle patterns with the combinatorics of the
regular hexagonal lattice (the honeycomb lattice). This means that
each elementary hexagon of the honeycomb lattice corresponds to a
circle, and each common vertex of two hexagons corresponds to an
intersection point of the corresponding circles. In particular,
each circle carries six intersection points with six neighboring
circles. Since at each vertex of the honeycomb lattice there meet
three elementary hexagons, there follows that at each intersection
point there meet three circles.

This class of hexagonal circle patterns is still too wide to be
manageable, but it includes several very interesting subclasses,
leading to integrable systems. For example, one can prescribe
intersection angles of the circles. This situation will be
considered in a subsequent publication. In the present one we
consider the following requirement: the six intersection points on
each circle have the multi-ratio equal to $-1$, where the
multi--ratio is a natural generalization of the notion of a
cross-ratio of four points on a plane.

We show that, adding to the intersection points of the circles
their centers, one embeds hexagonal circle patterns with the
multi-ratio property into an integrable system on the regular
triangular lattice. Each solution of this latter system describes
a peculiar geometrical construction: it consists of three
triangulations of the plane, such that the corresponding
elementary triangles in all three tilings are similar. Moreover,
given one such tiling, one can reconstruct the other two almost
uniquely (up to an affine transformation). If one of the tilings
comes from the hexagonal circle pattern, so do the other two. This
results are contained in Sect. \ref{Sect hex patterns}, \ref{Sect
fgh system}. In the intermediate Sect. \ref{Sect integrable
systems on graphs} we discuss a general notion of integrable
systems on graphs as flat connections with the values in loop
groups. It should be noticed that closely related integrable
equations (albeit on the standard grid ${\Bbb Z}^2$) were
previously introduced by Nijhoff \cite{N} in a totally different
context (discrete Bussinesq equation), see also similar results in
\cite{BK}. However, these results did not go beyond writing down
the equations: geometrical structures behind the equations were
not discussed in these papers.

Having included hexagonal circle patterns with the multi-ratio
property into the framework of the theory of integrable systems,
we get an opportunity of applying the immense machinery of the
latter to studying the properties of the former. This is
illustrated in Sect. \ref{Sect isomonodromic}, \ref{Sect
isomonodromic patterns}, where we introduce and study some
isomonodromic solutions of our integrable system on the triangular
lattice, as well as the corresponding circle patterns. Finally, in
Sect. \ref{Sect hexagonal z^a} we define a subclass of these
``isomonodromic circle patterns'' which are natural discrete
versions of the analytic functions $z^{\alpha}$, $\log z$. The
results of Sect. \ref{Sect isomonodromic}--\ref {Sect hexagonal
z^a} constitute an extension to the present, somewhat more
intricate, situation of the similar constructions for Schramm's
circle patterns with the combinatorics of the square grid
\cite{AB}.

\setcounter{equation}{0}
\section{Hexagonal circle patterns}\label{Sect hex patterns}
\begin{figure}[htbp]
  \begin{center}
    \includegraphics[width=0.3\hsize]{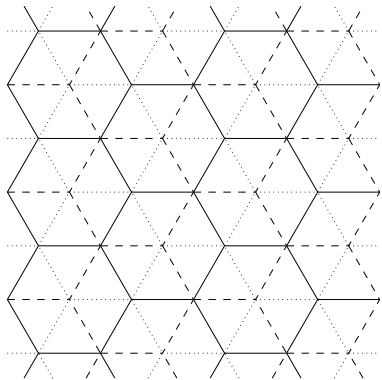}
    \caption{The regular triangular lattice with its hexagonal sublattices.}
    \label{fig:regularLattice}
  \end{center}
\end{figure}
First of all we define the {\itbf regular triangular lattice}
$\cT\cL$ as the cell complex whose vertices are
\begin{equation}
V(\cT\cL)=\Big\{\zz=k+\ell\omega+m\omega^2:\; k,\ell,m\in{\Bbb
Z}\Big\},\quad {\rm where}\quad \omega=\exp(2\pi i/3),
\end{equation}
whose edges are all non--ordered pairs
\begin{equation}
E(\cT\cL)=\Big\{[\zz_1,\zz_2]:\; \zz_1,\zz_2\in V(\cT\cL),\;
|\zz_1-\zz_2|=1 \Big\},
\end{equation}
and whose 2-cells are all regular triangles with the vertices in
$V(\cT\cL)$ and the edges in $E(\cT\cL)$. We shall use triples
$(k,\ell,m)\in {\Bbb Z}^3$ as coordinates of the vertices of the
regular triangular lattice, identifying two such triples iff they
differ by the vector $(n,n,n)$ with $n\in{\Bbb Z}$. We call two
points $\zz_1,\zz_2$ {\it neighbors in} $\cT\cL$, iff
$[\zz_1,\zz_2]\in E(\cT\cL)$.

To the complex $\cT\cL$ there correspond three {\itbf regular
hexagonal sublattices} $\cH\cL_j$, $j=0,1,2$. Each $\cH\cL_j$ is
the cell complex whose vertices are
\begin{equation}
V(\cH\cL_j)=\Big\{\zz=k+\ell\omega+m\omega^2:\; k,\ell,m\in{\Bbb
Z},\; k+\ell+m\not\equiv j\!\!\pmod 3\Big\},
\end{equation}
whose edges are
\begin{equation}
E(\cH\cL_j)=\Big\{[\zz_1,\zz_2]:\;\zz_1,\zz_2\in V(\cH\cL_j),\;
|\zz_1-\zz_2|=1 \Big\},
\end{equation}
and whose 2-cells are all regular hexagons with the vertices in
$V(\cH\cL_j)$ and the edges in $E(\cH\cL_j)$. Again, we call two
points $\zz_1,\zz_2$ {\it neighbors in} $\cH\cL_j$, iff
$[\zz_1,\zz_2]\in E(\cH\cL_j)$. Obviously, every point in
$V(\cH\cL_j)$ has three neighbors in $\cH\cL_j$, as well as three
neighbors in $\cT\cL$ which do not belong to $V(\cH\cL_j)$. The
centers of 2-cells of $\cH\cL_j$ are exactly the points of
$V(\cT\cL)\setminus V(\cH\cL_j)$, i.e. the points
$\zz'=k+\ell\omega+m\omega^2$ with $k+\ell+m\equiv j\!\!\pmod 3$.

In the following definition we consider only $\cH\cL_0$, since,
clearly, $\cH\cL_1$ and $\cH\cL_2$ are obtained from $\cH\cL_0$
via shifting all the corresponding objects by $\omega$, resp. by
$\omega^2$.

\begin{definition} \label{def hex pattern}
We say that a map $w:V(\cH\cL_0)\mapsto\hat{\Bbb C}$ defines a
{\itbf hexagonal circle pattern}, if the following condition is
satisfied:
\begin{itemize}
\item Let
\[
\zz_k=\zz'+\varepsilon^k\in V(\cH\cL_0), \quad
k=1,2,\ldots,6,\quad  where \quad\varepsilon=\exp(\pi i/3),
\]
 be the vertices of any elementary hexagon in $\cH\cL_0$ with the center
 $\zz'\in V(\cT\cL)\setminus V(\cH\cL_0)$. Then the points
 $w(\zz_1),w(\zz_2),\ldots,w(\zz_6)\in\hat{\Bbb C}$ lie on a circle, and their
 circular order is just the listed one. We denote the circle through
 the points $w(\zz_1),w(\zz_2),\ldots,w(\zz_6)$ by $C(\zz')$, thus putting it
 into a correspondence with the center $\zz'$ of the elementary hexagon above.
 \end{itemize}
\end{definition}
As a consequence of this condition, we see that if two elementary
hexagons of $\cH\cL_0$ with the centers in $\zz',\zz''\in
V(\cT\cL)\setminus V(\cH\cL_0)$ have a common edge
$[\zz_1,\zz_2]\in E(\cH\cL_0)$, then the circles $C(\zz')$ and
$C(\zz'')$ intersect in the points $w(\zz_1)$, $w(\zz_2)$.
Similarly, if three elementary hexagons of $\cH\cL_0$ with the
centers in $\zz',\zz'',\zz''' \in V(\cT\cL)\setminus V(\cH\cL_0)$
meet in one point $\zz_0\in V(\cH\cL_0)$, then the circles
$C(\zz')$, $C(\zz'')$ and $C(\zz''')$ also have a common
intersection point $w(\zz_0)$. (Note that in every point $\zz_0\in
V(\cH\cL_0)$ there  meet three distinct elementary hexagons of
$\cH\cL_0$). \vspace{2mm}

{\bf Remark.} Sometimes it will be convenient to consider circle
patterns defined not on the whole of $\cH\cL_0$, but rather on
some connected subgraph of the regular hexagonal lattice.
\vspace{2mm}

We shall study in this paper a subclass of hexagonal circle
patterns satisfying an additional condition. We need the following
generalization of the notion of cross-ratio.
\begin{definition}
Given a $(2p)$-tuple $(w_1,w_2,\ldots,w_{2p})\in{\Bbb C}^{2p}$ of
complex numbers, their {\itbf multi-ratio} is the following
number:
\begin{equation}
M(w_1,w_2,\ldots,w_{2p})=\frac{\prod_{j=1}^p(w_{2j-1}-w_{2j})}
{\prod_{j=1}^p (w_{2j}-w_{2j+1})},
\end{equation}
where it is agreed that $w_{2p+1}=w_1$.
\end{definition}
In particular,
\[
M(w_1,w_2,w_3,w_4)=\frac{(w_1-w_2)(w_3-w_4)}{(w_2-w_3)(w_4-w_1)}
\]
is the usual cross-ratio, while in the present paper we shall be
mainly dealing with
\[
M(w_1,w_2,\ldots,w_6)=\frac{(w_1-w_2)(w_3-w_4)(w_5-w_6)}
{(w_2-w_3)(w_4-w_5)(w_6-w_1)}.
\]
The following two obvious properties of the multi-ratio will be
important for us:
\begin{itemize}
\item[(i)] The multi-ratio $M(w_1,w_2,\ldots,w_{2p})$ is invariant with
respect to the action of an arbitrary M\"obius transformation
$w\mapsto (aw+b)/(cw+d)$ on all of its arguments.
\item[(ii)] The multi-ratio $M(w_1,w_2,\ldots,w_{2p})$ is a M\"obius
transformation with respect to each one of its arguments.
\end{itemize}
We shall need also the following, slightly less obvious, property:
\begin{itemize}
\item[(iii)] If the points $w_1,w_2,\ldots,w_{2p-1}$ lie on a circle $C\subset
\hat{\Bbb C}$, and the multi-ratio $M(w_1,w_2,\ldots,w_{2p})$ is
real, then also $w_{2p}\in C$.
\end{itemize}

\begin{definition}
We say that a map $w:V(\cH\cL_0)\mapsto\hat{\Bbb C}$ defines a
{\itbf hexagonal circle pattern with}
$\boldsymbol{M}\boldsymbol{R}\boldsymbol{=}
\boldsymbol{-}\boldsymbol{1}$, if in addition to the condition of
Definition \ref{def hex pattern} the following one is satisfied:
\begin{itemize}
\item For any elementary hexagon in $\cH\cL_0$ with the vertices
 $\zz_1,\zz_2,\ldots,\zz_6\in V(\cH\cL_0)$ (listed counterclockwise), the
 multi-ratio
\begin{equation}\label{spec cond}
 M(w_1,w_2,\ldots,w_6)=-1,
\end{equation}
where $w_k=w(\zz_k)$.
\end{itemize}
\end{definition}
Geometrically the condition (\ref{spec cond}) means that, first,
the lengths of the sides of the hexagon with the vertices
$w_1w_2\ldots w_6$ satisfy the condition
\[
|w_1-w_2|\cdot|w_3-w_4|\cdot|w_5-w_6|=
|w_2-w_3|\cdot|w_4-w_5|\cdot|w_6-w_1|,
\]
and, second, that the sum of the angles of the hexagon at the
vertices $w_1$, $w_3$, and $w_5$ is equal to $2\pi\!\!\pmod
{2\pi}$, as well as the sum of the angles at the vertices $w_2$,
$w_4$, and $w_6$. Notice that if a hexagon is inscribed in a
circle and satisfies (\ref{spec cond}), then it is {\it
conformally symmetric}, i.e. there exists a M\"obius
transformation mapping it onto a centrally symmetric hexagon.
Notice also that the regular hexagons satisfy this condition.

To demonstrate quickly the {\it existence} of hexagonal circle
patterns with $MR=-1$ we give their {\it construction} via solving
a suitable Cauchy problem.
\begin{lemma}
Consider a row of elementary hexagons of $\cH\cL_0$ running from
the north--west to the south-east, with the centers in the points
$\zz'_k=k-k\omega$. Let the map $w$ be defined in five vertices of
each hexagon -- in all except $\zz'_k+\varepsilon$. Suppose that
the five points $w(\zz'_k+\varepsilon^j)$, $j=2,3,\ldots,6$, lie
on the circles $C(\zz'_k)$. These data determine uniquely a map
$w:V(\cH\cL_0)\mapsto\hat{\Bbb C}$ yielding a hexagonal circle
pattern with $MR=-1$ on the whole lattice.
\end{lemma}
{\bf Proof.} Equation (\ref{spec cond}) determines the points
$w(\zz'_k+\varepsilon)$, which, according to the property above,
lie also on $C(\zz'_k)$. Now for every hexagon of the parallel row
next to north--east, with the centers in the points
$\zz''_k=\zz'_k+1+\varepsilon=(k+2)-(k-1)\omega$, we know the
value of the map $w$ in three vertices, namely in
\[
\zz''_k+\varepsilon^4=\zz'_k+1=\zz'_{k+1}+\varepsilon^2,\quad
\zz''_k+\varepsilon^3=\zz'_{k}+\varepsilon,  \quad
\zz''_k+\varepsilon^5=\zz'_{k+1}+\varepsilon^2.
\]
This uniquely defines the circle $C(\zz''_k)$, as the only circle
through three points $w(\zz''_k+\varepsilon^3)$,
$w(\zz''_k+\varepsilon^4)$ and $w(\zz''_k+\varepsilon^5)$. The
intersection points of these circles of the second row give us the
values of the map $w$ in the points $\zz''_k+\varepsilon^2$ and
$\zz''_k+\varepsilon^6$. Namely, $w(\zz''_k+\varepsilon^2)$ is the
intersection point of $C(\zz''_k)$ with $C(\zz''_{k-1})$,
different from $w(\zz''_k+\varepsilon^3)$, and
$w(\zz''_k+\varepsilon^6)$ is the intersection point of
$C(\zz''_k)$ with $C(\zz''_{k+1})$, different from
$w(\zz''_k+\varepsilon^5)$. Therefore we get the values of the map
$w$ in five vertices of each hexagon of the next parallel row --
in all except $\zz''_k+\varepsilon$. The induction allows to
continue the construction {\it ad infinitum}. \qed \vspace{3mm}

Now we show that, adding the centers of the circles of a hexagonal
pattern with $MR=-1$ to their intersection points, we come to a
new interesting notion.
\begin{theorem}\label{central extension}
Let the map $w:V(\cH\cL_0)\mapsto\hat{\Bbb C}$ define a hexagonal
circle pattern with $MR=-1$. Extend $w$ to the points of
$V(\cT\cL)\setminus V(\cH\cL_0)$ by the following rule. Fix some
point $P_{\infty}\in\hat{\Bbb C}$. Let $\zz'$ be a center of an
elementary hexagon of $\cH\cL_0$. Set $w(\zz')$ to be the
reflection of the point $P_{\infty}$ in the circle $C(\zz')$. Then
the condition (\ref{spec cond}) holds also for $w_k=w(\zz_k)$ in
the case when the points $\zz_1,\zz_2,\ldots,\zz_6$ are the
vertices of any elementary hexagon of the two complementary
hexagonal sublattices $\cH\cL_1$ and $\cH\cL_2$.
\end{theorem}
\begin{figure}[htbp]
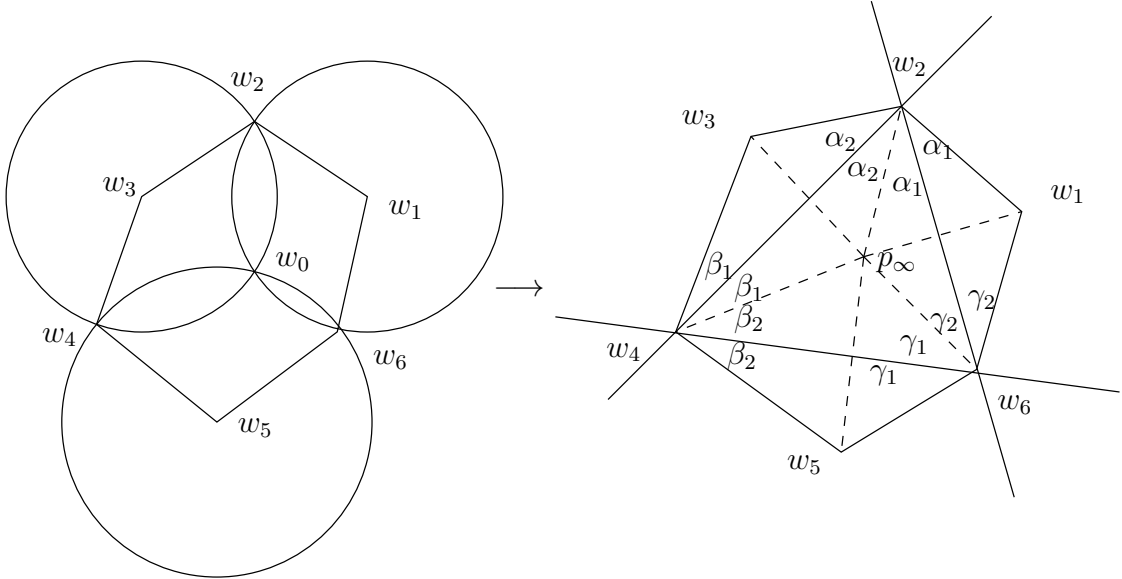

  \begin{center}
    \input a_proofMV.tex
    \caption{An elementary hexagon with its center point sent to $\infty$.}
    \label{fig:proof5}
  \end{center}
\end{figure}
{\bf Proof.} Consider the situation corresponding to an elementary
hexagon of the sublattice $\cH\cL_1$ or $\cH\cL_2$ (see
Fig.~\ref{fig:proof5}). The point $w_0$ is the intersection point
of the three circles $C(\zz_1)$, $C(\zz_3)$, and $C(\zz_5)$, the
points $w_1$, $w_3$, and $w_5$ are obtained by reflection of
$P_{\infty}$ in the corresponding circles, and the points $w_2$,
$w_4$, and $w_6$ are the pairwise intersection points of these
circles different from $w_0$. To simplify the geometry behind this
situation, perform a M\"obius transformation sending $w_0$ to
infinity. Then the circles $C(\zz_1)$, $C(\zz_3)$, and $C(\zz_5)$
become straight lines, and the points $w_1$, $w_3$, $w_5$ are the
reflections of $P_{\infty}$ in these lines (see
Fig.~\ref{fig:proof5}; for definiteness we suppose here that the
M\"obius image of $P_{\infty}$ lies in the interior of the
triangle formed by these straight lines). By construction, one
gets:
\[
|w_2-w_1|=|w_2-w_3|, \quad |w_4-w_3|=|w_4-w_5|, \quad
|w_6-w_5|=|w_6-w_1|;
\]
the angles by the vertices $w_2$, $w_4$, $w_6$ are equal to
$2(\alpha_1+\alpha_2)$, $2(\beta_1+\beta_2)$,
$2(\gamma_1+\gamma_2)$, respectively, so that their sum is equal
to
\[
2(\alpha_1+\alpha_2+\beta_1+\beta_2+\gamma_1+\gamma_2)=2\pi;
\]
the angles by the vertices $w_1$, $w_3$, $w_5$ are equal to
$\pi-(\alpha_1+\gamma_2)$, $\pi-(\beta_1+\alpha_2)$,
$\pi-(\gamma_1+\beta_2)$, respectively, so that their sum is equal
to
\[
3\pi-(\alpha_1+\alpha_2+\beta_1+\beta_2+\gamma_1+\gamma_2)=2\pi.
\]
This proves that the hexagon under consideration satisfies
(\ref{spec cond}).
 \qed \vspace{3mm}

A particular case of the construction of Theorem \ref{central
extension} is when $P_{\infty}=\infty$, so that the map $w$ is
extended by the {\it centers} of the corresponding circles. In any
case, this theorem suggests to consider the class of maps
described in the following definition.
\begin{definition}
We say that the map $w:V(\cT\cL)\mapsto\hat{\Bbb C}$ defines
{\itbf a triangular lattice with
$\boldsymbol{M}\boldsymbol{R}\boldsymbol{=}
\boldsymbol{-}\boldsymbol{1}$}, if the equation (\ref{spec cond})
holds for $w_k=w(\zz_k)$, whenever the points
$\zz_1,\zz_2,\ldots,\zz_6$ are the vertices (listed
counterclockwise) of any elementary hexagon of any of the
sublattices $\cH\cL_j$ $(j=0,1,2)$.
\end{definition}

In the next section we shall discuss an integrable system on the
regular triangular lattice, each solution of which delivers, in a
single construction, {\it three} different triangular lattices
with $MR=-1$. However, these three lattices are not independent:
given such a lattice, the two associated ones can be constructed
almost uniquely (up to an affine transformation $w\mapsto aw+b$).
It will turn out that if the original lattice comes from a
hexagonal circle pattern with $MR=-1$, then the two associated
ones do likewise.

\setcounter{equation}{0}
\section{Discrete flat connections on graphs}
\label{Sect integrable systems on graphs}

Let us describe a general construction of ``integrable systems''
on graphs which does not hang on the specific features of the
regular triangular lattice. This notion includes the following
ingredients:
\begin{itemize}
\item An {\it oriented graph} $\cG$; the set of its vertices will
be denoted $V(\cG)$, the set of its edges will be denoted
$E(\cG)$.
\item A loop group $G[\lambda]$, whose elements are functions from
${\Bbb C}$ into some group $G$. The complex argument $\lambda$ of
these functions is known in the theory of integrable systems as
the ``spectral parameter''.
\item A ``wave function'' $\Psi: V(\cG)\mapsto G[\lambda]$,
defined on the vertices of $\cG$.
\item A collection of ``transition matrices'' $L: E(\cG)\mapsto G[\lambda]$
defined on the edges of $\cG$.
\end{itemize}
It is supposed that for any oriented edge $\ee=(\zz_1,\zz_2)\in
E(\cG)$ the values of the wave functions in its ends are connected
via
\begin{equation}\label{wave function evol}
\Psi(\zz_2,\lambda)=L(\ee,\lambda)\Psi(\zz_1,\lambda).
\end{equation}
Therefore the following {\it discrete zero curvature condition} is
supposed to be satisfied. Consider any closed contour consisting
of a finite number of edges of $\cG$:
\[
\ee_1=(\zz_1,\zz_2),\quad \ee_2=(\zz_2,\zz_3),\quad \ldots,\quad
\ee_p=(\zz_p,\zz_1).
\]
Then
\begin{equation}\label{zero curv cond}
L(\ee_p,\lambda)\cdots L(\ee_2,\lambda)L(\ee_1,\lambda)=I.
\end{equation}
In particular, for any edge $\ee=(\zz_1,\zz_2)$, if
$\ee^{-1}=(\zz_2,\zz_1)$, then
\begin{equation}\label{zero curv cond inv}
L(\ee^{-1},\lambda)=\Big(L(\ee,\lambda)\Big)^{-1}.
\end{equation}

Actually, in applications the matrices $L(\ee,\lambda)$ depend
also on a point of some set $X$ (the ``phase space'' of an
integrable system), so that some elements $x(\ee)\in X$ are
attached to the edges $\ee$ of $\cG$. In this case the discrete
zero curvature condition (\ref{zero curv cond}) becomes equivalent
to the collection of equations relating the fields $x(\ee_1)$,
$\ldots$, $x(\ee_p)$ attached to the edges of each closed contour.
We say that this collection of equations admits a {\it zero
curvature representation}.

For an arbitrary graph, the analytical consequences of the zero
curvature representation for a given collection of equations are
not clear. However, in case of regular lattices, like $\cT\cL$,
such representation may be used to determine conserved quantities
for suitably defined Cauchy problems, as well as to apply powerful
analytical methods for finding concrete solutions. \vspace{3mm}

{\bf Remark.} The above construction of integrable systems on
graphs is not the only possible one. For example, in the
construction by Adler \cite{A} the fields are defined on the
vertices of a planar graph, and the equations relate the fields on
{\it stars} consisting  of the edges incident to each single
vertex, rather than the fields on closed contours. Examples
are given by discrete time systems of the relativistic Toda type.
In the corresponding zero curvature representation the wave
functions $\Psi$ naturally live on 2-cells rather than on
vertices. The transition matrices live on edges: the matrix
$L(\ee,\lambda)$ corresponds to the transition {\it across} $\ee$
and depends on the fields sitting on two ends of $\ee$.

\setcounter{equation}{0}
\section{An integrable system on the regular triangular lattice}
\label{Sect fgh system}

We now introduce an {\it orientation} of the edges of the regular
triangular lattice $\cT\cL$. Namely, we declare as positively
oriented all edges of the types
\[
(\zz,\zz+1),\quad (\zz,\zz+\omega), \quad (\zz,\zz+\omega^2).
\]
Correspondingly, all edges of the types
\[
(\zz,\zz-1),\quad (\zz,\zz-\omega), \quad (\zz,\zz-\omega^2)
\]
are negatively oriented. Thus all elementary triangles
become oriented. There are two types of elementary triangles:
those ``pointing upwards'' $(\zz,\zz+\omega,\zz-1)$ are oriented
counterclockwise, while those ``pointing downwards''
$(\zz,\zz+\omega^2,\zz-1)$ are oriented clockwise.

\subsection{Lax representation}

The group $G[\lambda]$ we use in our construction is the {\it
twisted loop group} over ${\rm SL}(3,{\Bbb C})$:
\begin{equation}\label{loop group}
\Big\{L:{\Bbb C}\mapsto {\rm SL}(3,{\Bbb C})\Big|\;
L(\omega\lambda)=\Omega L(\lambda)\Omega^{-1}\Big\},
\end{equation}
where $\Omega={\rm diag}(1,\omega,\omega^2)$. The elements of
$G[\lambda]$ we attach to every {\it positively oriented} edge of
$\cT\cL$ are of the form
\begin{equation}\label{L}
L(\lambda)=(1+\lambda^3)^{-1/3}\left(\begin{array}{ccc} 1 &
\lambda f & 0 \\ 0 & 1 & \lambda g \\\lambda h & 0 &
1\end{array}\right), \quad fgh=1.
\end{equation}
Hence, to each positively oriented edge we assign a triple of
complex numbers $(f,g,h)\in{\Bbb C}^3$ satisfying an additional
condition $fgh=1$. In other words, choosing $(f,g)$ (say) as the
basic variables, we can assume that the ``phase space'' $X$
mentioned in the previous section, is ${\Bbb C}_{*}\times {\Bbb
C}_{*}$. The scalar factor $(1+\lambda^3)^{-1/3}$ is not very
essential and assures merely that $\det L(\lambda)=1$.

It is obvious that the zero curvature condition (\ref{zero curv
cond}) is fulfilled for every closed contour in $\cT\cL$, if and
only if it holds for all elementary triangles.
\begin{theorem}\label{Equations of motion}
Let $\ee_1$, $\ee_2$, $\ee_3$ be the consecutive positively
oriented edges of an elementary triangle of $\cT\cL$. Then the
zero curvature condition
\[
L(\ee_3,\lambda)L(\ee_2,\lambda)L(\ee_1,\lambda)=I
\]
is equivalent to the following set of equations:
\begin{equation}\label{fields fact}
f_1+f_2+f_3=0,\qquad g_1+g_2+g_3=0,
\end{equation}
and
\begin{equation}\label{motion eq}
f_1g_1=f_3g_2\quad\Leftrightarrow\quad f_2g_2=f_1g_3
\quad\Leftrightarrow\quad f_3g_3=f_2g_1,
\end{equation}
with the understanding that $h_k=(f_kg_k)^{-1}$, $k=1,2,3$.
\end{theorem}
{\bf Proof.} An easy calculation shows that the matrix equation
$L_3L_2L_1=I$ consists of the following nine scalar equations:
\begin{equation}\label{motion eq aux1}
f_1+f_2+f_3=0,\qquad g_1+g_2+g_3=0,\qquad h_1+h_2+h_3=0,
\end{equation}
\begin{equation}\label{motion eq aux2}
f_3g_2h_1=1,\qquad g_3h_2f_1=1, \qquad h_3f_2g_1=1,
\end{equation}
\begin{equation}\label{motion eq aux3}
f_3g_2+f_3g_1+f_2g_1=0,\qquad g_3h_2+g_3h_1+g_2h_1=0,\qquad
h_3f_2+ h_3f_1+h_2f_1=0.
\end{equation}
It remains to isolate the independent ones among these nine
equations. First of all, equations (\ref{motion eq aux3}) are
equivalent to (\ref{motion eq aux2}), provided (\ref{motion eq
aux1}) and $f_kg_kh_k=1$ hold. For example:
\[
f_3(g_2+g_1)+f_2g_1=0\;\Leftrightarrow\;f_3g_3=f_2g_1\;\Leftrightarrow\;
h_3f_2g_1=1.
\]
Next, the conditions $f_kg_kh_k=1$ allow us to rewrite
(\ref{motion eq aux2}) as
\begin{equation}\label{motion eq aux4}
f_1g_1=f_3g_2,\qquad f_2g_2=f_1g_3, \qquad f_3g_3=f_2g_1.
\end{equation}
Further, all equations in (\ref{motion eq aux4}) are equivalent
provided (\ref{fields fact}) holds. For example:
\[
f_1g_1=f_3g_2\;\Rightarrow\;
(f_2+f_3)g_1=f_3(g_1+g_3)\;\Rightarrow f_2g_1=f_3g_3.
\]
Finally, $h_1+h_2+h_3=0$ follows from (\ref{fields fact}),
(\ref{motion eq}). Indeed,
\begin{eqnarray*}
h_1+h_2 & = & (f_1g_1)^{-1}+(f_2g_2)^{-1}=
(f_3g_2)^{-1}+(f_2g_2)^{-1}\\
& = & (f_2g_2)^{-1}(f_2+f_3)f_3^{-1}=-(f_2g_2)^{-1}f_1f_3^{-1}\\
& = & -(f_1g_3)^{-1}f_1f_3^{-1}=-(f_3g_3)^{-1}=-h_3.
\end{eqnarray*}
The theorem is proved. For want of a better name we shall call the
system of equations (\ref{fields fact}), (\ref{motion eq}) the
{\itbf fgh--system}. \qed \vspace{3mm}

The equations (\ref{motion eq aux1}) may be interpreted in the
following way: there exist functions $u,v,w:V(\cT\cL)\mapsto{\Bbb
C}$ such that for any positively oriented edge $\ee=(\zz_1,\zz_2)$
there holds:
\begin{equation}\label{fact}
f(\ee)=u(\zz_2)-u(\zz_1),\quad g(\ee)=v(\zz_2)-v(\zz_1),\quad
h(\ee)=w(\zz_2)-w(\zz_1).
\end{equation}
The function $u$ is determined by $f$ uniquely, up to an additive
constant, and similarly for the functions $v$, $w$. Having
introduced functions $u,v,w$ sitting in the vertices of $\cT\cL$,
we may reformulate the remaining equations (\ref{motion eq}) as
follows: let $\zz_1,\zz_2,\zz_3$ be the consecutive vertices of a
positively oriented elementary triangle, then
\begin{equation}\label{motion eq zw}
\frac{u(\zz_2)-u(\zz_1)}{u(\zz_3)-u(\zz_2)}=
\frac{v(\zz_3)-v(\zz_2)}{v(\zz_1)-v(\zz_3)}.
\end{equation}
The equations arising by cyclic permutations of indices
$(1,2,3)\mapsto(2,3,1)$ are equivalent to this one due to
(\ref{motion eq}). So, we have one equation pro elementary
triangle $\zz_1\zz_2\zz_3$. Its geometrical meaning is the
following: the triangle $u(\zz_1)u(\zz_2)u(\zz_3)$ is similar to
the triangle $v(\zz_2) v(\zz_3)v(\zz_1)$ (where the corresponding
vertices are listed on the corresponding places). Of course, these
two triangles are also similar to the third one,
$w(\zz_3)w(\zz_1)w(\zz_2)$.

\subsection{Cauchy problem}
We discuss now the Cauchy data which allow one to determine a
solution of the $fgh$--system. The key observation is the
following.
\begin{lemma}\label{lemma 4th point}
Given the values of two fields, say $u$ and $v$, in three points
$\zz_0$, $\zz_1=\zz_0+1$ and $\zz_2=\zz_0+\omega$, the equations
of the $fgh$--system determine uniquely the values of $u$ and $v$
in the point $\zz_3=\zz_0+1+\omega$:
\begin{equation}\label{induct aux1}
u_3-u_0=(u_1-u_0)\frac{v_1-v_0}{v_1-v_2}+(u_2-u_0)\frac{v_2-v_0}{v_2-v_1},
\end{equation}
\begin{equation}\label{induct aux2}
v_3-v_1=(v_1-v_0)\frac{u_1-u_0}{u_0-u_3}\quad \Leftrightarrow\quad
v_3-v_2=(v_2-v_0)\frac{u_2-u_0}{u_0-u_3}.
\end{equation}
\end{lemma}
{\bf Proof.} The formula (\ref{induct aux1}) follows by
eliminating $v_3$ from
\begin{equation}\label{induct aux0}
\frac{u_0-u_3}{u_1-u_0}=\frac{v_1-v_0}{v_3-v_1},\qquad
\frac{u_0-u_3}{u_2-u_0}=\frac{v_2-v_0}{v_3-v_2}.
\end{equation}
These equations yield then (\ref{induct aux2}). \qed \vspace{2mm}

This immediately yields the following statement.
\begin{proposition}\label{Cauchy data for fgh system}
\begin{itemize}
\item[{\rm a)}] The values of the fields $u$ and $v$ in the vertices of the
zig--zag line running from the north--west to the south--east,
\[
\Big\{\zz=k+\ell\omega: k+\ell=0,1\Big\},
\]
uniquely determine the functions $u,v:V(\cT\cL)\mapsto{\Bbb C}$ on
the whole lattice.
\item[{\rm b)}] The values of the fields $u$ and $v$ on the two positive
semi-axes,
\[
\Big\{\zz=k: k\ge 0\Big\}\cup\Big\{\zz=\ell\omega: \ell\ge
0\Big\},
\]
uniquely determine the functions $u,v$ on the whole sector
\[
\Big\{\zz=k+\ell\omega: k,\ell\ge 0\Big\}= \Big\{\zz\in V(\cT\cL):
0\le{\rm\arg}(\zz)\le2\pi/3\Big\}.
\]
\end{itemize}
\end{proposition}
{\bf Proof} follows by induction with the help of the formulas
(\ref{induct aux1}), (\ref{induct aux2}). \qed

\subsection{Sym formula and related results}

There holds the following result having many analogs in the
differential geometry described by integrable systems (``Sym
formula'', see, e.g., \cite{BP}).
\begin{proposition}\label{Sym formula}
Let $\Psi(\zz,\lambda)$ be the solution of (\ref{wave function
evol}) with the initial condition $\Psi(\zz_0,\lambda)=I$ for some
$\zz_0\in V(\cT\cL)$. Then the fields $u,v,w$ may be found as
\begin{equation}\label{Sym}
\left.\frac{d\Psi}{d\lambda}\right|_{\lambda=0}=
\left(\begin{array}{ccc} 0 & u & 0 \\ 0 & 0 & v \\ w & 0 & 0
\end{array} \right).
\end{equation}
\end{proposition}
{\bf Proof.} Note, first of all, that from $\Psi(\zz_0,0)=I$ and
$L(\ee,0)=I$ there follows that $\Psi(\zz,0)=I$ for all $\zz\in
V(\cT\cL)$. Consider an arbitrary positively oriented edge
$\ee=(\zz_1,\zz_2)$. From (\ref{wave function evol}) there
follows:
\[
\frac{d\Psi(\zz_2)}{d\lambda}-\frac{d\Psi(\zz_1)}{d\lambda} =
\left(\frac{dL(\ee)}{d\lambda}\Psi(\zz_1)+
L(\ee)\frac{d\Psi(\zz_1)}{d\lambda}\right)-\frac{d\Psi(\zz_1)}{d\lambda}
\]
At $\lambda=0$ we find:
\begin{eqnarray*}
\lefteqn{\left.\frac{d\Psi(\zz_2)}{d\lambda}\right|_{\lambda=0}
-\left.\frac{d\Psi(\zz_1)}{d\lambda}\right|_{\lambda=0}=
\left.\frac{dL(\ee)}{d\lambda}\right|_{\lambda=0}=
\left(\begin{array}{ccc} 0 & f(\ee) & 0 \\ 0 & 0 & g(\ee) \\
h(\ee) & 0 & 0
\end{array}\right)}\\
 & = & \left(\begin{array}{ccc} 0 & u(\zz_2)-u(\zz_1) & 0 \\
0 & 0 & v(\zz_2)-v(\zz_1) \\ w(\zz_2)-w(\zz_1) & 0 & 0
\end{array}\right).
\end{eqnarray*}
This proves the Proposition. \qed \vspace{2mm}

Next terms of the power series expansion of the wave function
$\Psi(\zz,\lambda)$ around $\lambda=0$ also deliver interesting
and important results.
\begin{proposition}\label{Closed forms}
Let $\Psi(\zz,\lambda)$ be the solution of (\ref{wave function
evol}) with the initial condition $\Psi(\zz_0,\lambda)=I$ for some
$\zz_0\in V(\cT\cL)$. Then
\begin{equation}\label{Sym2}
\frac{1}{2}\,\left.\frac{d^2\Psi}{d\lambda^2}\right|_{\lambda=0}=
\left(\begin{array}{ccc} 0 & 0 & a \\ b & 0 & 0 \\ 0 & c & 0
\end{array} \right),
\end{equation}
where the function $a:V(\cT\cL)\mapsto{\Bbb C}$ satisfies the
difference equation
\begin{equation}\label{eq for a}
a(\zz_2)-a(\zz_1)=v(\zz_1)\Big(u(\zz_2)-u(\zz_1)\Big),
\end{equation}
and similar equations hold for the functions
$b,c:V(\cT\cL)\mapsto{\Bbb C}$ (with the cyclic permutation
$(u,v,w)\mapsto(w,u,v)$).
\end{proposition}
{\bf Proof.} Proceeding as in the proof of Proposition \ref{Sym
formula}, we have:
\[
\frac{d^2\Psi(\zz_2)}{d\lambda^2}-\frac{d^2\Psi(\zz_1)}{d\lambda^2}
= \left(\frac{d^2L(\ee)}{d\lambda^2}\Psi(\zz_1)+
2\frac{dL(\ee)}{d\lambda}\frac{d\Psi(\zz_1)}{d\lambda}+
L(\ee)\frac{d^2\Psi(\zz_1)}{d\lambda^2}\right)-
\frac{d^2\Psi(\zz_1)}{d\lambda^2}
\]
Taking into account that $d^2L(\ee)/d\lambda^2|_{\lambda=0}=0$, we
find at $\lambda=0$:
\begin{eqnarray*}
\lefteqn{\left.\frac{d^2\Psi(\zz_2)}{d\lambda^2}\right|_{\lambda=0}
-\left.\frac{d^2\Psi(\zz_1)}{d\lambda^2}\right|_{\lambda=0}=
2\left.\frac{dL(\ee)}{d\lambda}\right|_{\lambda=0}
\left.\frac{d\Psi(\zz_1)}{d\lambda}\right|_{\lambda=0}=}\\ \nonumber\\
& = & 2\left(\begin{array}{ccc} 0 & f(\ee) & 0 \\ 0 & 0 & g(\ee)
\\ h(\ee) & 0 & 0
\end{array}\right)\left(\begin{array}{ccc} 0 & u(\zz_1) & 0 \\
0 & 0 & v(\zz_1) \\ w(\zz_1) & 0 & 0 \end{array}\right).
\end{eqnarray*}
This implies the statement of the proposition. \qed \vspace{2mm}

Notice that it is {\it \`a priori} not obvious that the equation
(\ref{eq for a}) admits a well--defined solution on $V(\cT\cL)$,
or, in other words, that its right--hand side defines a closed
form on $\cT\cL$. This fact might be proved by a direct
calculation, based upon the equations of the $fgh$--system, but
the above argument gives a more conceptual and a much shorter
proof.
\begin{corollary}
Under the conditions of Propositions \ref{Sym formula},\ref{Closed
forms}, we have:
\begin{equation}\label{Sym12}
-\frac{1}{2}\,\left.\frac{d^2\Psi}{d\lambda^2}\right|_{\lambda=0}+
\left(\frac{d\Psi}{d\lambda}\right)^2_{\lambda=0}=
\left(\begin{array}{ccc} 0 & 0 & a' \\ b' & 0 & 0 \\ 0 & c' & 0
\end{array} \right),
\end{equation}
where the function $a':V(\cT\cL)\mapsto{\Bbb C}$ satisfies the
difference equation
\begin{equation}\label{eq for a'}
a'(\zz_2)-a'(\zz_1)=u(\zz_2)\Big(v(\zz_2)-v(\zz_1)\Big),
\end{equation}
and similar equations hold for the functions
$b',c':V(\cT\cL)\mapsto{\Bbb C}$ (with the cyclic permutation
$(u,v,w)\mapsto(w,u,v)$).
\end{corollary}
Further examples of such exact forms may be obtained from the
values of higher derivatives of the wave function
$\Psi(\zz,\lambda)$ at $\lambda=0$.

\subsection{One--field equations}

We discuss now the equations satisfied by the field $u$ alone, as
well as by the field $v$ alone. In this point we make contact with
the geometric considerations of Sect. \ref{Sect hex patterns}.
\begin{theorem}\label{z to w}
\begin{enumerate}
\item Both maps $u,v:V(\cT\cL)\mapsto {\Bbb C}$ define triangular lattices
with $MR=-1$. In other words, if $\zz_1,\zz_2,\ldots,\zz_6$ are
the vertices (listed counterclockwise) of any elementary hexagon
of any of the hexagonal sublattices $\cH\cL_j$ $(j=0,1,2)$, and if
$u_k=u(\zz_k)$ and $v_k=v(\zz_k)$, then there hold both the
equations
\begin{equation}\label{hex eq z}
M(u_1,u_2,\ldots,u_6)=-1
\end{equation}
and
\begin{equation}\label{hex eq w}
M(v_1,v_2,\ldots,v_6)=-1.
\end{equation}
\item Given a triangular lattice $u:V(\cT\cL)\mapsto{\Bbb C}$ with $MR=-1$,
there exists a unique, up to an affine transformation $v\mapsto
av+b$, function $v:V(\cT\cL)\mapsto{\Bbb C}$ such that
(\ref{motion eq zw}) are satisfied everywhere. This function also
defines a triangular lattice with $MR=-1$.
\item Given a pair of complex--valued functions $(u,v)$ defined on
$V(\cT\cL)$ and satisfying the equation (\ref{motion eq zw})
everywhere, there exists a unique, up to an affine transformation,
function $w:V(\cT\cL)\mapsto{\Bbb C}$ such that the pairs $(v,w)$
and $(w,u)$ satisfy the same equation. The function $w$ also
defines a triangular lattice with $MR=-1$.
\end{enumerate}
\end{theorem}
{\bf Proof.} 1. To prove the first statement, we proceed as
follows. Let $\zz'\in V(\cT\cL)$, and let the vertices of an
elementary hexagonal with the center in $\zz'$ be enumerated as
$\zz_k=\zz'+\varepsilon^k$, $ k=1,2,\ldots,6$. Then the following
elementary triangles are positively oriented:
$(\zz_{2k},\zz_{2k-1}, \zz')$ and $(\zz_{2k},\zz_{2k+1},\zz')$ for
$k=1,2,3$ (with the agreement that $\zz_7=\zz_1$). According to
(\ref{motion eq zw}), we have:
\[
\frac{u_{2k-1}-u_{2k}}{u'-u_{2k-1}}=\frac{v'-v_{2k-1}}{v_{2k}-v'},\qquad
\frac{u_{2k+1}-u_{2k}}{u'-u_{2k+1}}=\frac{v'-v_{2k+1}}{v_{2k}-v'},\qquad
k=1,2,3.
\]
Dividing the first equation by the second one and taking the
product over $k=1,2,3$, we find:
\[
\prod_{k=1}^3\frac{u_{2k-1}-u_{2k}}{u_{2k+1}-u_{2k}}=1,
\]
which is nothing but (\ref{hex eq z}). The proof of (\ref{hex eq
w}) is similar.

2. As for the second statement, suppose we are given a function
$u$ on the whole of $V(\cT\cL)$. For an arbitrary elementary
triangle, if the values of $v$ in two vertices are known, the
equation (\ref{motion eq}) allows us to calculate the value of $v$
in the third vertex. Therefore, choosing arbitrarily the values of
$v$ in two neighboring vertices, we can extend this function on
the whole of $V(\cT\cL)$, provided this procedure is consistent.
It is easy to understand that it is enough to verify the
consistency in running once around a vertex. But this is assured
exactly by the equation (\ref{hex eq z}).

3. To prove the third statement, notice that the proof of Theorem
\ref{Equations of motion} shows that the formula
\begin{equation}
h(\ee)=w(\zz_2)-w(\zz_1)=\frac{1}{f(\ee)g(\ee)}=
\frac{1}{(u(\zz_2)-u(\zz_1))(v(\zz_2)-v(\zz_1))},
\end{equation}
valid for every edge $\ee=(\zz_1,\zz_2)$ of $\cT\cL$, correctly
defines the third field $h$ of the $fgh$--system. All affine
transformations of the field $w$ thus obtained, and only they,
lead to pairs $(v,w)$ and $(w,u)$ satisfying (\ref{motion eq zw}).
\qed \vspace{3mm}

{\bf Remark.} Notice that the above results remain valid in the
more general context, when the fields $f,g,h$ do not commute
anymore, e.g. when they take values in ${\Bbb H}$, the field of
quaternions. The formulation and the proof of Theorem
\ref{Equations of motion} hold in this case literally, while the
formula (\ref{hex eq z}) reads then as
\begin{equation}\label{hex eq z quat}
(u_1-u_2)(u_2-u_3)^{-1}(u_3-u_4)(u_4-u_5)^{-1}(u_5-u_6)(u_6-u_1)^{-1}=-1,
\end{equation}
and similarly for $v,w$.

\subsection{Circularity}

Recall that hexagonal circle patterns with $MR=-1$ lead to a
subclass of triangular lattices with $MR=-1$, namely those where
the points of one of the three hexagonal sublattices lie on
circles. We now prove a remarkable statement, assuring that this
subclass is stable with respect to the transformation $u\mapsto v$
described in Theorem \ref{z to w}.
\begin{theorem}\label{from circ to circ}
Let $u:V(\cH\cL_j)\mapsto{\Bbb C}$ define a hexagonal circle
pattern with $MR=-1$. Extend it with the centers of the circles to
$u:V(\cT\cL)\mapsto{\Bbb C}$, a triangular lattice with $MR=-1$.
Let $v:V(\cT\cL)\mapsto{\Bbb C}$ be the triangular lattice with
$MR=-1$ related to $u$ via (\ref{motion eq zw}). Then the
restriction of the map $v$ to the sublattice $\cH\cL_{j+1}$ also
defines a hexagonal circle pattern with $MR=-1$, while the points
$v$ corresponding to $\cT\cL\setminus\cH\cL_{j+1}$ are the centers
of the corresponding circles.
\end{theorem}
{\bf Proof} starts as the proof of Theorem \ref{z to w}. Let
$\zz'$ be a center of an arbitrary elementary hexagon of the
sublattice $\cH\cL_{j+1}$, i.e. $\zz'=k+\ell\omega+m\omega^2$ with
$k+\ell+m\equiv j+1\!\!\pmod 3$. Denote by
$\zz_k=\zz'+\varepsilon^k$, $k=1,2,\ldots,6$ the vertices of the
hexagon. As before, considering the positively oriented triangles
$(\zz_{2k}, \zz_{2k-1},\zz')$ and $(\zz_{2k},\zz_{2k+1},\zz')$,
$k=1,2,3$, surrounding the point $\zz'$, we come to the relations
\begin{equation}\label{circularity aux1}
\frac{u_{2k-1}-u_{2k}}{u'-u_{2k-1}}=\frac{v'-v_{2k-1}}{v_{2k}-v'}\,,\qquad
\frac{u_{2k+1}-u_{2k}}{u'-u_{2k+1}}=\frac{v'-v_{2k+1}}{v_{2k}-v'}\,,\qquad
k=1,2,3.
\end{equation}
But, obviously, $\zz_{2k-1}$ $(k=1,2,3)$ are centers of elementary
hexagons of the sublattice $\cH\cL_j$. By condition, the points
$u_{2k-2}$, $u_{2k}$ and $u'$ lie on a circle with the center in
$u_{2k-1}$. Therefore,
\begin{equation}\label{circularity aux2}
|u_{2k}-u_{2k-1}|=|u_{2k-2}-u_{2k-1}|=|u'-u_{2k-1}|,\qquad
k=1,2,3.
\end{equation}
So, the absolute values of the left--hand sides of all equations
in (\ref {circularity aux1}) are equal to 1. It follows that all
six points $v_1,v_2,\ldots,v_6$ lie on a circle with the center in
$v'$. \qed

\setcounter{equation}{0}
\section{Isomonodromic solutions}
\label{Sect isomonodromic}

Recall that we use triples $(k,\ell,m)\in{\Bbb Z}^3$ as
coordinates of the vertices $\zz=k+\ell\omega+m\omega^2$, and that
two such triples are identified iff they differ by the vector
$(n,n,n)$ with $n\in{\Bbb Z}$. By the $k$--axis we call the
straight line ${\Bbb R}\subset {\Bbb C}$, resp. by the
$\ell$--axis the straight line ${\Bbb R}\omega$, and by the
$m$--axis the straight line ${\Bbb R}\omega^2$.

It will be sometimes convenient to use the symbols
$\tilde{\cdot}$, $\hat{\cdot}$ and $\bar{\cdot}$ to denote the
shifts of various objects in the positive direction of the axes
$k$, $\ell$, $m$, respectively, and the symbols
$\undertilde{\cdot}$, $\underhat{\cdot}$, $\underline{\cdot}$ to
denote the shifts in the negative directions. This will apply to
vertices, edges and elementary triangles of $\cT\cL$, as well as
to various objects assigned to them. For example, if $\zz\in
V(\cT\cL)$, then
\[
\widetilde{\zz}=\zz+1, \quad \undertilde{\zz}=\zz-1, \quad
\widehat{\zz}=\zz+\omega, \quad \underhat{\zz}=\zz-\omega, \quad
\bar{\zz}=\zz+\omega^2,\quad \underline{\zz}=\zz-\omega^2.
\]
Similarly, if $\ee=(\zz_1,\zz_2)\in E(\cT\cL)$, then
\[
\widetilde{\ee}=(\zz_1+1,\zz_2+1), \quad
\widehat{\ee}=(\zz_1+\omega,\zz_2+\omega), \quad
\bar{\ee}=(\zz_1+\omega^2,\zz_2+\omega^2),\quad {\rm etc.}
\]

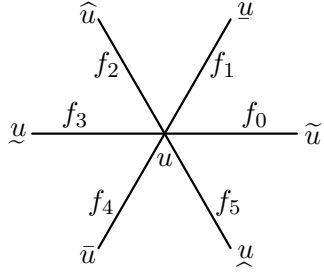
\begin{figure}[htbp]
  \begin{center}
    \setlength{\unitlength}{20pt}
\begin{picture}(6,5)
\thicklines \path(3,2.5)(5.5,2.5) \path(3,2.5)(0.5,2.5)
\path(3,2.5)(4.25,4.665) \path(3,2.5)(4.25,0.335)
\path(3,2.5)(1.75,4.665) \path(3,2.5)(1.75,0.335)
\put(5.8,2.5){\makebox(0,0){$\widetilde{u}$}}
\put(0.2,2.5){\makebox(0,0){$\undertilde{u}$}}
\put(4.5,0.2){\makebox(0,0){$\underhat{u}$}}
\put(1.55,4.8){\makebox(0,0){$\widehat{u}$}}
\put(4.5,4.8){\makebox(0,0){$\underbar{u}$}}
\put(1.55,0.2){\makebox(0,0){$\bar{u}$}}
\put(4.7,2.85){\makebox(0,0){$f_0$}}
\put(1.3,2.85){\makebox(0,0){$f_3$}}
\put(1.9,3.8){\makebox(0,0){$f_2$}}
\put(4.1,3.8){\makebox(0,0){$f_1$}}
\put(1.8,1.2){\makebox(0,0){$f_4$}}
\put(4.2,1.2){\makebox(0,0){$f_5$}} \put(3,2){\makebox(0,0){$u$}}
\end{picture}
    \caption{Notations for $u$ and $f$.}
    \label{fig:notationSec5}
  \end{center}
\end{figure}

A fundamental role in the subsequent presentation will be played
by a {\it non-autonomous constraint} for the solutions of the
$fgh$--system. This constraint consists of a pair of equations
which are formulated for every vertex $\zz\in V(\cT\cL)$ and
include the values of the fields on the edges incident to $\zz$,
i.e. on the {\it star} of this vertex. It will be convenient to
fix a numeration of these edges as follows:
\begin{eqnarray}
\ee_0=(\zz,\widetilde{\zz}),\quad \ee_2=(\zz,\widehat{\zz}),\quad
\ee_4=(\zz,\bar{\zz}), \\
\ee_1=(\underbar{\zz},\zz),\quad
\ee_3=(\undertilde{\zz},\zz),\quad \ee_5=(\underhat{\zz},\zz).
\end{eqnarray}
The notations $f_0,\ldots,f_6$ will refer to the values of the
field $f$ on these edges:
\begin{eqnarray}
f_0=\widetilde{u}-u,\quad  f_2=\widehat{u}-u,\quad  f_4=\bar{u}-u, \\
f_1=u-\underbar{u}, \quad f_3=u-\undertilde{u},\quad
f_5=u-\underhat{u},
\end{eqnarray}
and similarly for the fields $g$, $h$, see Fig.
\ref{fig:notationSec5}.

The constraint looks as follows:
\begin{equation}\label{constr 1}
\alpha u=k\frac{f_0g_0f_3}{f_0g_0+g_0f_3+f_3g_3}+
\ell\frac{f_2g_2f_5}{f_2g_2+g_2f_5+f_5g_5}+
m\frac{f_4g_4f_1}{f_4g_4+g_4f_1+f_1g_1},
\end{equation}
\begin{equation}\label{constr 2}
\beta v=k\frac{g_0f_3g_3}{f_0g_0+g_0f_3+f_3g_3}+
\ell\frac{g_2f_5g_5}{f_2g_2+g_2f_5+f_5g_5}+
m\frac{g_4f_1g_1}{f_4g_4+g_4f_1+f_1g_1}.
\end{equation}
These are supposed to be the equations for the vertex
$\zz=k+\ell\omega+m\omega^2$, and we use the notations $u=u(\zz)$,
$v=v(\zz)$. Since the fields $u$, $v$ are defined only up to an
affine transformation, one should replace the left--hand sides of
the above equations by $\alpha u+\phi$, $\beta v+\psi$,
respectively, with arbitrary constants $\phi$, $\psi$. In the form
we have choosen it is imposed that the fields $u$, $v$ are
normalized to vanish in the origin.

\begin{proposition}\label{constraint preliminary}
The equations (\ref{constr 1}), (\ref{constr 2}) are well defined
equations for the point $\zz\in V(\cT\cL)$, i.e. they are
invariant under the shift $(k,\ell,m)\mapsto(k+n,\ell+n,m+n)$,
provided the equations (\ref{motion eq zw}) hold.
\end{proposition}
{\bf Proof} is technical and is given in the Appendix
\ref{Appendix}. \qed \vspace{2mm}

We mention an important consequence of this proposition.
Apparently, the constraint (\ref{constr 1}), (\ref{constr 2})
relates the values of the fields $u$, $v$ in {\it seven} points
shown on Fig.~\ref{fig:notationSec5}. However, we are free to
choose any representative $(k,\ell,m)$ for $\zz$. In particular,
we can let vanish any one of the coordinates $k$, $\ell$, $m$. In
the corresponding representation the constraint relates the values
of the fields $u$, $v$ in {\it five} points, belonging to any one
of the three possible four--leg crosses through $\zz$.
\vspace{2mm}

An essential algebraic property of the constraint (\ref{constr
1}), (\ref{constr 2}) is given by the following statement.
\begin{proposition}\label{third constraint}
If the equations (\ref{motion eq zw}) hold, then the constraints
(\ref{constr 1}), (\ref{constr 2}) imply a similar equation for
the field $w$ (vanishing at $\zz=0$):
\begin{equation}\label{constr 3}
\gamma w=k\frac{1}{f_0g_0+g_0f_3+f_3g_3}+
\ell\frac{1}{f_2g_2+g_2f_5+f_5g_5}+
m\frac{1}{f_4g_4+g_4f_1+f_1g_1},
\end{equation}
where $\gamma=1-\alpha-\beta$.
\end{proposition}
{\bf Proof} is again based on calculations and is relegated to the
Appendix \ref{Appendix}. \qed \vspace{2mm}

{\bf Remark.} We notice that restoring the fields $h_k=1/(f_kg_k)$
allows us to rewrite the equations (\ref{constr 2}), (\ref{constr
3}) as
\begin{eqnarray}
\beta v & = & k\frac{g_0h_0g_3}{g_0h_0+h_0g_3+g_3h_3}+
\ell\frac{g_2h_2g_5}{g_2h_2+h_2g_5+g_5h_5}+
m\frac{g_4h_4g_1}{g_4h_4+h_4g_1+g_1h_1}, \label{constr 2 alt}\\
\gamma w & = & k\frac{h_0f_0h_3}{h_0f_0+f_0h_3+h_3f_3}+
\ell\frac{h_2f_2h_5}{h_2f_2+f_2h_5+h_5f_5}+
m\frac{h_4f_4h_1}{h_4f_4+f_4h_1+h_1f_1}, \label{constr 3 alt}
\end{eqnarray}
which coincides with (\ref{constr 1}) via a cyclic permutation of
fields $(f,g,h)\mapsto(g,h,f)$ performed once or twice,
respectively, and accompanied by changing $\alpha$ to $\beta$,
$\gamma$, respectively. \vspace{1mm}

Another similar remark: as it follows from the formulas
(\ref{constr 1 well aux2}), (\ref{constr 1 well aux3}) used in the
proof of Proposition \ref{constraint preliminary} (and their
analogs for the fields $g$, $h$), the constraints (\ref{constr
1}), (\ref{constr 2}), (\ref{constr 3}) may be rewritten as
equations for the single field $u$, resp. $v$, $w$:
\begin{eqnarray}
\alpha u & = & k\frac{f_0f_3(f_1+f_2)}{(f_0-f_2)(f_1-f_3)}+
\ell\frac{f_2f_5(f_3+f_4)}{(f_2-f_4)(f_3-f_5)}+
m\frac{f_4f_1(f_5+f_0)}{(f_4-f_0)(f_5-f_1)}, \label{constr 1 z}\\
\beta v & = & k\frac{g_0g_3(g_1+g_2)}{(g_0-g_2)(g_1-g_3)}+
\ell\frac{g_2g_5(g_3+g_4)}{(g_2-g_4)(g_3-g_5)}+
m\frac{g_4g_1(g_5+g_0)}{(g_4-g_0)(g_5-g_1)}, \label{constr 2 w}\\
\gamma w & = & k\frac{h_0h_3(h_1+h_2)}{(h_0-h_2)(h_1-h_3)}+
\ell\frac{h_2h_5(h_3+h_4)}{(h_2-h_4)(h_3-h_5)}+
m\frac{h_4h_1(h_5+h_0)}{(h_4-h_0)(h_5-h_1)}. \label{constr 3 v}
\end{eqnarray}
However, in this form, unlike the previous one, the terms attached
to the variable $k$ (say), contain not only the fields on two edges 
$\ee_0$, $\ee_3$ parallel to the $k$--axis. This form is therefore less
suited for the solution of the Cauchy problem for the constrained
$fgh$--system, which we discuss now.

\begin{theorem}\label{compatibility}
For arbitrary $\alpha,\beta\in{\Bbb C}$ the constraint
(\ref{constr 1}), (\ref{constr 2}) is compatible with the
equations (\ref{motion eq zw}).
\end{theorem}
{\bf Proof.} To prove this statement, one has to demonstrate the
solvability of a reasonably posed Cauchy problem for the
$fgh$--system constrained by (\ref{constr 1}), (\ref{constr 2}).
In this context, it is unnatural to assume that the fields $u$,
$v$ vanish at the origin, so that we replace (only in this proof)
the left--hand sides of (\ref{constr 1}), (\ref{constr 2}) by
$\alpha u+\phi$, $\beta v+\psi$, with arbitrary $\phi,\psi\in{\Bbb
C}$. We show that reasonable Cauchy data are given by the values
of two fields $u$, $v$ (say) in three points $\zz_0$,
$\zz_1=\zz_0+1$, $\zz_2=\zz_0+\omega$, where $\zz_0$ is arbitrary.
According to Lemma \ref{lemma 4th point}, these data yield via the
equations of the $fgh$--system the values of $u$, $v$ in
$\zz_3=\zz_0+1+\omega$. Further, these data together with the
constraint (\ref{constr 1}), (\ref{constr 2}) determine uniquely
the values of $u$, $v$ in $\zz_4=\zz_0+\omega^2$. Indeed, assign
$u(\zz_4)=\xi$, $v(\zz_4)=\eta$, where $\xi$, $\eta$ are two
arbitrary complex numbers. The constraint uniquely defines the
values of $u$, $v$ in the point $\zz_5=\zz_0-\omega$. The
requirement that these values agree with the ones obtained via
Lemma \ref{lemma 4th point} from the points $\zz_0$, $\zz_1$,
$\zz_4$, gives us two equations for $\xi$, $\eta$. It is shown by
a direct computation that these equations have a unique solution,
which is expressed via rational functions of the data at $\zz_0$,
$\zz_1$, $\zz_2$. It is also shown that the same solution is
obtained, if we work with $\zz_6=\zz_0-1$ instead of $\zz_5$.
Having found the fields $u$, $v$ at $\zz_4$, we determine
simultaneously $u$, $v$ at $\zz_5$, $\zz_6$. Now a similar
procedure allows us to determine $u$, $v$ at $\zz_7=\zz_0+2$ and
$\zz_8=\zz_0+2\omega$, using the constraint at the points $\zz_1$
and $\zz_2$, respectively. Simultaneously the values of $u$, $v$
are found at $\zz_9=\zz_0+2+\omega$ and $\zz_{10}=
\zz_0+1+2\omega$. A continuation of this procedure delivers the
values of $u$, $v$ on the both semiaxes
\[
\Big\{\zz=k: k\ge 0\Big\}\cup\Big\{\zz=\ell\omega: \ell\ge
0\Big\},
\]
using the condition that the constraint (\ref{constr 1}),
(\ref{constr 2}) is fulfilled on these semiaxes. As we know from
Proposition \ref{Cauchy data for fgh system}, these data are
enough to determine the solution of the $fgh$--system in the whole
sector
\[
\Big\{\zz=k+\ell\omega: k,\ell\ge 0\Big\}= \Big\{\zz\in V(\cT\cL):
0\le{\rm\arg}(\zz)\le2\pi/3\Big\}.
\]
It remains to prove that this solution fulfills also the
constraint (\ref{constr 1}), (\ref{constr 2}) in the whole sector.
This follows by induction from the following statement:
\begin{lemma}\label{lemma for compatibility}
If the constraint (\ref{constr 1}), (\ref{constr 2}) is satisfied
in $\zz_0$, $\zz_1$, $\zz_2$, then it is satisfied also in
$\zz_3$.
\end{lemma}
The constraint at $\zz_3$ includes the data at five points
$\zz_1$, $\zz_2$, $\zz_3$, $\zz_9$, $\zz_{10}$. As we have seen,
the data at $\zz_3$, $\zz_9$, $\zz_{10}$ are certain (complicated)
functions of the data at $\zz_0$, $\zz_1$, $\zz_2$. Therefore, to
check the constraint at $\zz_3$, one has to check that two
(complicated) equations for the values of $u$, $v$ at $\zz_0$,
$\zz_1$, $\zz_2$ are satisfied identically. This has been done
with the help of the Mathematica computer algebra system. \qed

\begin{figure}[htbp]
  \begin{center}
    \includegraphics[width = 0.5\hsize]{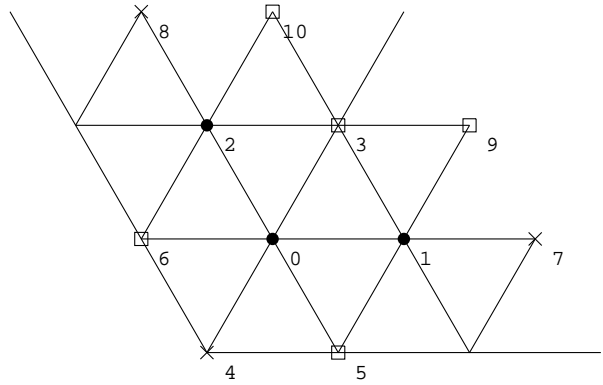}
    \caption{To the proof of Theorem~\ref{compatibility}: labelling of the
    points}
    \label{fig:compatibility}
  \end{center}
\end{figure}
Now we show how the constraint (\ref{constr 1}), (\ref{constr 2})
appears in the context of isomonodromic solutions of integrable
systems. In this context, the results look better with a different
gauge of the transition matrices for the $fgh$--system. Namely, we
conjugate them with the matrix ${\rm diag}(1,\lambda,\lambda^2)$,
and then multiply by $(1+\lambda^3)^{1/3}$ in order to get rid of
the normalization of the determinant. Writing then $\mu$ for
$\lambda^3$, we end up with the matrices
\begin{equation}
\cL(\mu)=\left(\begin{array}{ccc} 1 & f & 0 \\ 0 & 1 & g \\ \mu h
& 0 & 1\end{array}\right), \quad fgh=1.
\end{equation}
The zero curvature condition turns into
\begin{equation}\label{zero curv in mu}
\cL(\ee_3,\mu)\cL(\ee_2,\mu)\cL(\ee_1,\mu)=(1+\mu)I,
\end{equation}
$\ee_1$, $\ee_2$, $\ee_3$ being the consecutive positively
oriented edges of an elementary triangle of $\cT\cL$. This implies
some slight modifications also for the notion of the wave
function. Namely, the previous formula does not allow to define
the function $\Psi$ on $V(\cT\cL)$ such that
\[
\Psi(\zz_2,\mu)=\cL(\ee,\mu)\Psi(\zz_1,\mu)
\]
holds, whenever $\ee=(\zz_1,\zz_2)$. The way around this
difficulty is the following. We define the wave function $\Psi$ on
a covering of $V(\cT\cL)$. Namely, over each point
$\zz=k+\ell\omega+m\omega^2$ now sits a sequence
\begin{equation}\label{cover cond for psi}
\Psi_{k+n,\ell+n,m+n}(\mu)=(1+\mu)^n\Psi_{k,\ell,m}(\mu),\quad
n\in{\Bbb Z}.
\end{equation}
The values of these functions in neighboring vertices are related
by natural formulas
\begin{equation}\label{wave evolution in mu}
\left\{\begin{array}{l}
\Psi_{k+1,\ell,m}(\mu)=\cL(\ee_0,\mu)\Psi_{k,\ell,m}(\mu),\quad
\ee_0=(\zz,\zz+1),\\
\Psi_{k,\ell+1,m}(\mu)=\cL(\ee_2,\mu)\Psi_{k,\ell,m}(\mu),\quad
\ee_2=(\zz,\zz+\omega),\\
\Psi_{k,\ell,m+1}(\mu)=\cL(\ee_4,\mu)\Psi_{k,\ell,m}(\mu), \quad
\ee_4=(\zz,\zz+\omega^2).
\end{array}\right.
\end{equation}
We call a solution $(u,v):V(\cT\cL)\mapsto{\Bbb C}^2$ of the
equations (\ref{motion eq zw}) {\itbf isomonodromic} (cf.
\cite{I}), if there exists the wave function $\Psi:{\Bbb
Z}^3\mapsto {\rm GL}(3,{\Bbb C})[\mu]$ satisfying (\ref{wave
evolution in mu}) and some linear differential equation in $\mu$:
\begin{equation}\label{eq in mu}
\frac{d}{d\mu}\Psi_{k,\ell,m}(\mu)=\cA_{k,\ell,m}(\mu)\Psi_{k,\ell,m}(\mu),
\end{equation}
where $\cA_{k,\ell,m}(\mu)$ are $3\times 3$ matrices, meromorphic
in $\mu$, with the poles whose position and order do not depend on
$k,\ell,m$.

Obviously, due to (\ref{cover cond for psi}), the matrix $\cA$ has
to fulfill the condition
\begin{equation}\label{cover cond for A}
\cA_{k+n,\ell+n,m+n}(\mu)=\cA_{k,\ell,m}(\mu)+\frac{n}{1+\mu}I,\quad
n\in{\Bbb Z}.
\end{equation}
\begin{theorem}\label{monodromy}
 Solutions of the equations (\ref{motion eq zw}) satisfying
the constraints (\ref{constr 1}), (\ref{constr 2}) are
isomonodromic. The corresponding matrix $\cA_{k,\ell,m}$ is given
by the following formulas:
\begin{equation}\label{ans A}
\cA_{k,\ell,m}=\frac{C_{k,\ell,m}}{1+\mu}+\frac{D(\zz)}{\mu},
\end{equation}
where $C_{k,\ell,m}$ and $D(\zz)$ are $\mu$--independent matrices:
\begin{equation}\label{C}
C_{k,\ell,m}=kP_0(\zz)+\ell P_2(\zz)+mP_4(\zz),
\end{equation}
$P_{0,2,4}$ are rank 1 matrices
\begin{equation}\label{P}
P_j(\zz)=\frac{1}{f_jg_j+g_jf_{j+3}+f_{j+3}g_{j+3}}
\left(\begin{array}{ccc} f_jg_j & -f_jg_jf_{j+3} & f_jg_jf_{j+3}g_{j+3} \\
-g_j & g_jf_{j+3} & -g_jf_{j+3}g_{j+3} \\
1 & -f_{j+3} & f_{j+3}g_{j+3}   \end{array}\right),\quad j=0,2,4,
\end{equation}
and the matrix $D$ is well defined on $V(\cT\cL)$ and not only on
its covering ${\Bbb Z}^3$:
\begin{equation}\label{D}
D(\zz)=\left(\begin{array}{ccc}
-(2\alpha+\beta)/3 & \alpha u & \beta a-\alpha a' \\
0 & (\alpha-\beta)/3 & \beta v \\
0 & 0 & (2\beta+\alpha)/3 \end{array}\right),
\end{equation}
where the functions $a,a':V(\cT\cL)\mapsto{\Bbb C}$ are solutions
of the equations (\ref{eq for a}), (\ref{eq for a'}).
\end{theorem}
{\bf Proof} can be found in the Appendix \ref{Appendix}. \qed

\setcounter{equation}{0}
\section{Isomonodromic solutions and circle patterns}
\label{Sect isomonodromic patterns}

We now consider isomonodromic solutions of the $fgh$--system
satisfying the constraint (\ref{constr 1}), (\ref{constr 2}),
which are special in two respects:
\begin{itemize}
\item First, the constants $\alpha$ and $\beta$ in the constraint equations
are not arbitrary, but are {\it equal}: $\alpha=\beta$, so that
$\gamma=1-2\alpha$.
\item Second, the initial conditions will be choosen in a special way.
\end{itemize}
We will show that the resulting solutions lead to hexagonal circle
patterns.

First of all, we discuss the Cauchy data which allow one to
determine a solution of the $fgh$--system augmented by the
constraints (\ref{constr 1}), (\ref{constr 2}). Of course, the
fields $u$, $v$, $w$ have to vanish in the origin $\zz=0$. Next,
one sees easily that, given $u$ and $v$ in one of the points
neighboring to $0$, the constraint allows to calculate one after
another the values of $u$ and $v$ in all points of the
corresponding axis. For instance, fixing some values of $u(1)$ and
$v(1)$, we can calculate all $u(k)$ and $v(k)$ from the relations
\begin{equation}\label{contr u on axis}
\alpha
u(k)=k\frac{f(k)g(k)f(k-1)}{f(k)g(k)+g(k)f(k-1)+f(k-1)g(k-1)},
\end{equation}
\begin{equation}\label{contr v on axis}
\beta
v(k)=k\frac{g(k)f(k-1)g(k-1)}{f(k)g(k)+g(k)f(k-1)+f(k-1)g(k-1)},
\end{equation}
where we have set
\begin{equation}\label{fg uv}
f(k)=u(k+1)-u(k),\qquad g(k)=v(k+1)-v(k).
\end{equation}
Indeed, we start with $u(0)=0$, $v(0)=0$, $f(0)=u(1)$,
$g(0)=v(1)$, and continue via the recurrent formulas, which are
easily seen to be equivalent to (\ref{contr u on axis}),
(\ref{contr v on axis}), (\ref{fg uv}):
\begin{equation}\label{recur zw}
u(k)=u(k-1)+f(k-1),\qquad v(k)=v(k-1)+g(k-1),
\end{equation}
\begin{eqnarray}
f(k) & = & \frac{\alpha u(k)}{\beta v(k)}\,g(k-1), \label{recur f} \\
g(k) & = & \frac{\beta v(k)} {k-\displaystyle\frac{\alpha
u(k)}{f(k-1)}-\displaystyle\frac{\beta v(k)} {g(k-1)}}.
\label{recur g}
\end{eqnarray}
So, given the values of the fields $u$ and $v$ (and hence of $w$)
in the points $\zz=1$ and $\zz=\omega$, we get their values in all
points $\zz=k$ and $\zz=\ell\omega$ of the positive $k$- and
$\ell$-semiaxes. It is easy to see that $u(k)/u(1)$ and
$v(k)/v(1)$ do not depend on $u(1)$ and $v(1)$, respectively, so
that all points $u(k)$ lie on a straight line, and so do all
points $v(k)$. Similar statements hold also for all points
$u(\ell\omega)$ and for all points $v(\ell\omega)$. And, of
course, the third field $w$ behaves analogously.

So, we get the values of $u$ and $v$ in all points on the border
of the sector
\begin{equation}\label{sector}
S=\Big\{\zz\in V(\cT\cL): 0\le{\rm\arg}(\zz)\le 2\pi/3\Big\}
=\Big\{\zz=k+\ell\omega: k,\ell\ge 0\Big\}.
\end{equation}
Proposition \ref{Cauchy data for fgh system} assures that these
data determine the values of $u$ and $v$ in all points of $S$. By
Theorem \ref{compatibility} (more precisely, by Lemma \ref{lemma
for compatibility}) the solution thus obtained will satisfy the
constraint (\ref{constr 1}), (\ref{constr 2}) in the whole sector
$S$.

Now we are in a position to specify the above mentioned
isomonodromic solutions.
\begin{theorem}\label{circular}
Let $\beta=\alpha$. Let $u,v,w:S\mapsto{\Bbb C}$ be the solutions
of the $fgh$--system with the constraint (\ref{constr 1}),
(\ref{constr 2}), with the initial conditions
\begin{equation}\label{ini uv}
u(1)=v(1)=1, \quad u(\omega)=v(\omega)=\exp(i\theta),
\end{equation}
where $0<\theta<\pi$. Then all three maps $u,v,w$ define hexagonal
circle patterns with $MR=-1$ in the sector $S$. More precisely, if
$\zz_k=\zz'+\varepsilon^k$, $k=1,2,\ldots,6,$ are the vertices of
an elementary hexagon in this sector, then:
\begin{itemize}
\item $u(\zz_1),u(\zz_2),\ldots,u(\zz_6)$ lie on a circle with the center
in $u(\zz')$ whenever $\zz'\in S\setminus V(\cH\cL_1)$,
\item $v(\zz_1),v(\zz_2),\ldots,v(\zz_6)$ lie on a circle with the center
in $v(\zz')$ whenever $\zz'\in S\setminus V(\cH\cL_2)$,
\item $w(\zz_1),w(\zz_2),\ldots,w(\zz_6)$ lie on a circle with the center
in $w(\zz')$ whenever $\zz'\in S\setminus V(\cH\cL_0)$.
\end{itemize}
\end{theorem}

\begin{figure}[htbp]
  \begin{center}
         \setlength{\unitlength}{28.45pt}
\begin{picture}(15,5)
\thicklines \path(0,0.5)(3.9,0.5) \path(0,0.5)(1.2,4.2)
\path(0,0.5)(3.1,2.76) \path(1.2,4.2)(3.1,2.76)
\path(3.1,2.76)(3.9,0.5) \put(0.1,0.2){\makebox(0,0){$u_0$}}
\put(3.9,0.2){\makebox(0,0){$u_1$}}
\put(3.3,3){\makebox(0,0){$u_3$}}
\put(1.2,4.5){\makebox(0,0){$u_2$}} \linethickness{0.7pt}
\qbezier(0.5,0.5)(0.45,0.85)(0.15,0.923)
\put(0.6,0.6){\makebox(1,0.2){$2\pi\!\!-\!\!2\vartheta$}}
\path(7,0.5)(4.4,2.4) \path(7,0.5)(10.2,0.5) \path(7,0.5)(7.6,2.4)
\path(10.2,0.5)(7.6,2.4) \path(4.4,2.4)(7.6,2.4)
\put(7,0.2){\makebox(0,0){$v_0$}}
\put(10.2,0.2){\makebox(0,0){$v_1$}}
\put(7.6,2.75){\makebox(0,0){$v_3$}}
\put(4.4,2.75){\makebox(0,0){$v_2$}} \linethickness{0.7pt}
\qbezier(7.4,0.5)(7.1,0.8)(6.677,0.736)
\put(7.25,0.7){\makebox(0.5,0.2){$\vartheta$}}
\path(11.7,0.5)(13.8,0.5) \path(11.7,0.5)(12.75,3.8)
\path(11.7,0.5)(10,1.8) \path(10,1.8)(12.75,3.8)
\path(12.75,3.8)(13.8,0.5) \put(11.7,0.2){\makebox(0,0){$w_0$}}
\put(13.8,0.2){\makebox(0,0){$w_1$}}
\put(12.75,4.1){\makebox(0,0){$w_3$}}
\put(9.8,2.1){\makebox(0,0){$w_2$}}
\end{picture}
    \caption{To Lemma \ref{geometry}: elementary triangles for
    $u$, $v$, and $w$ are isosceles.}
    \label{fig:statement17}
  \end{center}
\end{figure}
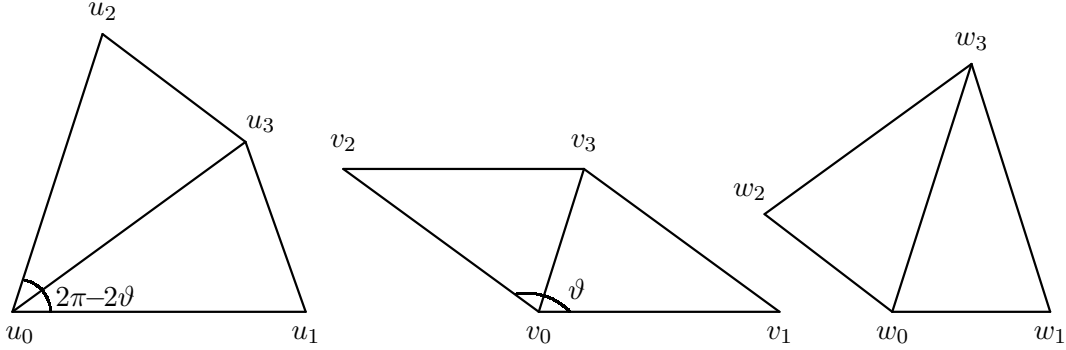

{\bf Proof} follows from the above inductive construction with the
help of two lemmas. The first one shows that if $\beta=\alpha$
then the constraint yields a very special property of the
sequences of the values of the fields $u$, $v$, $w$ in the points
of the $k$- and $\ell$-axes.
\begin{lemma}\label{equidist}
If $\beta=\alpha$, then for $k,\ell\ge 1$:
\begin{eqnarray}
|u(3k-1)-u(3k-2)| & = & |u(3k-2)-u(3k-3)|, \label{u=u on k}\\
|v(3k)-v(3k-1)| & = & |v(3k-1)-v(3k-2)|,\label{v=v on k}\\
|w(3k+1)-w(3k)| & = & |w(3k)-w(3k-1)|,\label{w=w on k}\\ \nonumber\\
|u((3\ell-1)\omega)-u((3\ell-2)\omega)| & = & |u((3\ell-2)\omega)-
u((3\ell-3)\omega)|,  \label{u=u on l}\\
|v(3\ell\omega)-v((3\ell-1)\omega)| & = &
|v((3\ell-1)\omega)-v((3\ell-2)\omega)|,  \label{v=v on l}\\
|w((3\ell+1)\omega)-w(3\ell\omega)| & = &
|w(3\ell\omega)-w((3\ell-1)\omega)|.  \label{w=w on l}
\end{eqnarray}
\end{lemma}
The second one allows to extend inductively these special
properties to the whole sector (\ref{sector}).
\begin{lemma}\label{geometry}
Consider two elementary triangles with the vertices $\zz_0$,
$\zz_1=\zz_0+1$, $\zz_2=\zz_0+\omega$, and $\zz_3=\zz_0+1+\omega$.
Suppose that
\begin{itemize}
\item[{\rm (i)}] $|u(\zz_1)-u(\zz_0)|=|u(\zz_2)-u(\zz_0)|$;
\item[{\rm (ii)}] $\measuredangle v(\zz_1)v(\zz_0)v(\zz_2)=\vartheta\;\;$ and
$\;\;\measuredangle u(\zz_1)u(\zz_0)u(\zz_2)=2\pi-2\vartheta\;\;$
for some $\vartheta$.
\end{itemize}
Then
\begin{equation}
|u(\zz_3)-u(\zz_0)|=|u(\zz_1)-u(\zz_0)|=|u(\zz_2)-u(\zz_0)|,
\end{equation}
and hence
\begin{equation}
|v(\zz_3)-v(\zz_1)|=|v(\zz_0)-v(\zz_1)|, \qquad
|v(\zz_3)-v(\zz_2)|=|v(\zz_0)-v(\zz_2)|
\end{equation}
and
\begin{equation}
|w(\zz_3)-w(\zz_1)|=|w(\zz_3)-w(\zz_2)|=|w(\zz_3)-w(\zz_0)|
\end{equation}
\end{lemma}
The assertion of this lemma is illustrated on Fig.
\ref{fig:statement17}.

First of all, we show how do these lemmas work towards the proof
of Theorem \ref{circular}. The initial conditions (\ref{ini uv})
imply:
\begin{equation}\label{ini w}
w(1)=1, \quad w(\omega)=\exp(-2i\theta)=\exp(i(2\pi-2\theta)).
\end{equation}
Therefore, the conditions of Lemma \ref{geometry} are fulfilled in
the point $\zz_0=0$ with the fields $(w,u,v)$ instead of
$(u,v,w)$. From this Lemma it follows that
\begin{itemize}
\item[$({\rm a}_0)$] The points $w(1)$, $w(\omega)$, $w(1+\omega)$ are
equidistant from $w(0)$;
\item[$({\rm b}_0)$] The points $v(0)$, $v(1)$, $v(\omega)$ are equidistant
from $v(1+\omega)$;
\item[$({\rm c}_0)$] The points $u(1+\omega)$, $u(0)$ are equidistant from
$u(\omega)$;
\item[$({\rm d}_0)$] The points $u(1+\omega)$, $u(0)$ are equidistant from
$u(1)$.
\end{itemize}
\begin{figure}[htbp]
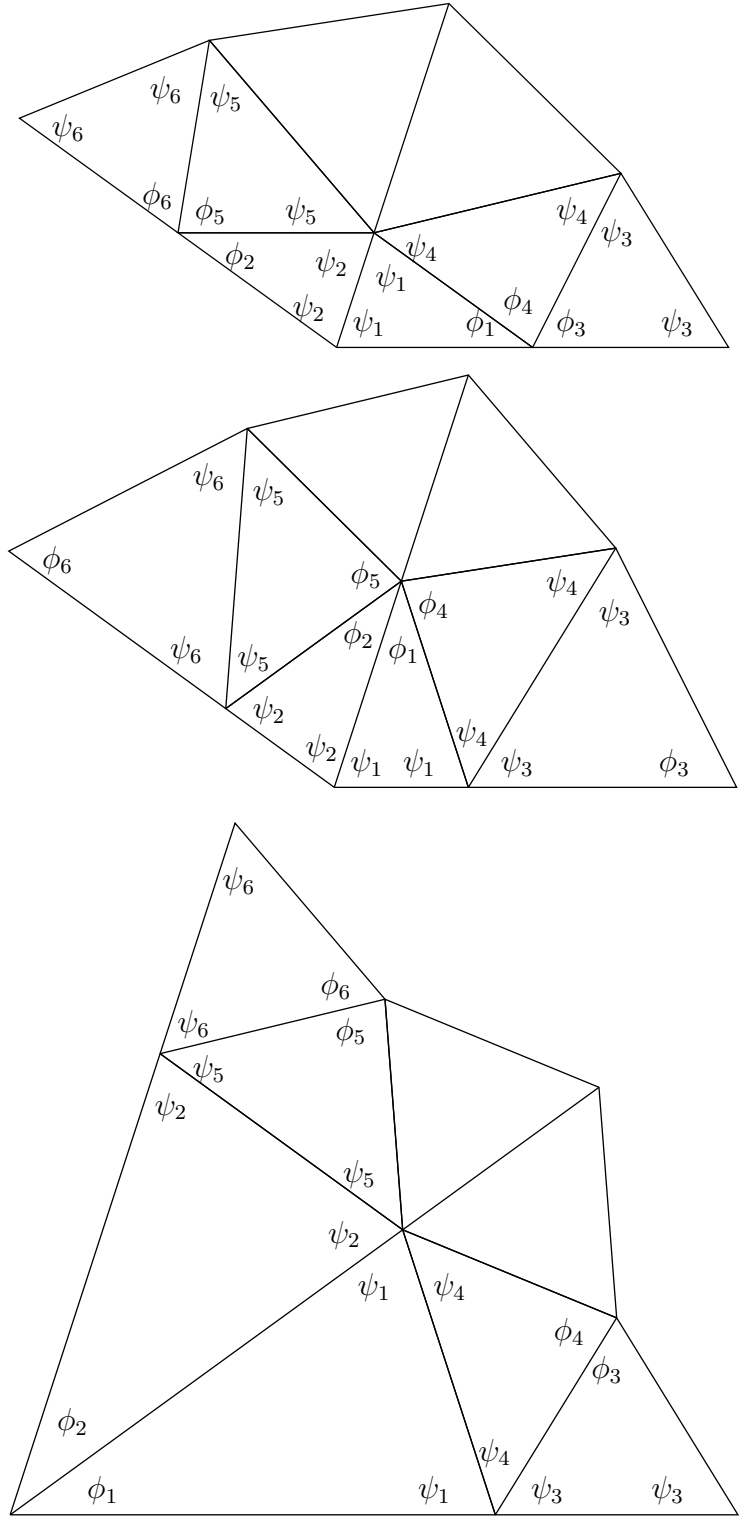

  \begin{center}
    \input a_proof15u.tex
    \input a_proof15v.tex
    \input a_proof15w.tex
    \caption{To the proof of Theorem \ref{circular}: similar isosceles
    triangles for $u$, $v$, and $w$}
    \label{fig:proof15}
  \end{center}
\end{figure}
Since, by Lemma \ref{equidist}, we have $|u(0)-u(1)|=|u(2)-u(1)|$,
there follows from $({\rm d}_0)$ that
$|u(1+\omega)-u(1)|=|u(2)-u(1)|$. Finally, from Lemma
\ref{geometry} there follows that (see Fig.~\ref{fig:proof15})
\[
\measuredangle v(2)v(1)v(1+\omega)=\pi-\psi_1,\qquad
\measuredangle
u(2)u(1)u(1+\omega)=\pi-\phi_1=2\psi_1=2\pi-2(\pi-\psi_1).
\]
Therefore the conditions of Lemma \ref{geometry} are fulfilled in
the point $\zz_0=1$ with the fields $(u,v,w)$. We deduce that
\begin{itemize}
\item[$({\rm a}_1)$] The points $u(2)$, $u(1+\omega)$, $u(2+\omega)$ are
equidistant from $u(1)$;
\item[$({\rm b}_1)$] The points $w(1)$, $w(2)$, $w(1+\omega)$ are equidistant
from $w(2+\omega)$;
\item[$({\rm c}_1)$] The points $v(2+\omega)$, $v(1)$ are equidistant
from $v(1+\omega)$, which adds the point $v(2+\omega)$ to the list
of equidistant neighbors of $v(1+\omega)$ from the conclusion
$({\rm b}_0)$ above; and
\item[$({\rm d}_1)$] The points $v(2+\omega)$, $v(1)$ are equidistant from
$v(2)$.
\end{itemize}
By Lemma \ref{equidist}, we have $|v(1)-v(2)|=|v(3)-v(2)|$, and
there follows from $({\rm d}_1)$ that
$|v(2+\omega)-v(2)|=|v(3)-v(2)|$. Finally, from Lemma
\ref{geometry} there follows that (see Fig.~\ref{fig:proof15})
\[
\measuredangle w(3)w(2)w(2+\omega)=\pi-\psi_3,\qquad
\measuredangle
v(3)v(2)v(2+\omega)=\pi-\phi_3=2\psi_3=2\pi-2(\pi-\psi_3).
\]
Hence, the conditions of Lemma \ref{geometry} are again fulfilled
in the point $\zz_0=2$ with the fields $(v,w,u)$.

These arguments may be continued by induction along the $k$-axis,
and, by symmetry, along the $\ell$-axis. This delivers all the
necessary relations which involve the points $\zz=k+\ell\omega$
with $k\le 1$ or $\ell\le 1$. We call them the relations of the
level 1.

The arguments of the level 2 start with the pair of fields $(v,w)$
at the point $\zz=1+\omega$. We have the level 1 relation
\[
|v(2+\omega)-v(1+\omega)|=|v(1+2\omega)-v(1+\omega)|.
\]
For the angles, we have from the level 1 (see Fig.
\ref{fig:proof15}):
\begin{eqnarray*}
\measuredangle
w(2+\omega)w(1+\omega)w(1+2\omega) & = & 2\pi-(\psi_1+\psi_2+\psi_4+\psi_5), \\
\measuredangle v(2+\omega)v(1+\omega)v(1+2\omega) & = &
2\pi-(\phi_1+\phi_2+\phi_4+\phi_5)=
2\pi-2(2\pi-\psi_1-\psi_2-\psi_4-\psi_5).
\end{eqnarray*}
So, the conditions of Lemma \ref{geometry} are again satisfied in
the point $\zz_0=1+\omega$ for the fields $(v,w,u)$. Continuing
this sort of arguments, we prove all the necessary relations which
involve the points $\zz=k+\ell\omega$ with $k\le 2$ or $\ell\le
2$, and which will be called the relations of the level 2. The
induction with respect to the level finishes the proof. \qed
\vspace{2mm}

It remains to prove Lemmas \ref{equidist} and \ref{geometry}
above.

What concerns the key Lemma \ref{geometry}, it might be
instructive to give two proofs for it, an analytic and a geometric
ones. The first one is shorter, but the second one seems to
provide more insight into the geometry. \vspace{1mm}

{\bf Analytic proof of Lemma \ref{geometry}.} We rewrite the
assumptions of the lemma as
\[
u_2-u_0=(u_1-u_0)e^{2i(\pi-\vartheta)}=(u_1-u_0)e^{-2i\vartheta}
\]
and
\[
v_2-v_0=c(v_1-v_0)e^{i\vartheta},\quad c>0.
\]

{\bf Geometric proof of Lemma \ref{geometry}.} The equations of
the $fgh$--system imply that the triangles $u_0u_1u_3$ and
$v_1v_3v_0$ are similar, and the triangles $u_0u_2u_3$ and
$v_2v_3v_0$ are similar. Therefore,
\[
\frac{|v_1-v_0|}{|u_0-u_3|}=\frac{|v_1-v_3|}{|u_0-u_1|},\qquad
\frac{|v_2-v_0|}{|u_0-u_3|}=\frac{|v_2-v_3|}{|u_0-u_2|}.
\]
From $|u_0-u_1|=|u_0-u_2|$ there follows now
\begin{equation}\label{lemma aux1}
\frac{|v_1-v_0|}{|v_1-v_3|}=\frac{|v_2-v_0|}{|v_2-v_3|}.
\end{equation} 
\begin{figure}[htbp]
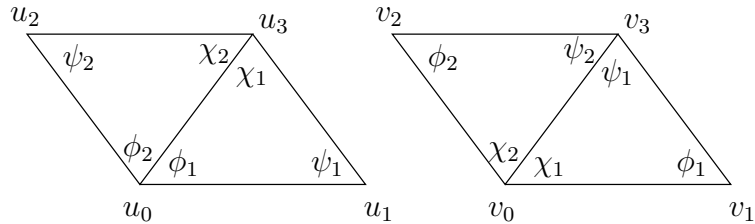

  \begin{center}
    \input a_proof17u.tex
    \input a_proof17v.tex
    \caption{To the proof of Lemma \ref{geometry}}
    \label{fig:proof17}
  \end{center}
\end{figure}

Denoting the angles as on Fig. \ref{fig:proof17}, we have:
\[
\chi_1+\chi_2=\vartheta,\qquad \phi_1+\phi_2=2\pi-2\vartheta,
\]
hence
\[
\psi_1+\psi_2=2\pi-(\phi_1+\phi_2)-(\chi_1+\chi_2)=\vartheta=\chi_1+\chi_2.
\]
In other words,
\begin{equation}\label{lemma aux2}
\measuredangle v_1v_3v_2=\measuredangle v_1v_0v_2.
\end{equation}
The relations (\ref{lemma aux1}), (\ref{lemma aux2}) yield that
the triangles $v_1v_3v_2$ and $v_1v_0v_2$ are similar. But they
have a common edge $[v_1,v_2]$, therefore they are congruent
(symmetric with respect to this edge). This implies that the
triangles $v_0v_2v_3$ and $v_0v_1v_3$ are isosceles, so that
$\chi_1=\psi_1$ and $\chi_2=\psi_2$, and
\[
|v_0-v_1|=|v_3-v_1|,\qquad |v_0-v_2|=|v_3-v_2|.
\]
Therefore
\[
|u_3-u_0|=|u_1-u_0|=|u_2-u_0|.
\]
Lemma is proved. \qed \vspace{2mm}

As for Lemma \ref{equidist}, its statement is a small part of the
following theorem and its corollary.
\begin{theorem}\label{Th circular zalpha}
If $\beta=\alpha$, then the recurrent relations (\ref{recur zw}),
(\ref{recur f}), (\ref{recur g}) with $u(1)=v(1)=1$ can be solved
for $u(k)$, $v(k)$, $f(k)$, $g(k)$ $(k\ge 0)$ in a closed form:
\begin{equation}\label{u3k}
u(3k)=\frac{2k}{k+2\alpha}\,\Pi_1(k),\qquad
u(3k+1)=\frac{2k+2\alpha}{k+2\alpha}\,\Pi_1(k),\qquad
u(3k+2)=2\,\Pi_1(k),
\end{equation}
\begin{equation}\label{f for zalpha}
f(3k-1)=f(3k)=f(3k+1)=\frac{2\alpha}{k+2\alpha}\,\Pi_1(k),
\end{equation}
and
\begin{equation}\label{v3k}
v(3k-1)=\frac{k-\alpha}{k+\alpha}\,\Pi_2(k),\qquad
v(3k)=\frac{k}{k+\alpha}\,\Pi_2(k),\qquad v(3k+1)=\Pi_2(k),
\end{equation}
\begin{equation}\label{g for zalpha}
g(3k-2)=g(3k-1)=g(3k)=\frac{\alpha}{k+\alpha}\,\Pi_2(k),
\end{equation}
where
\begin{equation}\label{Pi}
\Pi_1(k)=\frac{(1+2\alpha)(2+2\alpha)\ldots(k+2\alpha)}
{(1-\alpha)(2-\alpha)\ldots(k-\alpha)},\qquad
\Pi_2(k)=\frac{(1+\alpha)(2+\alpha)\ldots(k+\alpha)}
{(1-2\alpha)(2-2\alpha)\ldots(k-2\alpha)}.
\end{equation}
\end{theorem}
{\bf Proof.} Elementary calculations show that the expressions
above satisfy the recurrent relations (\ref{recur zw}),
(\ref{recur f}), (\ref{recur g}) with $\beta=\alpha$, as well as
the initial conditions. The uniqueness of the solution yields the
statement. We remark that similar formulas can be found also in
the general case $\alpha\neq\beta$, however, the property
formulated in Lemma \ref{equidist} fails to hold in general. \qed
\begin{corollary}\label{Cor circular zalpha}
If $\beta=\alpha$, and $u(1)=v(1)=1$, then for the third field
$w(k)$, $h(k)$ $(k\ge 0)$ we have:
\begin{equation}\label{w3k}
w(3k-1)=\frac{k-1+2\alpha}{1-2\alpha}\,\Pi_3(k),\qquad
w(3k)=\frac{k}{1-2\alpha}\,\Pi_3(k),\qquad
w(3k+1)=\frac{k+1-2\alpha}{1-2\alpha}\,\Pi_3(k),
\end{equation}
\begin{equation}\label{h3k-1 for zalpha}
h(3k-1)=h(3k)=\Pi_3(k),\qquad
h(3k+1)=\frac{k+1-2\alpha}{k+\alpha}\,\Pi_3(k),
\end{equation}
where
\begin{equation}
\Pi_3(k)=\frac{(1-\alpha)(2-\alpha)\ldots(k-\alpha)}
{\alpha(1+\alpha)\ldots(k-1+\alpha)}\cdot
\frac{(1-2\alpha)(2-2\alpha)\ldots(k-2\alpha)}
{2\alpha(1+2\alpha)\ldots(k-1+2\alpha)}.
\end{equation}
\end{corollary}
{\bf Proof.} The formulas for $h(k)=(f(k)g(k))^{-1}$ follow from
(\ref{f for zalpha}), (\ref{g for zalpha}). The formulas for
$w(k)= w(k-1)+h(k-1)$ with $w(0)=0$ follow by induction. \qed

\setcounter{equation}{0}
\section{Discrete hexagonal $z^{\alpha}$ and $\log z$}
\label{Sect hexagonal z^a}

\begin{figure}[htbp]
  \begin{center}
    \includegraphics{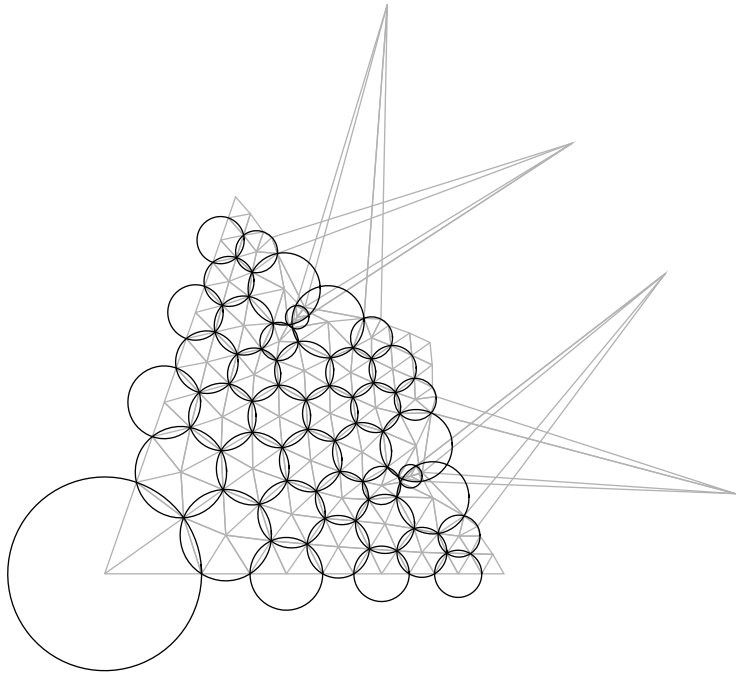}
    \caption{A non--immersed pattern with $\theta\neq 2\pi\alpha$.}
    \label{fig:nonReg}
  \end{center}
\end{figure}

Although the construction of the previous section always delivers
hexagonal circle patterns with $MR=-1$, these do not always behave
regularly. As a rule, they are not embedded (i.e. some elementary
triangles overlap), and even not immersed (i.e. some {\it
neighboring} triangles overlap), cf. Fig. \ref{fig:nonReg}).
However, there exists a choice of the initial values (i.e. of
$\theta$ in Theorem \ref{circular}) which assures that this is not
the case.

\begin{definition}
Let $0<\alpha=\beta<\frac{1}{2}$, so that $0<\gamma=1-2\alpha<1$.
Set $\theta=2\pi\alpha$. Then the hexagonal circle patterns of
Theorem \ref{circular} are called:
\begin{itemize}
\item[$u,v:$] the hexagonal $z^{3\alpha}$ with an intersection point
at the origin;
\item[$w:$] the hexagonal $z^{3\gamma}$ with a circle at the origin.
\end{itemize}
\end{definition}
In other words, for the hexagonal $z^{3\alpha}$ the opening angle
of the image of the sector (\ref{sector}) is equal to
$2\pi\alpha$, exactly as for the analytic function $z\mapsto
z^{3\alpha}$.
\begin{conjecture}\label{Conj z^a}
For $0<\alpha<\frac{1}{2}$ the hexagonal circle patterns
$z^{3\alpha}$ with an intersection point at the origin and
$z^{3\gamma}$ with a circle at the origin are embedded.
\end{conjecture}
For the proof of a similar statement for $z^{\alpha}$ circle
patterns with the combinatorics of the square grid see \cite{AB},
where it is proven that they are immersed. \vspace{2mm}

{\bf Remark.} Actually, the $u$ and $v$ versions of the hexagonal
$z^{3\alpha}$ with an intersection point at the origin are not
essentially different. Indeed, it is not difficult to see that the
half--sector of the $u$ pattern, corresponding to $0\le {\rm
arg}(\zz)\le \pi/3$, being rotated by $\pi\alpha$, coincides with
the half--sector of the $v$ pattern, corresponding to $\pi/3\le
{\rm arg}(\zz)\le 2\pi/3$, and vice versa. For the $w$ pattern, both
sectors are identical (up to the rotation by $\pi\gamma$). So, for every
$0<\alpha<\frac{1}{2}$ we have {\it two} essentially different hexagonal
pattrens $z^{3\alpha}$. \vspace{2mm}

\begin{figure}[htbp]
  \begin{center}
    \includegraphics{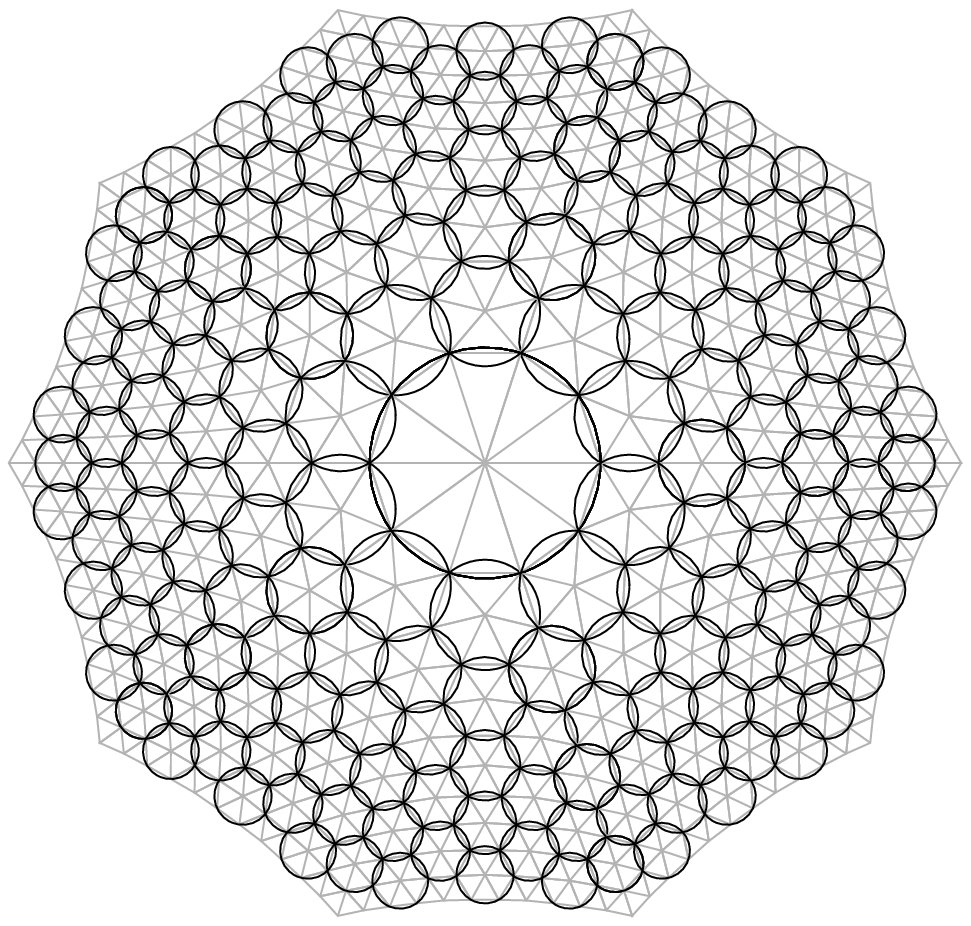}
    \includegraphics{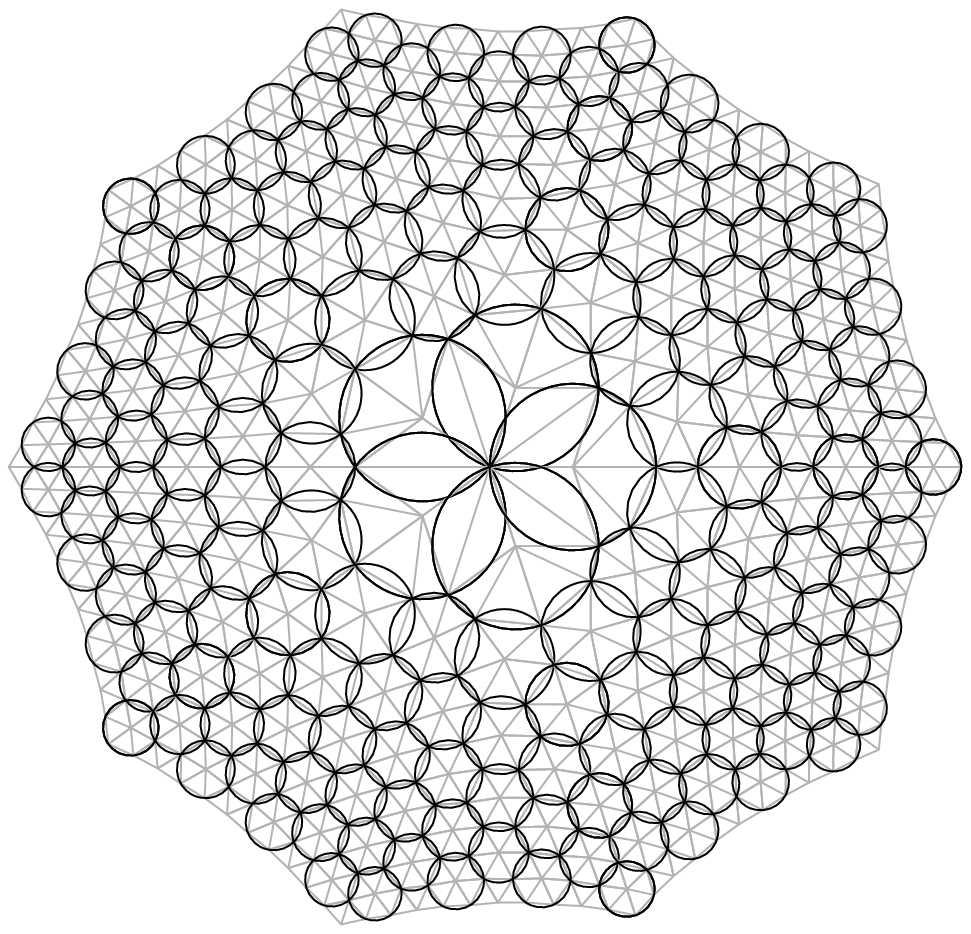}
    \caption{The hexagonal patterns $z^{3/5}$ with a circle at the
    origin and with an
    intersection point at the origin.}
    \label{fig:gamma15and25}
  \end{center}
\end{figure}

It is important to notice the peculiarity of the case when
$\alpha=n/N$ with $n,N\in{\Bbb N}$. Then one can attach to the
$u,v$--images of the sector $S$ its $N$ copies, rotated each time
by the angle $2\pi\alpha=2\pi n/N$.  The resulting object will
satisfy the conditions for the hexagonal circle pattern everywhere
except the origin $\zz=0$, which will be an intersection point of
$M=nN$ circles. Similarly, if $\gamma/2=n'/N'$, and we attach to
the $w$--image of the sector $S$ its $N'$ copies, rotated each
time by the angle $2\pi\gamma=4\pi n'/N'$, then the origin $\zz=0$
will be the center of a circle intersecting with $M'=n'N'$
neighboring circles. See Fig. \ref{fig:gamma15and25} for the
examples of the $w$--pattern with $\gamma=1/5$ and the $u$--pattern
with $\alpha=1/5$.

\vspace{2mm}

\begin{figure}[p]
  \begin{center}
    \centerline{\hspace{-0.9cm}\includegraphics[width = 0.4\hsize]{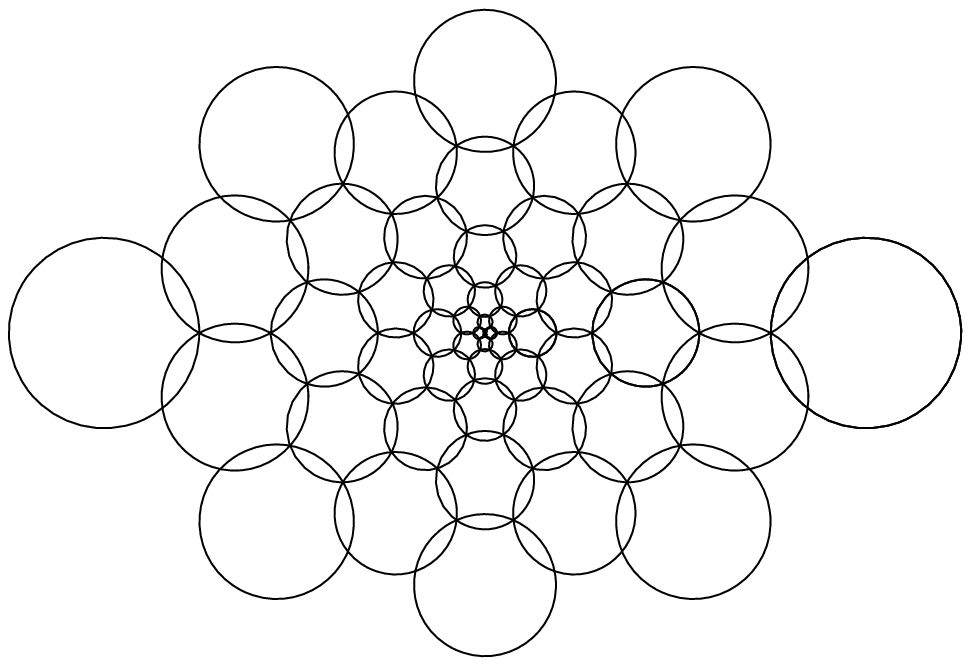}%
    \hspace{0.9cm}\raise 5cm\hbox{\includegraphics[angle=-90,width = 0.4\hsize]{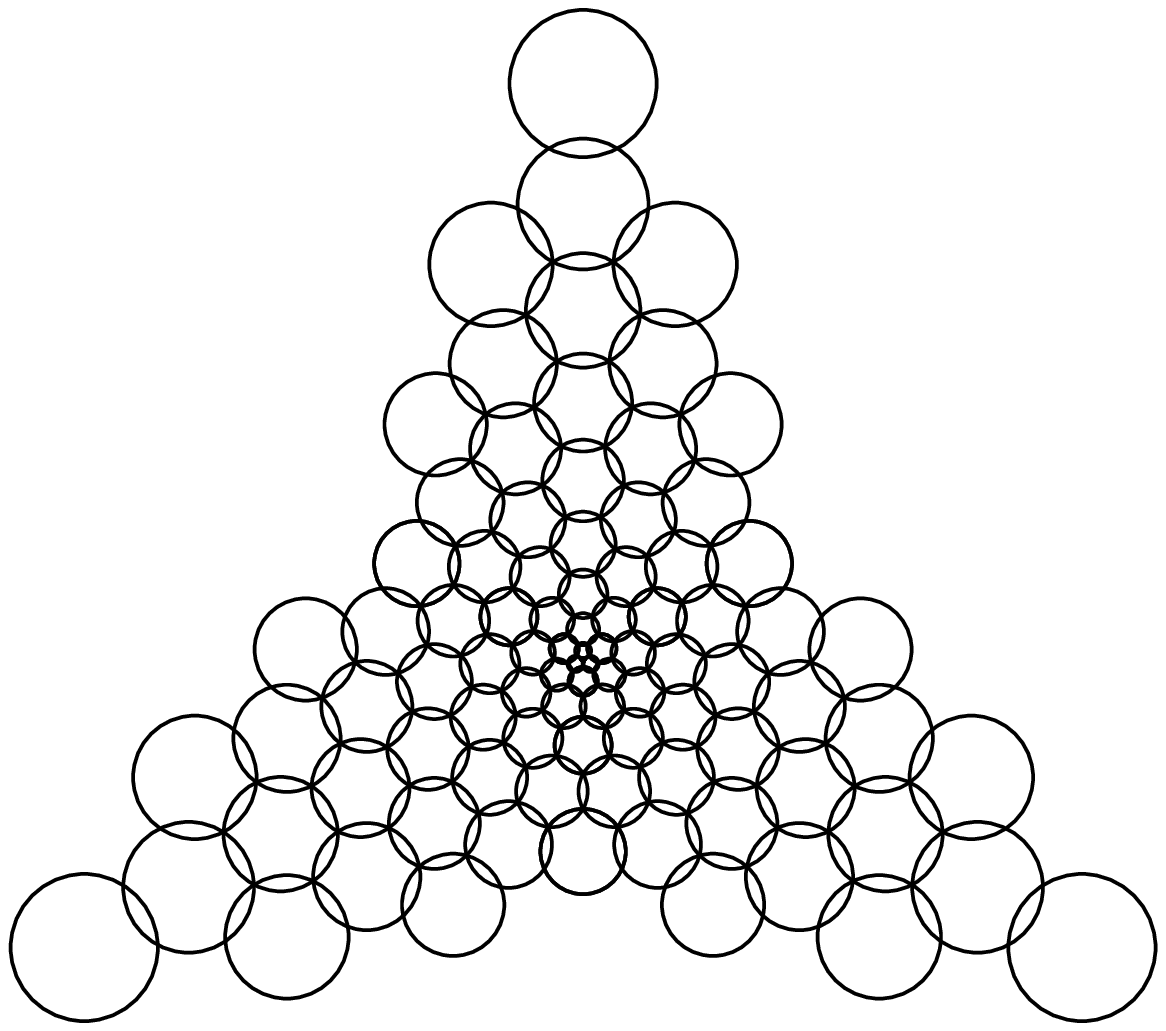}}}
\vspace{-2cm}
    \centerline{\includegraphics[width = 0.4\hsize]{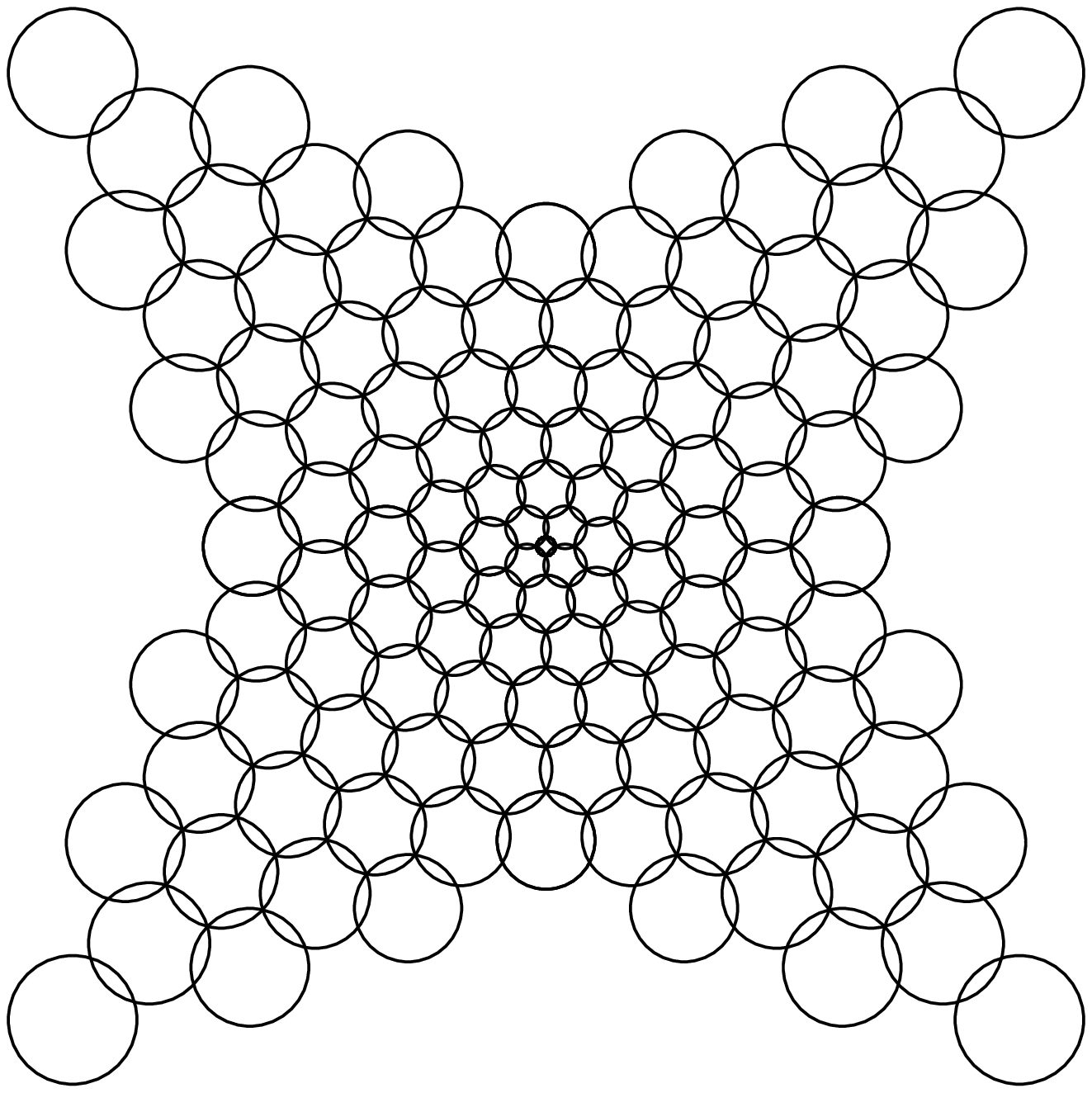}%
    \hbox{\raise 1.6cm\hbox{\includegraphics[width =
        0.4\hsize,angle=-18]{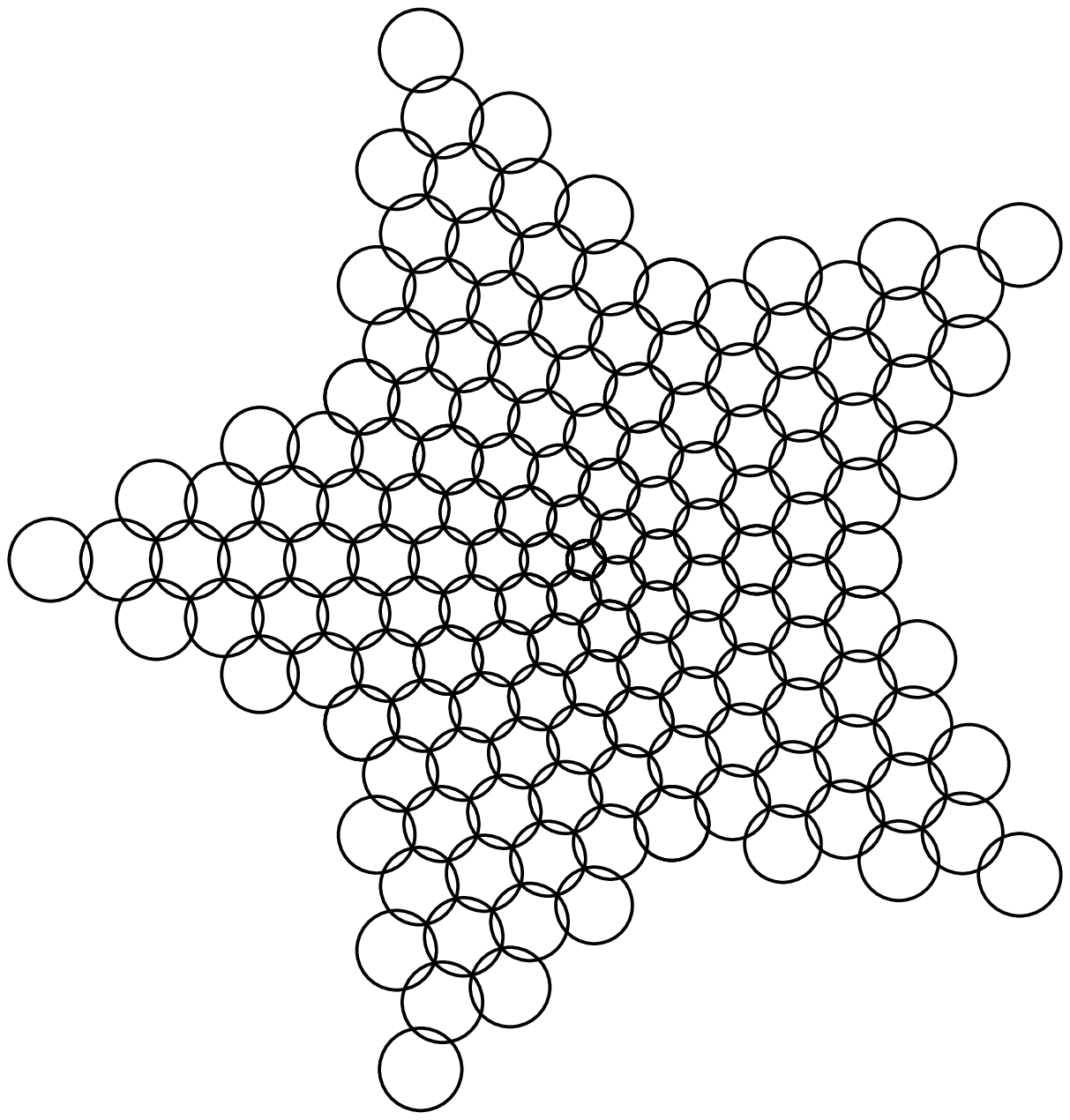}}}}%
\vspace{-2cm}
    \centerline{\includegraphics[width = 0.4\hsize]{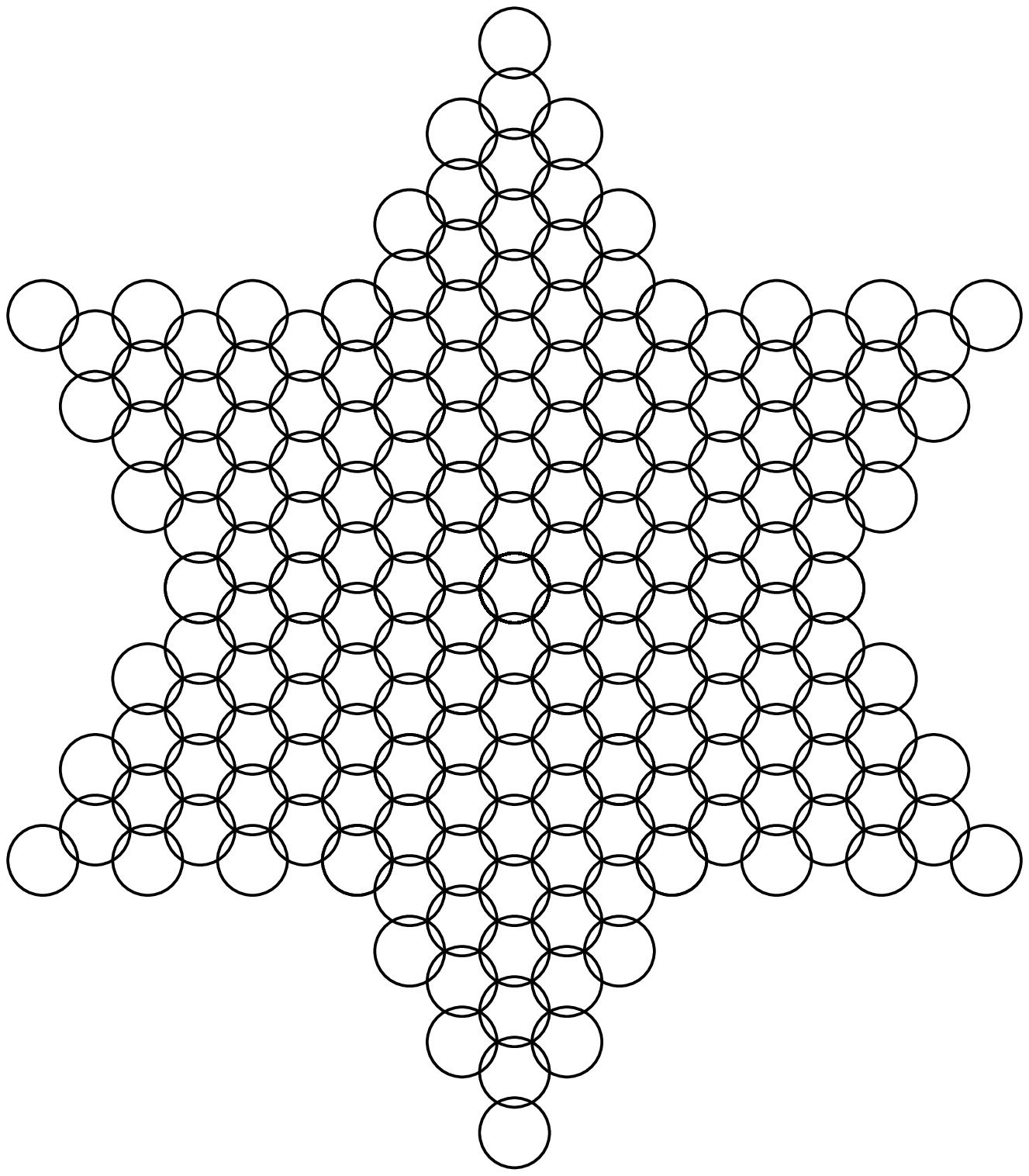}
    \includegraphics[width = 0.4\hsize,angle=12.8571]{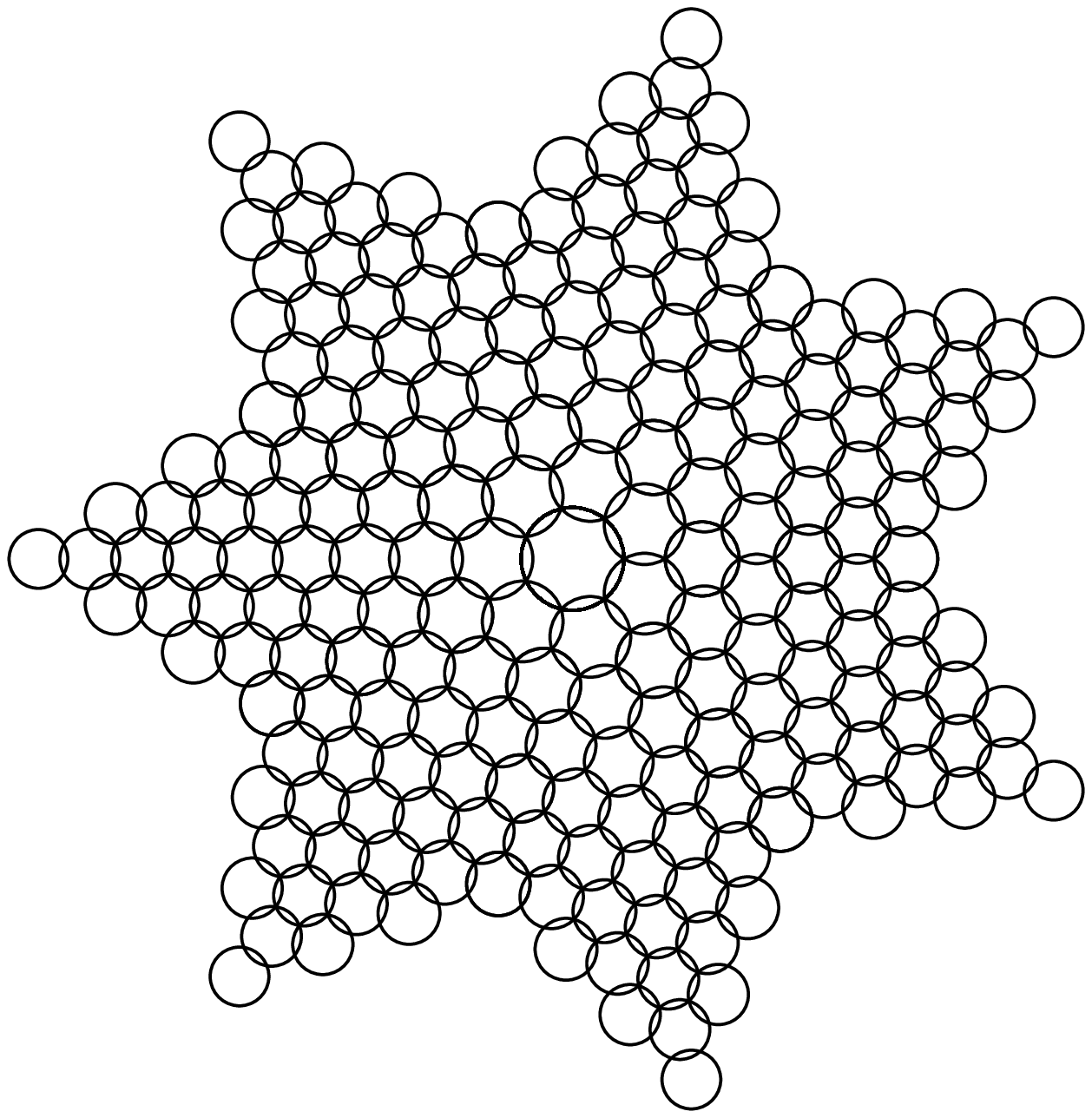}}%
    \caption{Some examples of $w$--pattern: $\gamma = 1, 2/3, 1/2,
    2/5, 1/3, 2/7$.}
    \label{fig:nice}
  \end{center}
\end{figure}
Now we turn our attention to the limiting cases $\alpha=1/2$ and
$\alpha=0$.

\subsection{ Case $\alpha=\frac{1}{2}$, $\gamma=0$: hexagonal $z^{3/2}$
and $\log z$} It is easy to see that the quantities $g(k)$, $k\ge
1$, and $v(k)$, $k\ge 2$, become singular as
$\alpha\to\frac{1}{2}$ (see (\ref{g for zalpha}) and (\ref{v3k})).
As a compensation, the quantities $h(k)$, $k\ge 1$, vanish with
$\alpha\to\frac{1}{2}$, so that $w(k)\to w(1)=1$ for all $k\ge 2$.
Similar effects hold for the $\ell$--axis, where $v(\ell\omega)$,
$\ell\ge 2$, become singular, and $w(\ell\omega)\to 1$ for all
$\ell\ge 1$. (Recall that for the $w$ pattern we have: 
$w(\omega)=e^{2\pi i\gamma}\to 1$). 
These observations suggest the following rescaling:
\begin{equation}\label{rescaling zalpha1}
u=\overset{\circ}{u},\qquad
v=\overset{\circ}{v}/(1-2\alpha),\qquad
w=1+(1-2\alpha)\overset{\circ}{w}.
\end{equation}
In order to be able to go to the limit $\alpha\to\frac{1}{2}$, we
have to calculate the values of our fields in several lattice
points next to $\zz=0$. Applying formulas (\ref{induct aux0}),
(\ref{induct aux1}), we find:
\begin{eqnarray}
u(0)=0, \quad u(1)=1, & &  u(\omega)=e^{2\pi i\alpha}, \qquad
u(1+\omega)=1+e^{2\pi i\alpha}, \label{u init}\\
v(0)=0,\quad v(1)=1, & &  v(\omega)=e^{2\pi i\alpha}, \qquad
v(1+\omega)=\frac{e^{2\pi i\alpha}}{1+e^{2\pi i\alpha}},  \label{v init}\\
w(0)=0,\quad w(1)=1, & &  w(\omega)=e^{2\pi i(1-2\alpha)}, \quad
w(1+\omega)=e^{\pi i(1-2\alpha)}.  \label{w init}
\end{eqnarray}
For the rescaled variables $\overset{\circ}{u}$,
$\overset{\circ}{v}$, $\overset{\circ}{w}$ in the limit
$\alpha\to\frac{1}{2}$ we find:
\begin{eqnarray}
\overset{\circ}{u}(0)=0, \quad \overset{\circ}{u}(1)=1, & &
\overset{\circ}{u}(\omega)=-1, \quad
\overset{\circ}{u}(1+\omega)=0,
\label{zalpha1 u init}\\
\overset{\circ}{v}(0)=0,\quad \overset{\circ}{v}(1)=0, & &
\overset{\circ}{v}(\omega)=0, \qquad
\overset{\circ}{v}(1+\omega)=\frac{i}{\pi},
\label{zalpha1 v init}\\
\overset{\circ}{w}(0)=\infty,\quad \overset{\circ}{w}(1)=0,  & &
\overset{\circ}{w}(\omega)=2\pi i, \quad
\overset{\circ}{w}(1+\omega)=\pi i. \label{zalpha1 w init}
\end{eqnarray}
These initial values have to be supplemented by the values in all
further points of the $k$-- and $\ell$--axes. From the formulas of
Theorem \ref{Th circular zalpha} there follows:
\begin{equation}\label{zalpha1 u}
\overset{\circ}{u}(3k)=\frac{2^k k!}{(2k-1)!!}\cdot (2k), \quad
\overset{\circ}{u}(3k+1)=\frac{2^k k!}{(2k-1)!!}\cdot (2k+1),
\quad \overset{\circ}{u}(3k+2)=\frac{2^k k!}{(2k-1)!!}\cdot
(2k+2),
\end{equation}
\begin{equation}\label{f for zalpha1}
\overset{\circ}{f}(3k-1)=\overset{\circ}{f}(3k)=\overset{\circ}{f}(3k+1)=
\frac{2^k k!}{(2k-1)!!},
\end{equation}
and
\begin{equation}\label{zalpha1 v}
\overset{\circ}{v}(3k-1)=\frac{(2k-1)!!}{2^k (k-1)!}\cdot (2k-1),
\quad \overset{\circ}{v}(3k)=\frac{(2k-1)!!}{2^k (k-1)!}\cdot
(2k), \quad \overset{\circ}{v}(3k+1)=\frac{(2k-1)!!}{2^k
(k-1)!}\cdot (2k+1),
\end{equation}
\begin{equation}\label{g for zalpha1}
\overset{\circ}{g}(3k-2)=\overset{\circ}{g}(3k-1)=\overset{\circ}{g}(3k)=
\frac{(2k-1)!!}{2^k (k-1)!},
\end{equation}
which have to be augmented by $\overset{\circ}{u}(k\omega)=
-\overset{\circ}{u}(k)$,
$\overset{\circ}{v}(k\omega)=-\overset{\circ}{v}(k)$.
From Corollary \ref{Cor circular zalpha} there follow the formulas for
the edges of the $\overset{\circ}{w}$ lattice:
\begin{eqnarray}
\overset{\circ}{h}(3k-1)=\overset{\circ}{h}(3k) & = & 
\overset{\circ}{h}((3k-1)\omega)=\overset{\circ}{h}(3k\omega)
\;\;=\;\; \frac{1}{k},\qquad k\ge 1,
\label{h3k-1 for zalpha1}\\
\overset{\circ}{h}(3k+1) & = & \overset{\circ}{h}((3k+1)\omega)\;\;=\;\;
\frac{1}{k+1/2}, \qquad k\ge 0.
\label{h3k+1 for zalpha1}
\end{eqnarray}
\begin{definition}
The hexagonal circle patterns corresponding to the solutions of the 
$fgh$--system in the sector (\ref{sector}) defined by
the boundary values (\ref{zalpha1 u init})--(\ref{h3k+1 for zalpha1})
 are called:
\begin{itemize}
\item[$\overset{\circ}{u}$, $\overset{\circ}{v}:$] the hexagonal $z^{3/2}$
with an intersection point at the origin;
\item[$\overset{\circ}{w}:$] the symmetric hexagonal $\log z$.
\end{itemize}
\end{definition}
Alternatively, one could define the lattices $\overset{\circ}{u}$,
$\overset{\circ}{v}$, $\overset{\circ}{w}$ as the solutions of the 
$fgh$--system with the initial values 
(\ref{zalpha1 u init})--(\ref{zalpha1 w init}), satisfying the
constraint (\ref{constr 1}), (\ref{constr 2}) with $\alpha=\beta=1/2$.
In this appoach the values (\ref{zalpha1 u})--(\ref{h3k+1 for zalpha1})
would be derived from the constraint. Notice also that the formulas 
(\ref{constr 3}), (\ref{constr 3 alt}) in this case turns into
\begin{eqnarray}
1 & = & k\frac{1}{\overset{\circ}{f}_0\overset{\circ}{g}_0+
\overset{\circ}{g}_0\overset{\circ}{f}_3+
\overset{\circ}{f}_3\overset{\circ}{g}_3}+
\ell\frac{1}{\overset{\circ}{f}_2\overset{\circ}{g}_2+
\overset{\circ}{g}_2\overset{\circ}{f}_5+
\overset{\circ}{f}_5\overset{\circ}{g}_5}+
m\frac{1}{\overset{\circ}{f}_4\overset{\circ}{g}_4+
\overset{\circ}{g}_4\overset{\circ}{f}_1+
\overset{\circ}{f}_1\overset{\circ}{g}_1}\\
& = & k\frac{\overset{\circ}{h}_0\overset{\circ}{f}_0\overset{\circ}{h}_3}
{\overset{\circ}{h}_0\overset{\circ}{f}_0+
\overset{\circ}{f}_0\overset{\circ}{h}_3+
\overset{\circ}{h}_3\overset{\circ}{f}_3}+
      \ell\frac{\overset{\circ}{h}_2\overset{\circ}{f}_2\overset{\circ}{h}_5}
{\overset{\circ}{h}_2\overset{\circ}{f}_2+
\overset{\circ}{f}_2\overset{\circ}{h}_5+
\overset{\circ}{h}_5\overset{\circ}{f}_5}+
       m\frac{\overset{\circ}{h}_4\overset{\circ}{f}_4\overset{\circ}{h}_1}
{\overset{\circ}{h}_4\overset{\circ}{f}_4+
\overset{\circ}{f}_4\overset{\circ}{h}_1+
\overset{\circ}{h}_1\overset{\circ}{f}_1}.
\end{eqnarray}

\begin{figure}[htbp]
  \begin{center}
    \includegraphics{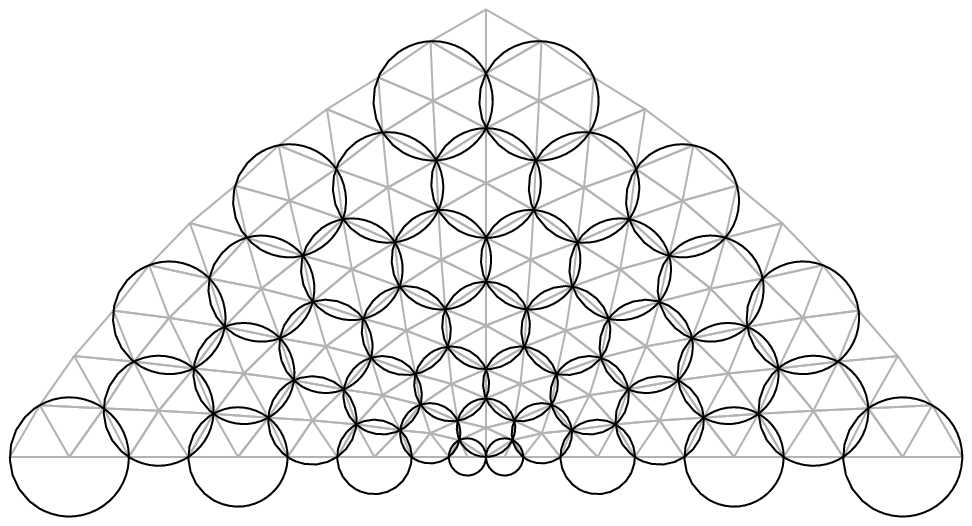}
    \includegraphics{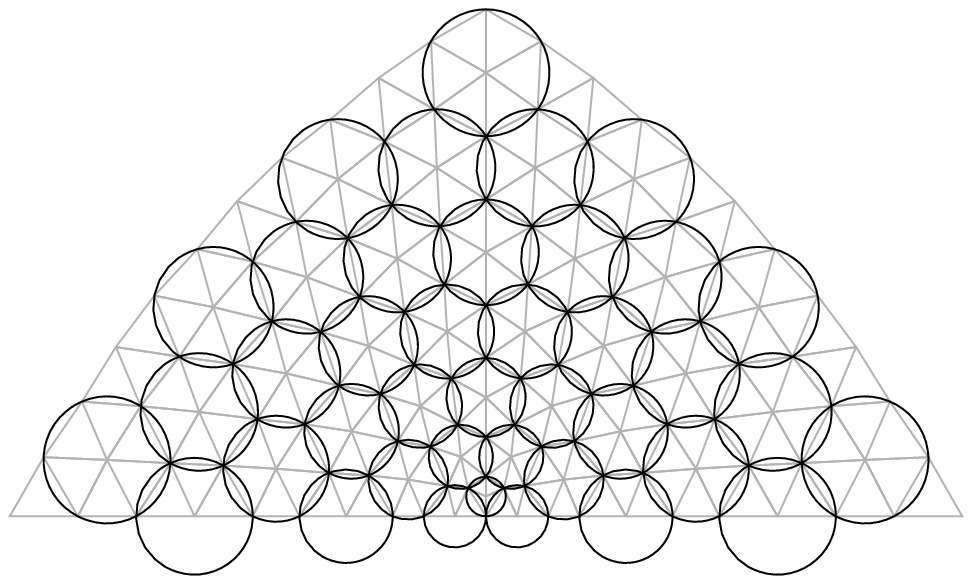}
    \includegraphics{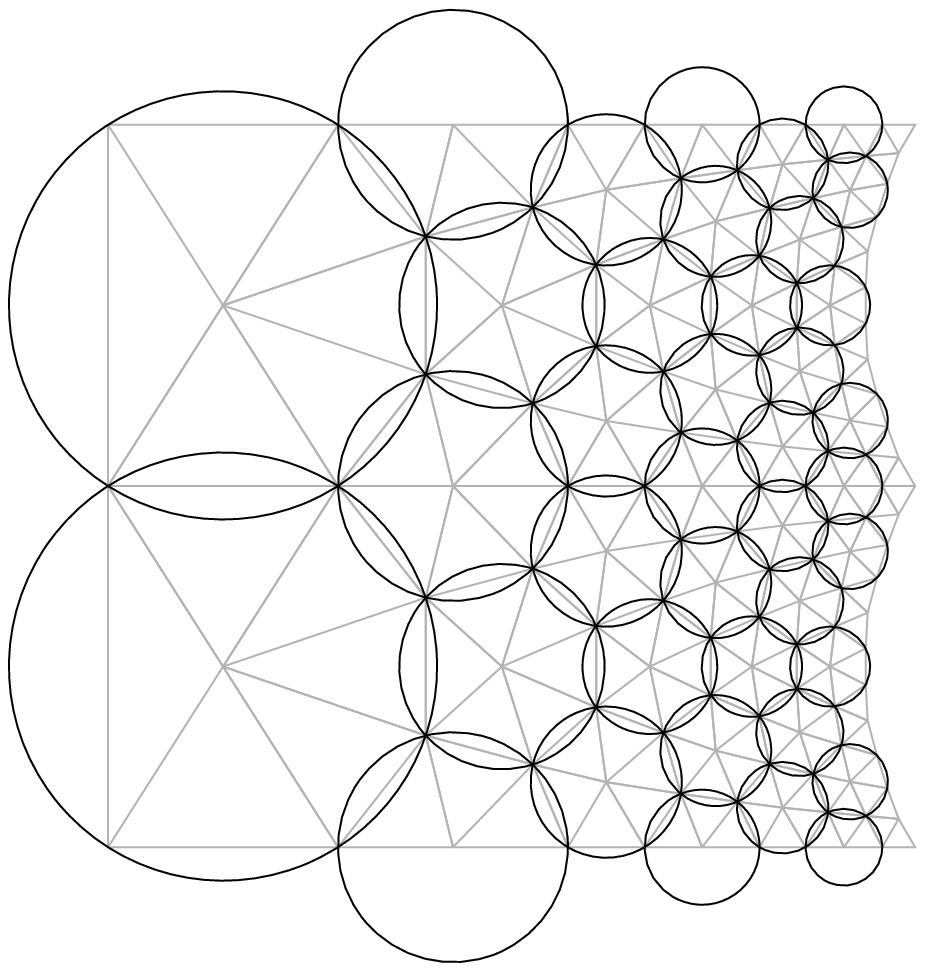}
    \caption{The patterns $z^{3/2}$ with an intersection point at the origin,
    and the symmetric hexagonal $\log z$; the second pattern coincides
    with the first one upon the rotation by $\pi/2$}
    \label{fig:gamma0}
  \end{center}
\end{figure}

\subsection{ Case $\alpha=0$, $\gamma=1$: hexagonal $\log z$ and $z^3$}

Considerations similar to those of the previous subsection show
that, as $\alpha\to 0$, the quantities $h(k)$, $k\ge 1$, and
$w(k)$, $k\ge 2$, become singular (see (\ref{h3k-1 for zalpha})
and (\ref{w3k})). As a compensation, the quantities $f(k)$, $k\ge
2$, and $g(k)$, $k\ge 1$, vanish with $\alpha\to 0$, so that
$u(k)\to u(2)=2$ for all $k\ge 3$, and $v(k)\to v(1)=1$ for all
$k\ge 2$. Similar effects hold for the $\ell$--axis. These
observations suggest the following rescaling:
\begin{equation}\label{rescaling zalpha2}
u=2+2\alpha\overset{\circ}{u},\qquad
v=1+\alpha\overset{\circ}{v},\qquad
w=\overset{\circ}{w}/(2\alpha^2).
\end{equation}
It turns out that in this case we need to calculate the values of
these functions in a larger number of lattice points in the
vicinity of $\zz=0$. To this end, we add to (\ref{u
init})--(\ref{w init}) the following values, which are obtained by
a direct calculation:
\begin{eqnarray}
u(2)=2, \qquad u(2\omega)=2e^{2\pi i\alpha},  & &
u(2+\omega)=\frac{1+e^{2\pi i\alpha}}{1+\alpha(e^{2\pi
i\alpha}-1)},
\label{u init begin}\\
u(1+2\omega)=\frac{1+e^{2\pi i\alpha}}{1+\alpha(e^{-2\pi
i\alpha}-1)}, & &
u(2+2\omega)=\frac{1-\alpha}{1-2\alpha}\,(1+e^{2\pi i\alpha}),
\end{eqnarray}
\begin{eqnarray}
v(2)=\frac{1-\alpha}{1-2\alpha},\qquad
v(2\omega)=\frac{1-\alpha}{1-2\alpha}\,e^{2\pi i\alpha},  & &
v(2+\omega)=\frac{1}{1+\alpha(e^{-2\pi i\alpha}-1)}, \\
v(1+2\omega)=\frac{e^{2\pi i\alpha}}{1+\alpha(e^{2\pi
i\alpha}-1)}, & & v(2+2\omega)=\frac{2e^{2\pi i\alpha}}{1+e^{2\pi
i\alpha}}
\end{eqnarray}
\begin{eqnarray}
w(2)=\frac{1-\alpha}{\alpha}, \qquad
w(2\omega)=\frac{1-\alpha}{\alpha}\,
e^{-2\pi i\alpha}, & &    w(2+\omega)=-\frac{1}{\alpha(e^{2\pi i\alpha}-1)}, \\
w(1+2\omega)=\frac{e^{-2\pi i\alpha}}{\alpha(e^{2\pi i\alpha}-1)},
& & w(2+2\omega)=-\frac{1-\alpha}{\alpha}\,e^{-2\pi i\alpha} .
\label{w init end}
\end{eqnarray}
From (\ref{u init})--(\ref{w init}) and 
(\ref{u init begin})--(\ref{w init end}) we obtain in the limit $\alpha\to 0$
under the rescaling (\ref{rescaling zalpha2}) the following initial values:
\begin{eqnarray}
\overset{\circ}{u}(0)=\infty, \quad \overset{\circ}{u}(1)=\infty, \quad
\overset{\circ}{u}(\omega)=\infty, \quad \overset{\circ}{u}(2)=0, 
\quad \overset{\circ}{u}(2\omega)=2\pi i,
\label{zalpha2 u init1}\\
\overset{\circ}{u}(1+\omega)=\pi i,  \quad
\overset{\circ}{u}(2+\omega)=\pi i,\quad \overset{\circ}{u}(1+2\omega)=\pi i,
\quad  \overset{\circ}{u}(2+2\omega)=1+\pi i,
\label{zalpha2 u init2}
\end{eqnarray}
\begin{eqnarray}
\overset{\circ}{v}(0)=\infty, \quad \overset{\circ}{v}(1)=0, \quad
\overset{\circ}{v}(\omega)=2\pi i, \quad \overset{\circ}{v}(2)=1, 
\quad \overset{\circ}{v}(2\omega)=1+2\pi i,
\label{zalpha2 v init1}\\
\overset{\circ}{v}(1+\omega)=\infty,  \quad
\overset{\circ}{v}(2+\omega)=0, \quad \overset{\circ}{v}(1+2\omega)=2\pi i,
\quad  \overset{\circ}{v}(2+2\omega)=\pi i,
\label{zalpha2 v init2}
\end{eqnarray}
\begin{eqnarray}
\overset{\circ}{w}(0)=0, \quad \overset{\circ}{w}(1)=0, \quad 
\overset{\circ}{w}(\omega)=0, \quad \overset{\circ}{w}(2)=0, 
\quad \overset{\circ}{w}(2\omega)=0,
\label{zalpha2 w init1}\\
\overset{\circ}{w}(1+\omega)=0,  \quad 
\overset{\circ}{w}(2+\omega)=\frac{i}{\pi}, 
\quad \overset{\circ}{w}(1+2\omega)=-\frac{i}{\pi},
\quad  \overset{\circ}{w}(2+2\omega)=0.
\label{zalpha2 w init2}
\end{eqnarray}
These initial values have to be supplemented by the values in all further
points of the $k$-- and $\ell$--axes. From the formulas of Theorem
\ref{Th circular zalpha} there follow the expressions for the edges 
of the lattices $\overset{\circ}{u}$, $\overset{\circ}{v}$:
\begin{eqnarray}
\overset{\circ}{f}(3k-1)=\overset{\circ}{f}(3k)=\overset{\circ}{f}(3k+1)
\;\;=\;\; \overset{\circ}{f}((3k-1)\omega)=\overset{\circ}{f}(3k\omega)
=\overset{\circ}{f}((3k+1)\omega) & = & 
\frac{1}{k}, \qquad k\ge 1,\nonumber\\\label{f for zalpha2}\\
\overset{\circ}{g}(3k-2)=\overset{\circ}{g}(3k-1)=\overset{\circ}{g}(3k)
\;\;=\;\; \overset{\circ}{g}((3k-2)\omega)=\overset{\circ}{g}((3k-1)\omega)=
\overset{\circ}{g}(3k\omega) & = & 
\frac{1}{k}, \qquad k\ge 1.\nonumber\\ \label{g for zalpha2}
\end{eqnarray}
The formulas of Corollary \ref{Cor circular zalpha} yield the results for
the lattice $\overset{\circ}{w}$:
\begin{equation}\label{zalpha2 w}
\overset{\circ}{w}(3k)=k^3, \quad 
\overset{\circ}{w}(3k+1)=k^2(k+1),  \quad 
\overset{\circ}{w}(3k+2)=k(k+1)^2, \quad k\ge 1, 
\end{equation}
so that
\begin{equation}\label{h for zalpha2}
\overset{\circ}{h}(3k-1)=\overset{\circ}{h}(3k)=k^2,\quad
\overset{\circ}{h}(3k+1)=k(k+1), \quad k\ge 1. 
\end{equation}
Of course, one has also $\overset{\circ}{w}(k\omega)=\overset{\circ}{w}(k)$.

\begin{figure}[htbp]
  \begin{center}
    \includegraphics[width=0.45\hsize]{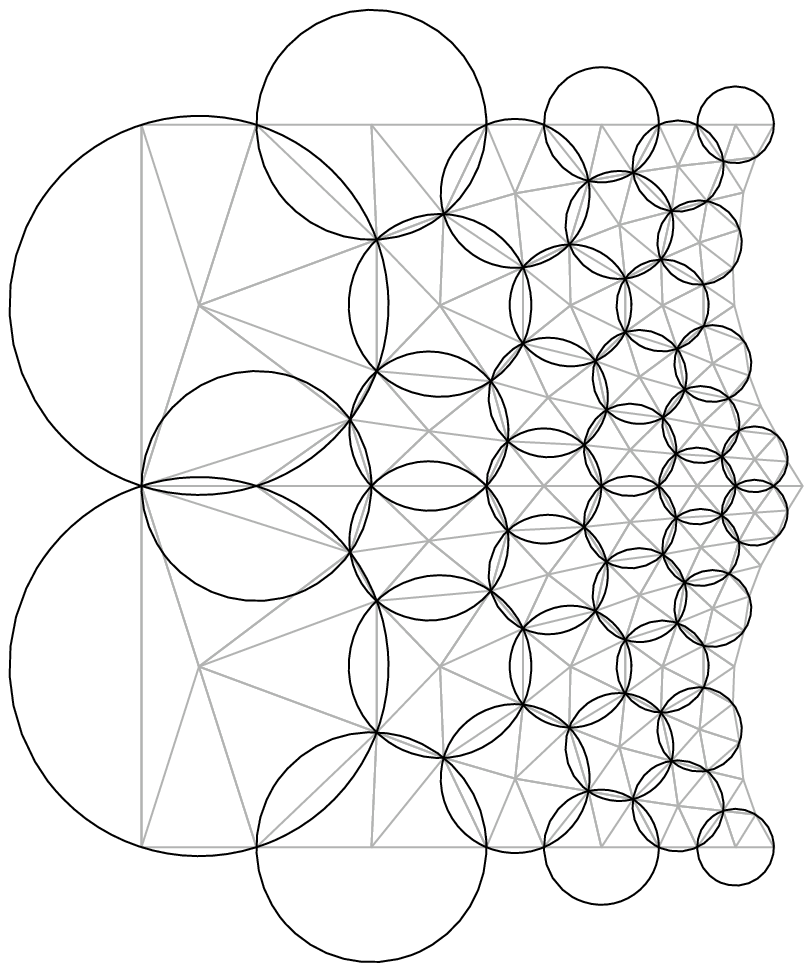}
    \includegraphics[width=0.45\hsize]{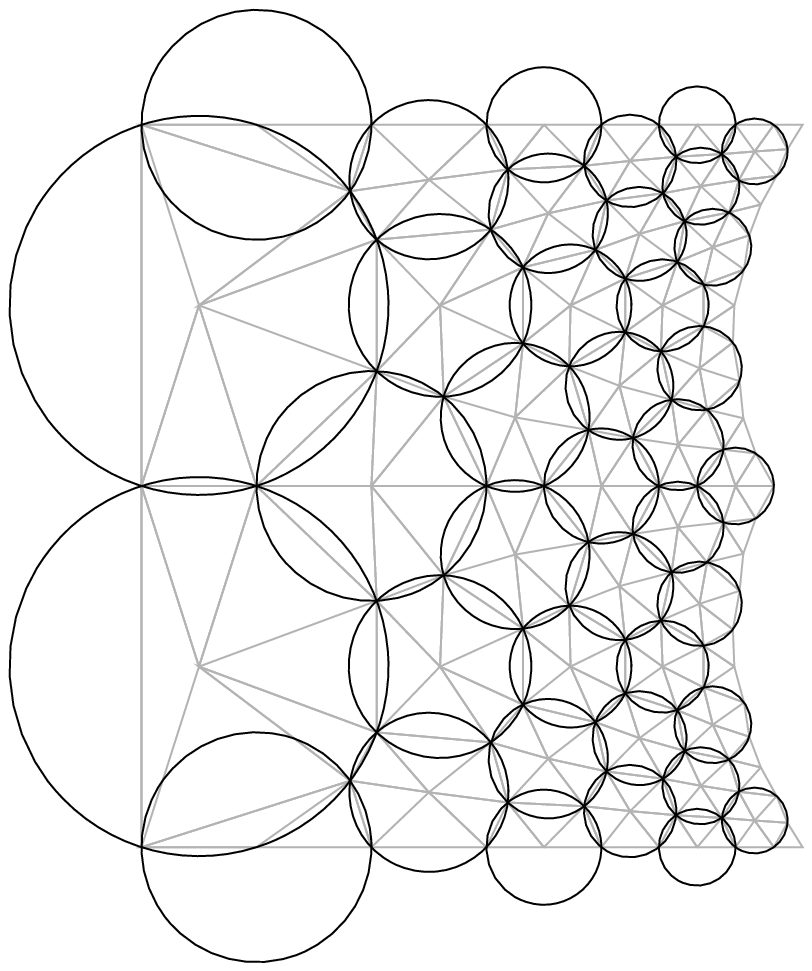}
    \includegraphics{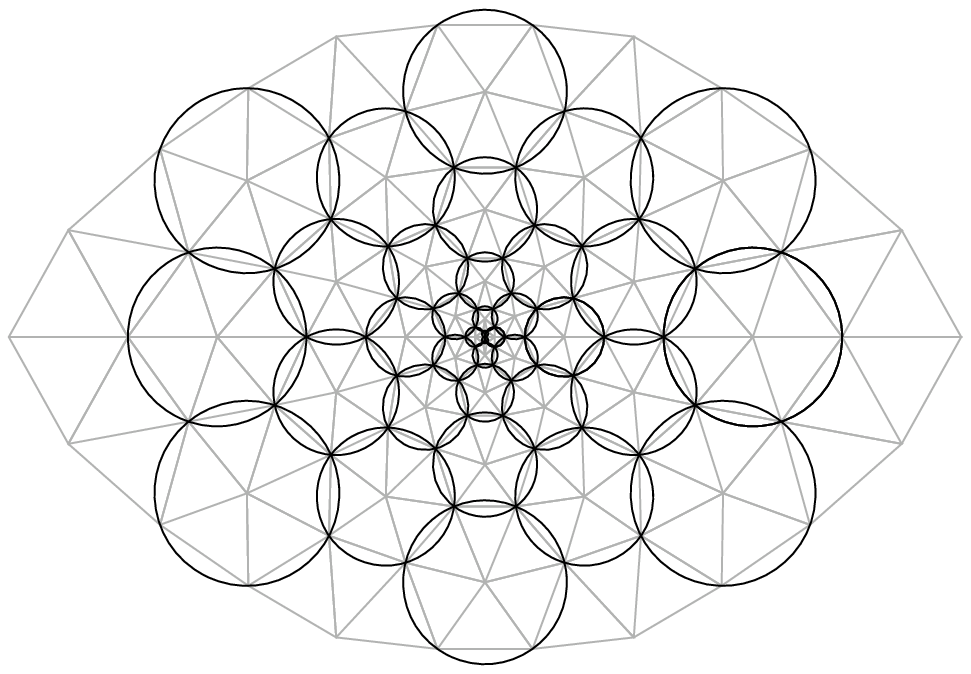}
    \caption{The asymmetric patterns $\log z$ and the hexagonal
    pattern $z^3$ with a circle at the origin; the upper half of the first 
    pattern coincides with the lower half of the second one, and vice versa}
    \label{fig:gamma1}
  \end{center}
\end{figure}

\begin{definition}
The hexagonal circle patterns corresponding to the solutions of
the $fgh$--system in the sector (\ref{sector}) defined by the
boundary values (\ref{zalpha2 u init1})--(\ref{h for zalpha2}) are
called:
\begin{itemize}
\item[$\overset{\circ}{u}$, $\overset{\circ}{v}:$] the  asymmetric
 hexagonal $\log z$;
\item[$\overset{\circ}{w}:$] the hexagonal $z^3$ with a (degenerate) circle
at the origin.
\end{itemize}
\end{definition}
It is meant that the $u$--image of the half-sector $0\le {\rm arg}(\zz)\le 
\pi/3$ is not symmetric with respect to the line $\Im(u)=\pi i/2$ (the image
of ${\rm arg}(\zz)=\pi/6$, and the same for $v$. Instead, this symmetry
interchanges the $u$ pattern and the $v$ pattern, see Fig. \ref{fig:gamma1}.

Alternatively, one can define these lattices as the solutions of
the $fgh$--system with the initial values (\ref{zalpha2 u
init1})--(\ref{zalpha2 w init2}), satisfying the constraint
(\ref{constr 1}), (\ref{constr 2}), which in the present situation
degenerates into
\begin{eqnarray}
1 & = &
k\frac{\overset{\circ}{f}_0\overset{\circ}{g}_0\overset{\circ}{f}_3}
{\overset{\circ}{f}_0\overset{\circ}{g}_0+
\overset{\circ}{g}_0\overset{\circ}{f}_3+
\overset{\circ}{f}_3\overset{\circ}{g}_3}+
      \ell\frac{\overset{\circ}{f}_2\overset{\circ}{g}_2\overset{\circ}{f}_5}
{\overset{\circ}{f}_2\overset{\circ}{g}_2+
\overset{\circ}{g}_2\overset{\circ}{f}_5+
\overset{\circ}{f}_5\overset{\circ}{g}_5}+
       m\frac{\overset{\circ}{f}_4\overset{\circ}{g}_4\overset{\circ}{f}_1}
{\overset{\circ}{f}_4\overset{\circ}{g}_4+
\overset{\circ}{g}_4\overset{\circ}{f}_1+
\overset{\circ}{f}_1\overset{\circ}{g}_1},\\
1 & = &
k\frac{\overset{\circ}{g}_0\overset{\circ}{f}_3\overset{\circ}{g}_3}
{\overset{\circ}{f}_0\overset{\circ}{g}_0+
\overset{\circ}{g}_0\overset{\circ}{f}_3+
\overset{\circ}{f}_3\overset{\circ}{g}_3}+
      \ell\frac{\overset{\circ}{g}_2\overset{\circ}{f}_5\overset{\circ}{g}_5}
{\overset{\circ}{f}_2\overset{\circ}{g}_2+
\overset{\circ}{g}_2\overset{\circ}{f}_5+
\overset{\circ}{f}_5\overset{\circ}{g}_5}+
       m\frac{\overset{\circ}{g}_4\overset{\circ}{f}_1\overset{\circ}{g}_1}
{\overset{\circ}{f}_4\overset{\circ}{g}_4+
\overset{\circ}{g}_4\overset{\circ}{f}_1+
\overset{\circ}{f}_1\overset{\circ}{g}_1}.
\end{eqnarray}
Just as in the non--degenerate case, these formulas allow one to
calculate inductively the values of $\overset{\circ}{u}$,
$\overset{\circ}{v}$ on the $k$-- and $\ell$--axes. The formulas
(\ref{constr 3}), (\ref{constr 3 alt}) hold literally with
$\gamma=1$.

\section{Conclusions}

In this paper we introduced the notion of hexagonal circle
patterns, and studied in some detail a subclass consisting of
circle patterns with the property that six intersection points on
each circle have the multi-ratio $-1$. We established the
connection of this subclass with integrable systems on the regular
triangular lattice, and used this connection to describe some
B\"acklund--like transformations of hexagonal circle patterns
(transformation $u\mapsto v\mapsto w$, see Theorems \ref{z to w},
\ref{from circ to circ}), and to find discrete analogs of the
functions $z^{\alpha}$, $\log z$. Of course, this is only the
beginning of the story of hexagonal circle patterns. In a
subsequent publication we shall demonstrate that there exists
another subclass related to integrable systems, namely the
patterns with fixed intersection angles. The intersection of both
subclasses constitute conformally symmetric patterns, including
analogs of Doyle's spirals (cf. \cite{BH}).

A very interesting question is, what part of the theory of
integrable circle patterns can be applied to hexagonal circle
packings. This also will be a subject of our investigation.

This research was financially supported by DFG
(Sonderforschungsbereich 288 ``Differential Geometry and Quantum
Physics'').


\renewcommand{\thesection}{\Alph{section}}

\setcounter{section}{0}

\setcounter{equation}{0}
\section{Appendix: Square lattice version of the $fgh$--system}
\label{square grid version}

Dropping all edges of $E(\cT\cL)$ parallel to the $m$--axis, we
end up with the cell complex isomorphic to the regular square
lattice: its vertices $\zz=k+\ell\omega$ may be identified with
$(k,\ell)\in{\Bbb Z}^2$, its edges are then identified with those
pairs $[(k_1,\ell_1),(k_2,\ell_2)]$ for which
$|k_1-k_2|+|\ell_1-\ell_2|=1$, and its 2-cells (parallelograms)
are identified with the elementary squares of the square lattice.
Hence, flat connections on $\cT\cL$ form a subclass of flat
connections on the square lattice. A natural question is, whether
this inclusion is strict, i.e. whether there exist flat
connections on the square lattice which cannot be extended to flat
connections on $\cT\cL$. At least for the $fgh$--system, the
answer is negative: denote by $\cM\subset{\rm
SL}(3,{\Bbb C}) [\lambda]$ the set of matrices (\ref{L}), then
flat connections on the regular square grid with values in $\cM$
are essentially in a one-to-one correspondence with flat
connections on $\cT\cL$ with values in $\cM$, i.e. with solutions
of the $fgh$--system. This is a consequence of the following
statement dealing with an elementary square of the regular square
lattice: a flat connection on such an elementary square with
values in $\cM$ can be extended by an element of $\cM$ sitting on
its diagonal without violating the flatness property. More
precisely:
\begin{lemma}\label{lemma 6to4}
Let
\[
L_1L_2=L_3L_4, \quad where \quad L_i\in\cM\;\;(i=1,2,3,4),
\]
and let the off--diagonal parts of $L_1$, $L_2$ be componentwise
distinct from the off--diagonal parts of $L_3$, $L_4$,
respectively. Then there exists $L_0\in\cM$ such that
\[
L_0L_1L_2=L_0L_3L_4=I.
\]
\end{lemma}
%
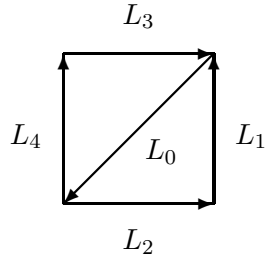
\begin{figure}[htbp]
  \begin{center}
\setlength{\unitlength}{1cm}
\begin{picture}(4,4)
\thicklines \put(1,1){\vector(1,0){2}} \put(3,1){\vector(0,1){2}}
\put(1,1){\vector(0,1){2}} \put(1,3){\vector(1,0){2}}
\put(3,3){\vector(-1,-1){2}} \put(3.5,1.9){\makebox(0,0){$L_1$}}
\put(0.5,1.9){\makebox(0,0){$L_4$}}
\put(2,0.5){\makebox(0,0){$L_2$}}
\put(2,3.5){\makebox(0,0){$L_3$}}
\put(2.3,1.7){\makebox(0,0){$L_0$}}
\end{picture}
    \caption{To Lemma \ref{lemma 6to4}.}
    \label{fig:aboutNijhoff}
  \end{center}
\end{figure}
{\bf Proof.} We have to prove that
$(L_1L_2)^{-1}=(L_3L_4)^{-1}\in\cM$. It is easy to see that it is
necessary and sufficient to prove that the entries 13, 21, 32 of
this matrix vanish, i.e. that there holds
\begin{equation}\label{lemma square to prove}
f_1g_1+f_2g_1+f_2g_2=f_3g_3+f_4g_3+f_4g_4=0,
\end{equation}
as well as two similar equations resulting by two successive
permutations $(f,g,h)\mapsto (g,h,f)$. We are given the relations
$f_ig_ih_i=1$ and
\begin{equation}\label{lemma square have 1}
f_1+f_2=f_3+f_4,\qquad g_1+g_2=g_3+g_4,\qquad h_1+h_2=h_3+h_4,
\end{equation}
\begin{equation}\label{lemma square have 2}
f_1g_2=f_3g_4,\qquad g_1h_2=g_3h_4,\qquad h_1f_2=h_3f_4.
\end{equation}
In order to prove (\ref{lemma square to prove}), we start with the
third equation in (\ref{lemma square have 1}):
\begin{equation}\label{lemma square aux1}
h_1\left(1-\frac{h_3}{h_1}\right)=h_2\left(\frac{h_4}{h_2}-1\right).
\end{equation}
Using $f_ig_ih_i=1$ and (\ref{lemma square have 2}), we find:
\begin{equation}\label{lemma square aux2}
\frac{h_3}{h_1}=\frac{f_1g_1}{f_3g_3}=\frac{g_4g_1}{g_2g_3},\qquad
\frac{h_4}{h_2}=\frac{g_1}{g_3}.
\end{equation}
Plugging this into (\ref{lemma square aux1}), we get:
\begin{equation}\label{lemma square aux3}
\frac{g_2g_3-g_1g_4}{f_1g_1g_2g_3}=\frac{g_1-g_3}{f_2g_2g_3}.
\end{equation}
Now, due to the second equation in (\ref{lemma square have 1}), we
find:
\begin{equation}\label{lemma square aux4}
g_2g_3-g_1g_4=g_2(g_3-g_1)+g_1(g_2-g_4)=(g_1+g_2)(g_3-g_1).
\end{equation}
Substituting this into (\ref{lemma square aux3}), we come to the
equation:
\begin{equation}
(g_3-g_1)\left(\frac{g_1+g_2}{f_1g_1}+\frac{1}{f_2}\right)=0.
\end{equation}
Since, by condition, $g_1\neq g_3$, we obtain
$f_2(g_1+g_2)+f_1g_1=0$, which is the equation (\ref{lemma square
to prove}). \qed \vspace{2mm}

This result shows that the $fgh$--system could be alternatively
studied in a more common framework of integrable systems on a
square lattice. However, such an approach would hide a rich and
interesting geometric structures immanently connected with the
triangular lattice. It should be said at this point that the
one--field equation (\ref{hex eq z}) was first found, under the
name of the ``Schwarzian lattice Bussinesq equation'' by Nijhoff
in \cite{N} using a (different) Lax representation on the square
lattice. The same holds for the one--field form of the constraint
(\ref{constr 1 z}).

\setcounter{equation}{0}
\section{Appendix: Proofs of statements of Sect. \ref{Sect isomonodromic}}
\label{Appendix}

{\bf Proof of Proposition \ref{constraint preliminary}.} The
arguments are similar for both equations (\ref{constr 1}),
(\ref{constr 2}). For instance, for the first one we have to
demonstrate that
\begin{eqnarray}\label{constr 1 well aux1}
\frac{f_0g_0f_3}{f_0g_0+g_0f_3+f_3g_3}+
\frac{f_2g_2f_5}{f_2g_2+g_2f_5+f_5g_5}+
\frac{f_4g_4f_1}{f_4g_4+g_4f_1+f_1g_1} & = & \nonumber\\
\frac{f_0f_3}{f_0+f_3+f_3g_3/g_0}+
\frac{f_2f_5}{f_2+f_5+f_5g_5/g_2}+
\frac{f_4f_1}{f_4+f_1+f_1g_1/g_4} & = & 0.
\end{eqnarray}
To eliminate the fields $g$ from this equation, consider six
elementary triangles surrounding the vertex $\zz$. The equations
(\ref{motion eq}) imply:
\begin{eqnarray*}
\frac{g_1}{g_0}=-\frac{f_0+f_1}{f_1},\quad
\frac{g_2}{g_1}=-\frac{f_1}{f_1+f_2},\quad
\frac{g_3}{g_2}=-\frac{f_2+f_3}{f_3}, \\
\frac{g_5}{g_0}=-\frac{f_5+f_0}{f_5},\quad
\frac{g_4}{g_5}=-\frac{f_5}{f_4+f_5},\quad
\frac{g_3}{g_2}=-\frac{f_3+f_4}{f_3}.
\end{eqnarray*}
Therefore,
\begin{eqnarray}
f_0+f_3+f_3\frac{g_3}{g_0} & = &
f_0+f_3-\frac{(f_0+f_1)(f_2+f_3)}{f_1+f_2}=
\frac{(f_0-f_2)(f_1-f_3)}{f_1+f_2}\label{constr 1 well aux2}\\
& = & f_0+f_3-\frac{(f_5+f_0)(f_3+f_4)}{f_4+f_5}=
\frac{(f_4-f_0)(f_3-f_5)}{f_4+f_5}.\label{constr 1 well aux3}
\end{eqnarray}
By the way, this again yields the property $MR=-1$ of the lattice
$u$, which can be written now as
\begin{equation}\label{constr 1 well H}
(f_0+f_1)(f_2+f_3)(f_4+f_5)=(f_1+f_2)(f_3+f_4)(f_5+f_0),
\end{equation}
Using (\ref{constr 1 well aux2}), an analogous expression along
the $\ell$--axis, and an expression analogous to (\ref{constr 1
well aux3}) along the $m$--axis, we rewrite (\ref{constr 1 well
aux1}) as
\begin{equation}
\frac{f_0f_3(f_1+f_2)}{(f_0-f_2)(f_1-f_3)}+
\frac{f_2f_5(f_3+f_4)}{(f_2-f_4)(f_3-f_5)}+
\frac{f_4f_1(f_2+f_3)}{(f_2-f_4)(f_1-f_3)}=0.
\end{equation}
Clearing denominators, we put it in the equivalent form
\begin{eqnarray*}
f_0f_3(f_1+f_2)(f_2-f_4)(f_3-f_5)+
f_2f_5(f_3+f_4)(f_0-f_2)(f_1-f_3) & & \\
+f_4f_1(f_2+f_3)(f_0-f_2)(f_3-f_5) & = & 0.
\end{eqnarray*}
But the polynomial on the left--hand side of the last formula is
equal to
\[
f_2f_3\Big((f_1+f_2)(f_3+f_4)(f_5+f_0)-(f_0+f_1)(f_2+f_3)(f_4+f_5)\Big),
\]
and hence vanishes in virtue of (\ref{constr 1 well H}). \qed
\vspace{2mm}

{\bf Proof of Proposition \ref{third constraint}.} Denote the
right--hand sides of (\ref{constr 1}), (\ref{constr 2}),
(\ref{constr 3}) through $U(\zz)$, $V(\zz)$, $W(\zz)$,
respectively. In order to prove (\ref{constr 3}), i.e. $\gamma
w=W(\zz)$, it is necessary and sufficient to demonstrate that
\[
\gamma h_0=W(\widetilde{\zz})-W(\zz),\quad \gamma
h_2=W(\widehat{\zz})-W(\zz), \quad \gamma h_4=W(\bar{\zz})-W(\zz),
\]
(or, actually, any two of these three equations). We perform the
proof for the first one only, since for the other two everything
is similar. In dealing with our constraints we are free to choose
any representative $(k,\ell,m)$ for $\zz$. In order to keep things
shorter, we always assume in this proof that $m=0$. Writing the
formula
\[
\gamma=\frac{1}{h_0}\Big(W(\widetilde{\zz})-W(\zz)\Big)
\]
in long hand, we have to prove that
\begin{eqnarray}\label{constr 3 to prove}
\gamma=1-\alpha-\beta & = &
(k+1)\frac{1/h_0}{\widetilde{f}_0\widetilde{g}_0+
\widetilde{g}_0\widetilde{f}_3+\widetilde{f}_3\widetilde{g}_3}-
k\frac{1/h_0}{f_0g_0+g_0f_3+f_3g_3}\nonumber\\
& & +\ell \frac{1/h_0}{\widetilde{f}_2\widetilde{g}_2+
\widetilde{g}_2\widetilde{f}_5+\widetilde{f}_5\widetilde{g}_5}-
\ell\frac{1/h_0}{f_2g_2+g_2f_5+f_5g_5}.
\end{eqnarray}
Assuming that (\ref{constr 1}) and (\ref{constr 2}) hold, we have:
\[
\alpha+\beta=\frac{1}{f_0}\Big(U(\widetilde{\zz})-U(\zz)\Big)+
\frac{1}{g_0}\Big(V(\widetilde{\zz})-V(\zz)\Big).
\]
Taking into account that $\widetilde{f}_3=f_0$,
$\widetilde{g}_3=g_0$, we find:
\begin{eqnarray*}
\alpha+\beta & = & (k+1)\frac{\widetilde{f}_0\widetilde{g}_0+
\widetilde{g}_0\widetilde{f}_3}{\widetilde{f}_0\widetilde{g}_0+
\widetilde{g}_0\widetilde{f}_3+\widetilde{f}_3\widetilde{g}_3}-
k\frac{g_0f_3+f_3g_3}{f_0g_0+g_0f_3+f_3g_3}\nonumber\\
& & +\ell\left(
\frac{\widetilde{f}_2\widetilde{g}_2\widetilde{f}_5/f_0+
\widetilde{g}_2\widetilde{f}_5\widetilde{g}_5/g_0}
{\widetilde{f}_2\widetilde{g}_2+
\widetilde{g}_2\widetilde{f}_5+\widetilde{f}_5\widetilde{g}_5}-
\frac{f_2g_2f_5/f_0+g_2f_5g_5/g_0}{f_2g_2+g_2f_5+f_5g_5}\right),
\end{eqnarray*}
or, equivalently,
\begin{eqnarray}
\gamma=1-\alpha-\beta & = &
(k+1)\frac{\widetilde{f}_3\widetilde{g}_3}
{\widetilde{f}_0\widetilde{g}_0+
\widetilde{g}_0\widetilde{f}_3+\widetilde{f}_3\widetilde{g}_3}-
k\frac{f_0g_0}{f_0g_0+g_0f_3+f_3g_3}\nonumber\\
& & -\ell\left(
\frac{\widetilde{f}_2\widetilde{g}_2\widetilde{f}_5/f_0+
\widetilde{g}_2\widetilde{f}_5\widetilde{g}_5/g_0}
{\widetilde{f}_2\widetilde{g}_2+
\widetilde{g}_2\widetilde{f}_5+\widetilde{f}_5\widetilde{g}_5}-
\frac{f_2g_2f_5/f_0+g_2f_5g_5/g_0}{f_2g_2+g_2f_5+f_5g_5}\right).
\end{eqnarray}
The first two terms on the right--hand side already have the
required form, since
$\widetilde{f}_3\widetilde{g}_3=f_0g_0=1/h_0$. So, it remains to
prove that
\begin{equation}\label{third constraint aux0}
-\frac{\widetilde{f}_2\widetilde{g}_2\widetilde{f}_5/f_0+
\widetilde{g}_2\widetilde{f}_5\widetilde{g}_5/g_0}
{\widetilde{f}_2\widetilde{g}_2+
\widetilde{g}_2\widetilde{f}_5+\widetilde{f}_5\widetilde{g}_5}+
\frac{f_2g_2f_5/f_0+g_2f_5g_5/g_0}{f_2g_2+g_2f_5+f_5g_5}=
\frac{1/h_0}{\widetilde{f}_2\widetilde{g}_2+
\widetilde{g}_2\widetilde{f}_5+\widetilde{f}_5\widetilde{g}_5}-
\frac{1/h_0}{f_2g_2+g_2f_5+f_5g_5}.
\end{equation}
The most direct and unambiguous way to do this is to notice that
everything here may be expressed with the help of the
$fgh$--equations in terms of a single field $h$. After
straightforward calculations one obtains:
\begin{eqnarray}
\widetilde{f}_2\widetilde{g}_2\frac{\widetilde{f}_5}{f_0}+
\widetilde{f}_5\widetilde{g}_5\frac{\widetilde{g}_2}{g_0} & = &
-\frac{1}{\widetilde{h}_5}+\frac{h_0(h_0-\widetilde{h}_5)}{\widetilde{h}_2
\widetilde{h}_5h_5},  \label{third constraint aux1}\\
f_2g_2\frac{f_5}{f_0}+f_5g_5\frac{g_2}{g_0} & = &
-\frac{1}{h_2}+\frac{h_0(h_0-h_2)}{\widetilde{h}_2h_2h_5},
\label{third constraint aux2}\\
\widetilde{f}_2\widetilde{g}_2+\widetilde{g}_2\widetilde{f}_5+
\widetilde{f}_5\widetilde{g}_5 & = &
\frac{(h_0-\widetilde{h}_5)(\widetilde{h}_4-\widetilde{h}_2)}
{\widetilde{h}_2\widetilde{h}_5h_5},  \label{third constraint aux3}\\
f_2g_2+g_2f_5+f_5g_5 & = &
\frac{(h_0-h_2)(h_1-h_5)}{\widetilde{h}_2h_2h_5}. \label{third
constraint aux4}
\end{eqnarray}
Taking into account that
$\widetilde{h}_4-\widetilde{h}_2=h_1-h_5$, we see that (\ref{third
constraint aux0}) and Proposition \ref{third constraint} are
proved. \qed \vspace{2mm}

{\bf Proof of Theorem \ref{monodromy}.} In order for the
isomonodromy property to hold, the following compatibility
conditions of (\ref{wave evolution in mu}) with (\ref{eq in mu})
are necessary and sufficient: (\ref{zero curv in mu}) and
\renewcommand{\arraystretch}{1.9}
\begin{equation}
\left\{\begin{array}{l}
\displaystyle\frac{d}{d\mu}\cL(\ee_0,\mu)=\cA_{k+1,\ell,m}\cL(\ee_0,\mu)-
\cL(\ee_0,\mu)\cA_{k,\ell,m},\\
\displaystyle\frac{d}{d\mu}\cL(\ee_2,\mu)=\cA_{k,\ell+1,m}\cL(\ee_2,\mu)-
\cL(\ee_2,\mu)\cA_{k,\ell,m},\\
\displaystyle\frac{d}{d\mu}\cL(\ee_4,\mu)=\cA_{k,\ell,m+1}\cL(\ee_4,\mu)-
\cL(\ee_4,\mu)\cA_{k,\ell,m}.
\end{array}\right.
\end{equation}
Substituting the ansatz (\ref{ans A}) and calculating the residues
at $\mu=-1$, $\mu=0$ and $\mu=\infty$, we see that the above
system is equivalent to the following nine matrix equations:
\begin{eqnarray}
C_{k+1,\ell,m}\cL(\ee_0,-1) & = & \cL(\ee_0,-1)C_{k,\ell,m},  \label{eq C1}\\
C_{k,\ell+1,m}\cL(\ee_2,-1) & = & \cL(\ee_2,-1)C_{k,\ell,m},  \label{eq C2}\\
C_{k,\ell,m+1}\cL(\ee_4,-1) & = & \cL(\ee_4,-1)C_{k,\ell,m},
\label{eq C3}
\end{eqnarray}
\begin{eqnarray}
D(\widetilde{\zz})\cL(\ee_0,0) & = & \cL(\ee_0,0)D(\zz),  \label{eq D1}\\
D(\widehat{\zz})\cL(\ee_2,0) & = & \cL(\ee_2,0)D(\zz),  \label{eq D2}\\
D(\bar{\zz})\cL(\ee_4,0) & = & \cL(\ee_4,0)D(\zz),  \label{eq D3}
\end{eqnarray}
\begin{eqnarray}
\Big(C_{k+1,\ell,m}+D(\widetilde{\zz})\Big)Q-Q\Big(C_{k,\ell,m}+D(\zz)\Big)
& =& Q, \label{eq CD1}\\
\Big(C_{k,\ell+1,m}+D(\widehat{\zz})\Big)Q-Q\Big(C_{k,\ell,m}+D(\zz)\Big)
& = & Q, \label{eq CD2}\\
\Big(C_{k,\ell,m+1}+D(\bar{\zz})\Big)Q-Q\Big(C_{k,\ell,m}+D(\zz)\Big)
& = & Q, \label{eq CD3}
\end{eqnarray}
where
\renewcommand{\arraystretch}{1}
\begin{equation}\label{matr Q}
Q=\left(\begin{array}{ccc} 0 & 0 & 0 \\ 0 & 0 & 0 \\ 1 & 0 & 0
\end{array} \right).
\end{equation}
We do not aim at solving these equations completely, but rather at
finding a {\it certain} solution leading to the constraint
(\ref{constr 1}), (\ref{constr 2}). The subsequent reasoning will
be divided into several steps.

{\bf Step 1. Consistency of the ansatz for} $C_{k,\ell,m}$. First
of all, we have to convince ourselves that the ansatz (\ref{C}),
(\ref{P}) does not violate the necessary condition (\ref{cover
cond for A}), i.e. that
\begin{equation}\label{projector wonder}
P_2+P_4+P_6=I.
\end{equation}
Notice that the entries 12 and 23 of this matrix equation are
nothing but the content of Proposition \ref{constraint
preliminary}. Upon the cyclic permutation of the fields
$(f,g,h)\mapsto(g,h,f)$ this gives also the entry 31. To check the
entry 21, we proceed as in the proof of Proposition
\ref{constraint preliminary}. We have to prove that
\begin{eqnarray*}
\frac{g_0}{f_0g_0+g_0f_3+f_3g_3}+
\frac{g_2}{f_2g_2+g_2f_5+f_5g_5}+
\frac{g_4}{f_4g_4+g_4f_1+f_1g_1} & = & \nonumber\\
\frac{1}{f_0+f_3+f_3g_3/g_0}+ \frac{1}{f_2+f_5+f_5g_5/g_2}+
\frac{1}{f_4+f_1+f_1g_1/g_4} & = & \nonumber\\
\frac{f_1+f_2}{(f_0-f_2)(f_1-f_3)}+
\frac{f_3+f_4}{(f_2-f_4)(f_3-f_5)}+
\frac{f_2+f_3}{(f_2-f_4)(f_1-f_3)} & = & 0.
\end{eqnarray*}
Clearing denominators, we put it in the equivalent form
\[
(f_1+f_2)(f_2-f_4)(f_3-f_5)+ (f_3+f_4)(f_0-f_2)(f_1-f_3)+
(f_2+f_3)(f_0-f_2)(f_3-f_5)=0.
\]
But the polynomial on the left--hand side is equal to
\[
(f_1+f_2)(f_3+f_4)(f_5+f_0)-(f_0+f_1)(f_2+f_3)(f_4+f_5),
\]
and vanishes due to (\ref{constr 1 well H}). Via the cyclic
permutation of fields this proves also the entries 32 and 13 of
the matrix identity (\ref{projector wonder}). Finally, turning to
the diagonal entries, we consider, for the sake of definiteness,
the entry 22. We have to prove that
\begin{eqnarray*}
\frac{f_3g_0}{f_0g_0+g_0f_3+f_3g_3}+
\frac{f_5g_2}{f_2g_2+g_2f_5+f_5g_5}+
\frac{f_1g_4}{f_4g_4+g_4f_1+f_1g_1} & = & \nonumber\\
\frac{f_3}{f_0+f_3+f_3g_3/g_0}+ \frac{f_5}{f_2+f_5+f_5g_5/g_2}+
\frac{f_1}{f_4+f_1+f_1g_1/g_4} & = & \nonumber\\
\frac{f_3(f_1+f_2)}{(f_0-f_2)(f_1-f_3)}+
\frac{f_5(f_3+f_4)}{(f_2-f_4)(f_3-f_5)}+
\frac{f_1(f_2+f_3)}{(f_2-f_4)(f_1-f_3)} & = & 1,
\end{eqnarray*}
or
\begin{eqnarray*}
f_3(f_1+f_2)(f_2-f_4)(f_3-f_5)+ f_5(f_3+f_4)(f_0-f_2)(f_1-f_3)
 & + & \nonumber\\
f_1(f_2+f_3)(f_0-f_2)(f_3-f_5)-(f_0-f_2)(f_1-f_3)(f_2-f_4)(f_3-f_5)
& = & 0.
\end{eqnarray*}
Again, the polynomial on the left--hand side is equal to
\[
f_3\Big((f_1+f_2)(f_3+f_4)(f_5+f_0)-(f_0+f_1)(f_2+f_3)(f_4+f_5)\Big),
\]
and vanishes due to (\ref{constr 1 well H}). The formula
(\ref{projector wonder}) is proved.

{\bf Step 2. Checking the equations for the matrix}
$C_{k,\ell,m}$. Next, we have to show that the ansatz (\ref{C}),
(\ref{P}) verifies (\ref{eq C1})--(\ref{eq C3}). Notice that the
matrices
\[
\cL(\ee,-1)=\left(\begin{array}{ccc} 1 & f & 0 \\ 0 & 1 & g \\ -h
& 0 & 1
\end{array}\right)
\]
are degenerate, and that
\[
\xi=\left(\begin{array}{c} fg \\ -g \\ 1\end{array}\right) \quad
{\rm and} \quad \eta^{\rm T}=\Big(1,\;\; -f,\;\; fg\Big)
\]
are the right null--vector and the left null--vector of
$\cL(\ee,-1)$, respectively. In terms of these vectors one can
write the projectors $P_{0,2,4}$ as
\[
P_j=\frac{1}{\langle \xi_j,\eta_{j+3}\rangle}\xi_j\eta_{j+3}^{\rm
T},\quad j=0,2,4.
\]
Therefore we have:
\begin{eqnarray}
P_0(\widetilde{\zz})\cL(\ee_0,-1)=\cL(\ee_0,-1)P_0(\zz) & = & 0,\\
P_2(\widehat{\zz})\cL(\ee_2,-1)=\cL(\ee_2,-1)P_2(\zz) & = & 0,\\
P_4(\bar{\zz})\cL(\ee_4,-1)=\cL(\ee_4,-1)P_4(\zz) & = & 0.
\end{eqnarray}
In order to demonstrate (\ref{eq C1})--(\ref{eq C3}) it is
sufficient to prove that
\begin{eqnarray}
P_2(\widetilde{\zz})\cL(\ee_0,-1)=\cL(\ee_0,-1)P_2(\zz), & \quad &
P_4(\widetilde{\zz})\cL(\ee_0,-1)=\cL(\ee_0,-1)P_4(\zz)=0,\\
P_4(\widehat{\zz})\cL(\ee_2,-1)=\cL(\ee_2,-1)P_4(\zz), & \quad &
P_0(\widehat{\zz})\cL(\ee_2,-1)=\cL(\ee_2,-1)P_0(\zz)=0,\\
P_0(\bar{\zz})\cL(\ee_4,-1)=\cL(\ee_4,-1)P_0(\zz), & \quad &
P_2(\bar{\zz})\cL(\ee_4,-1)=\cL(\ee_4,-1)P_2(\zz)=0.
\end{eqnarray}
All these equations are verified in a similar manner, therefore we
restrict ourselves to the first one.
\[
\frac{1}{\langle \widetilde{\xi}_2,\widetilde{\eta}_5\rangle}
\widetilde{\xi}_2\widetilde{\eta}_5^{\rm T}\cL(\ee_0,-1)=
\frac{1}{\langle \xi_2,\eta_5\rangle}\cL(\ee_0,-1)\xi_2\eta_5^{\rm
T},
\]
or, in long hand,
\begin{equation}\label{C ansatz aux0}
\frac{1}{\widetilde{f}_2\widetilde{g}_2+\widetilde{g}_2\widetilde{f}_5+
\widetilde{f}_5\widetilde{g}_5}\left(\begin{array}{c}
\widetilde{f}_2\widetilde{g}_2 \\ -\widetilde{g}_2 \\
1\end{array}\right)\!
\left(\begin{array}{c} 1-h_0/\widetilde{h}_5 \\
f_0-\widetilde{f}_5 \\
\widetilde{f}_5(\widetilde{g}_5-g_0)\end{array} \right)^{\rm
T}\!\!= \frac{1}{f_2g_2+g_2f_5+f_5g_5} \left(\begin{array}{c}
(f_2-f_0)g_2 \\ g_0-g_2 \\ 1-h_0/h_2\end{array}\right)
\!\left(\begin{array}{c} 1 \\ -f_5 \\ f_5g_5\end{array}
\right)^{\rm T}\!\!.
\end{equation}
To prove this we have, first, to check that these two rank one
matrices are proportional, and then to check that their entries 31
(say) coincide. The second of these claims reads:
\begin{equation}\label{C ansatz aux1}
\frac{1-h_0/\widetilde{h}_5}
{\widetilde{f}_2\widetilde{g}_2+\widetilde{g}_2\widetilde{f}_5+
\widetilde{f}_5\widetilde{g}_5}=
\frac{1-h_0/h_2}{f_2g_2+g_2f_5+f_5g_5},
\end{equation}
and follows from (\ref{third constraint aux3}), (\ref{third
constraint aux4}). The first claim above is equivalent to:
\[
\left(\begin{array}{c} \widetilde{f}_2\widetilde{g}_2 \\
-\widetilde{g}_2 \\ 1\end{array}\right) \sim
\left(\begin{array}{c} (f_2-f_0)g_2 \\ g_0-g_2 \\
1-h_0/h_2\end{array}\right) \quad{\rm and}\quad
\left(\begin{array}{c} 1-h_0/\widetilde{h}_5 \\
f_0-\widetilde{f}_5 \\
\widetilde{f}_5(\widetilde{g}_5-g_0)\end{array} \right)\sim
\left(\begin{array}{c} 1 \\ -f_5 \\ f_5g_5\end{array} \right),
\]
which, in turn, is equivalent to:
\begin{equation}\label{C ansatz aux2}
\widetilde{f}_2=g_2\frac{f_0-f_2}{g_0-g_2},\quad
\widetilde{g}_2=h_2\frac{g_0-g_2}{h_0-h_2},
\end{equation}
and
\begin{equation}\label{C ansatz aux3}
\widetilde{h}_5=f_5\frac{h_0-\widetilde{h}_5}{f_0-\widetilde{f}_5},\quad
\widetilde{f}_5=g_5\frac{f_0-\widetilde{f}_5}{g_0-\widetilde{g}_5}.
\end{equation}
All these relations easily follow from the equations of the
$fgh$--system. For instance, to check the first equation in
(\ref{C ansatz aux2}), one has to consider the two elementary
positively oriented triangles $(\zz,\zz+\omega,\zz+\varepsilon)$
and $(\zz,\zz+1,\zz+\varepsilon)$. Denoting the edge
$\ee_{12}=(\zz+\omega,\zz+\varepsilon)$, we have:
\[
f_2+f_{12}=f_0+\widetilde{f}_2\;(=-f_1),\qquad
f_{12}g_2=\widetilde{f}_2g_0\;(=f_1g_1).
\]
Eliminating $f_{12}$ from these two equations, we end up with the
desired one. This finishes the proof of (\ref{eq C1})--(\ref{eq
C3}).

{\bf Step 3. Checking the equations for the matrix} $D(\zz)$.
Notice that the matrices
\[
L(\ee,0)=\left(\begin{array}{ccc} 1 & f & 0 \\ 0 & 1 & g \\ 0 & 0
& 1
\end{array}\right)
\]
are upper triangular. We require that the matrices $D(\zz)$ are
also upper triangular:
\[
D=\left(\begin{array}{ccc} d_{11} & d_{12} & d_{13} \\ 0 & d_{22} & d_{23} \\
0 & 0 & d_{33}\end{array}\right).
\]
It is immediately seen that the diagonal entries are constants. By
multiplying the wave function $\Psi_{k,\ell,m}(\mu)$ from the
right by a constant ($\mu$--dependent) matrix one can arrange that
the matrices $D(\zz)$ are traceless. Hence the diagonal part of
$D$ is parameterized by two arbitrary numbers. It will be
convenient to choose this parametrization as
\[
(d_{11},d_{22},d_{33})=\Big(-(2\alpha+\beta)/3,\;
(\alpha-\beta)/3, \; (2\beta+\alpha)/3\Big).
\]
Equating the entries 12 and 23 in (\ref{eq D1})--(\ref{eq D3}), we
find for an arbitrary positively oriented edge
$\ee=(\zz_1,\zz_2)\in E(\cT\cL)$:
\begin{eqnarray*}
d_{12}(\zz_2)-d_{12}(\zz_1) & = & (d_{22}-d_{11})f=
\alpha\Big(u(\zz_2)-u(\zz_1)\Big),\\
d_{23}(\zz_2)-d_{23}(\zz_1) & = & (d_{33}-d_{22})g=
\beta\Big(v(\zz_2)-v(\zz_1)\Big).
\end{eqnarray*}
Obviously, a solution (unique up to an additive constant) is given
by
\[
d_{12}=\alpha u,\quad d_{23}=\beta v.
\]
Finally, equating in (\ref{eq D1})--(\ref{eq D3}) the entries 13,
we find:
\begin{eqnarray}
d_{13}(\zz_2)-d_{13}(\zz_1) & = & d_{23}(\zz_1)f-d_{12}(\zz_2)g
\label{eq for d13}\\
& = & \beta v(\zz_1)\Big(u(\zz_2)-u(\zz_1)\Big)- \alpha
u(\zz_2)\Big(v(\zz_2)-v(\zz_1)\Big).
\end{eqnarray}
Comparing this with (\ref{eq for a}), (\ref{eq for a'}), we see
that (\ref{D}) is proved.

{\bf Step 4. Equations relating the matrices} $C_{k,\ell,m}$ {\bf
and} $D(\zz)$. It remains to consider the equations (\ref{eq
CD1})--(\ref{eq CD3}). Denoting entries of the matrix $C$ by
$c_{ij}$, we see that these matrix equations are equivalent to the
following scalar ones:
\begin{eqnarray}
c_{12}+d_{12} & = & 0, \label{cd 12}\\
c_{23}+d_{23} & = & 0, \label{cd 23}\\
c_{13}+d_{13} & = & 0,\label{cd 13}
\end{eqnarray}
\begin{eqnarray}
(c_{33})_{k+1,\ell,m}-(c_{11})_{k,\ell,m}+d_{33}-d_{11} & = & 1,
\label{cd 11 33 k} \\
(c_{33})_{k,\ell+1,m}-(c_{11})_{k,\ell,m}+d_{33}-d_{11} & = & 1,
\label{cd 11 33 l}\\
(c_{33})_{k,\ell,m+1}-(c_{11})_{k,\ell,m}+d_{33}-d_{11} & = & 1.
\label{cd 11 33 m}
\end{eqnarray}
(In the last three equations we took into account that $d_{11}$,
$d_{33}$ are constants.) It is easy to see that the equations
(\ref{cd 12}), (\ref{cd 23}) are nothing but the constraint
equations (\ref{constr 1}), (\ref{constr 2}), respectively. We
show now that the remaining equations (\ref{cd 13})--(\ref{cd 11
33 m}) are not independent, but rather follow from the equations
of the $fgh$--system and the constraints (\ref{cd 12}), (\ref{cd
23}). We start with the last three equations, and prove the claim
for (\ref{cd 11 33 k}), since for other two everything is similar.
As in the proof of Proposition \ref{third constraint}, we write
the formulas here with $m=0$. Writing (\ref{cd 11 33 k}) in long
hand, using the ans\"atze (\ref{C}), (\ref{P}), (\ref{D}), we see
that it is equivalent to
\begin{eqnarray}\label{cd 11 33 to prove}
1-\alpha-\beta & = &
(k+1)\frac{1/h_0}{\widetilde{f}_0\widetilde{g}_0+
\widetilde{g}_0\widetilde{f}_3+\widetilde{f}_3\widetilde{g}_3}-
k\frac{1/h_0}{f_0g_0+g_0f_3+f_3g_3}\nonumber\\
& & +\ell \frac{1/\widetilde{h}_5}{\widetilde{f}_2\widetilde{g}_2+
\widetilde{g}_2\widetilde{f}_5+\widetilde{f}_5\widetilde{g}_5}-
\ell\frac{1/h_2}{f_2g_2+g_2f_5+f_5g_5}.
\end{eqnarray}
But this follows immediately from (\ref{constr 3 to prove}),
(\ref{C ansatz aux1}). Finally, we turn to (\ref{cd 13}).
Actually, since the entry 13 of the matrix $D$ is defined only up
to an additive constant, this equation is equivalent to the system
of the following three ones:
\begin{eqnarray}
(c_{13})_{k+1,\ell,m}-(c_{13})_{k,\ell,m}+\widetilde{d}_{13}-d_{13}
& = & 0,
\label{cd 13 k}\\
(c_{13})_{k,\ell+1,m}-(c_{13})_{k,\ell,m}+\widehat{d}_{13}-d_{13}
& = & 0,
\label{cd 13 l}\\
(c_{13})_{k,\ell,m+1}-(c_{11})_{k,\ell,m}+\bar{d}_{13}-d_{13} & =
& 0. \label{cd 13 m}
\end{eqnarray}
As usual, we restrict ourselves to the first one. Upon using the
equation (\ref{eq for d13}) and the constraints (\ref{cd 12}),
(\ref{cd 23}), we see that it is equivalent to
\begin{equation}\label{cd 13 to prove}
(c_{13})_{k+1,\ell,m}-(c_{13})_{k,\ell,m}+g_0(c_{12})_{k+1,\ell,m}
-f_0(c_{23})_{k,\ell,m}=0.
\end{equation}
Writing in long hand, in the representation with $m=0$, we see
that the terms proportional to $k+1$ and $k$ vanish identically,
while the vanishing of the terms proportional to $\ell$ is
equivalent to:
\[
\frac{1}{\widetilde{h}_2\widetilde{h}_5}\cdot
\frac{1-g_0/\widetilde{g}_5}
{\widetilde{f}_2\widetilde{g}_2+\widetilde{g}_2\widetilde{f}_5+
\widetilde{f}_5\widetilde{g}_5}=
\frac{1}{h_2h_5}\cdot\frac{1-f_0/f_2}{f_2g_2+g_2f_5+f_5g_5}.
\]
But this follows immediately from (\ref{C ansatz aux1}) and the
formulas
\[
\widetilde{g}_5=h_5\frac{g_0-\widetilde{g}_5}{h_0-\widetilde{h}_5},
\qquad \widetilde{g}_2=h_2\frac{g_0-g_2}{h_0-h_2},
\]
which are similar to (and follow from) the equations (\ref{C
ansatz aux3}), (\ref{C ansatz aux2}).

This finishes the proof of Theorem \ref{monodromy}. \qed


\end{document}

%% file: a_proofMV.tex
\begin{picture}(0,0)%
\epsfig{file=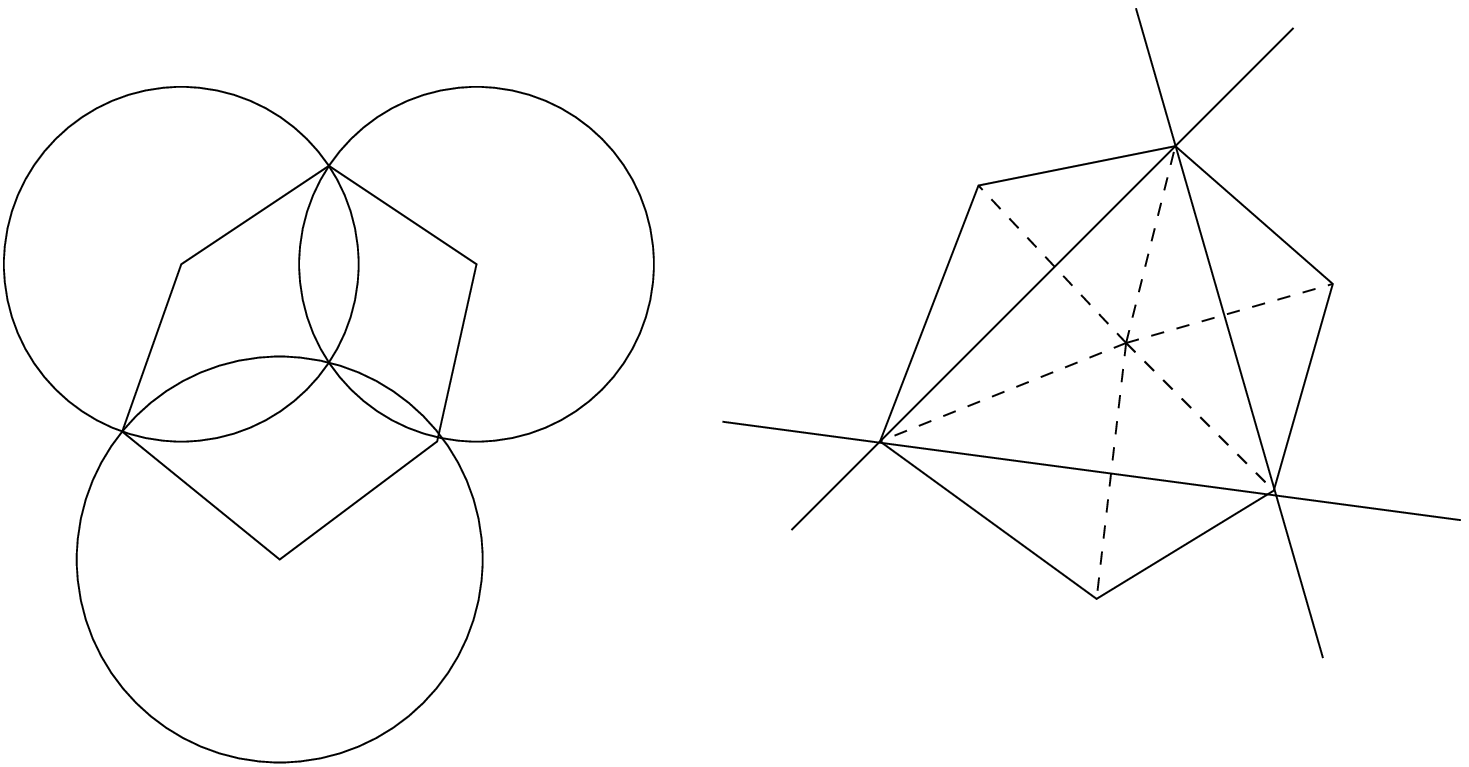}%
\end{picture}%
\setlength{\unitlength}{4144sp}%
\begingroup\makeatletter\ifx\SetFigFont\undefined%
\gdef\SetFigFont#1#2#3#4#5{%
  \reset@font\fontsize{#1}{#2pt}%
  \fontfamily{#3}\fontseries{#4}\fontshape{#5}%
  \selectfont}%
\fi\endgroup%
\begin{picture}(6681,3467)(757,-3021)
\put(3691,-1321){\makebox(0,0)[lb]{\smash{\SetFigFont{12}{14.4}{\rmdefault}{\mddefault}{\updefault}$\longrightarrow$}}}
\put(6526,-1366){\makebox(0,0)[lb]{\smash{\SetFigFont{12}{14.4}{\rmdefault}{\mddefault}{\updefault}$\gamma_2$}}}
\put(5941,-1816){\makebox(0,0)[lb]{\smash{\SetFigFont{12}{14.4}{\rmdefault}{\mddefault}{\updefault}$\gamma_1$}}}
\put(6301,-1501){\makebox(0,0)[lb]{\smash{\SetFigFont{12}{14.4}{\rmdefault}{\mddefault}{\updefault}$\gamma_2$}}}
\put(6076,-691){\makebox(0,0)[lb]{\smash{\SetFigFont{12}{14.4}{\rmdefault}{\mddefault}{\updefault}$\alpha_1$}}}
\put(6256,-466){\makebox(0,0)[lb]{\smash{\SetFigFont{12}{14.4}{\rmdefault}{\mddefault}{\updefault}$\alpha_1$}}}
\put(5671,-421){\makebox(0,0)[lb]{\smash{\SetFigFont{12}{14.4}{\rmdefault}{\mddefault}{\updefault}$\alpha_2$}}}
\put(5806,-601){\makebox(0,0)[lb]{\smash{\SetFigFont{12}{14.4}{\rmdefault}{\mddefault}{\updefault}$\alpha_2$}}}
\put(4951,-1186){\makebox(0,0)[lb]{\smash{\SetFigFont{12}{14.4}{\rmdefault}{\mddefault}{\updefault}$\beta_1$}}}
\put(5131,-1321){\makebox(0,0)[lb]{\smash{\SetFigFont{12}{14.4}{\rmdefault}{\mddefault}{\updefault}$\beta_1$}}}
\put(5086,-1726){\makebox(0,0)[lb]{\smash{\SetFigFont{12}{14.4}{\rmdefault}{\mddefault}{\updefault}$\beta_2$}}}
\put(5986,-1141){\makebox(0,0)[lb]{\smash{\SetFigFont{12}{14.4}{\rmdefault}{\mddefault}{\updefault}$p_\infty$}}}
\put(6076, 29){\makebox(0,0)[lb]{\smash{\SetFigFont{12}{14.4}{\rmdefault}{\mddefault}{\updefault}$w_2$}}}
\put(6121,-1636){\makebox(0,0)[lb]{\smash{\SetFigFont{12}{14.4}{\rmdefault}{\mddefault}{\updefault}$\gamma_1$}}}
\put(5131,-1501){\makebox(0,0)[lb]{\smash{\SetFigFont{12}{14.4}{\rmdefault}{\mddefault}{\updefault}$\beta_2$}}}
\put(4816,-286){\makebox(0,0)[lb]{\smash{\SetFigFont{12}{14.4}{\rmdefault}{\mddefault}{\updefault}$w_3$}}}
\put(2386,-1141){\makebox(0,0)[lb]{\smash{\SetFigFont{12}{14.4}{\rmdefault}{\mddefault}{\updefault}$w_0$}}}
\put(3061,-826){\makebox(0,0)[lb]{\smash{\SetFigFont{12}{14.4}{\rmdefault}{\mddefault}{\updefault}$w_1$}}}
\put(2971,-1726){\makebox(0,0)[lb]{\smash{\SetFigFont{12}{14.4}{\rmdefault}{\mddefault}{\updefault}$w_6$}}}
\put(2161,-2131){\makebox(0,0)[lb]{\smash{\SetFigFont{12}{14.4}{\rmdefault}{\mddefault}{\updefault}$w_5$}}}
\put(2116,-61){\makebox(0,0)[lb]{\smash{\SetFigFont{12}{14.4}{\rmdefault}{\mddefault}{\updefault}$w_2$}}}
\put(1351,-691){\makebox(0,0)[lb]{\smash{\SetFigFont{12}{14.4}{\rmdefault}{\mddefault}{\updefault}$w_3$}}}
\put(991,-1591){\makebox(0,0)[lb]{\smash{\SetFigFont{12}{14.4}{\rmdefault}{\mddefault}{\updefault}$w_4$}}}
\put(4366,-1681){\makebox(0,0)[lb]{\smash{\SetFigFont{12}{14.4}{\rmdefault}{\mddefault}{\updefault}$w_4$}}}
\put(6706,-1996){\makebox(0,0)[lb]{\smash{\SetFigFont{12}{14.4}{\rmdefault}{\mddefault}{\updefault}$w_6$}}}
\put(7021,-736){\makebox(0,0)[lb]{\smash{\SetFigFont{12}{14.4}{\rmdefault}{\mddefault}{\updefault}$w_1$}}}
\put(5446,-2356){\makebox(0,0)[lb]{\smash{\SetFigFont{12}{14.4}{\rmdefault}{\mddefault}{\updefault}$w_5$}}}
\end{picture}

%% file: a_proof15u.tex
\begin{picture}(0,0)%
\epsfig{file=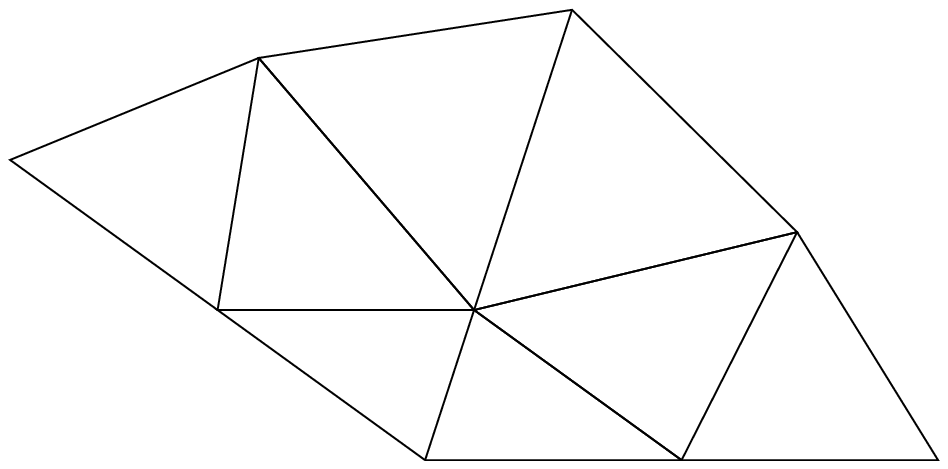}%
\end{picture}%
\setlength{\unitlength}{4144sp}%
\begingroup\makeatletter\ifx\SetFigFont\undefined%
\gdef\SetFigFont#1#2#3#4#5{%
  \reset@font\fontsize{#1}{#2pt}%
  \fontfamily{#3}\fontseries{#4}\fontshape{#5}%
  \selectfont}%
\fi\endgroup%
\begin{picture}(4479,2211)(1159,-4510)
\put(3061,-3976){\makebox(0,0)[lb]{\smash{\SetFigFont{12}{14.4}{\rmdefault}{\mddefault}{\updefault}$\psi_2$}}}
\put(3421,-4066){\makebox(0,0)[lb]{\smash{\SetFigFont{12}{14.4}{\rmdefault}{\mddefault}{\updefault}$\psi_1$}}}
\put(3601,-3886){\makebox(0,0)[lb]{\smash{\SetFigFont{12}{14.4}{\rmdefault}{\mddefault}{\updefault}$\psi_4$}}}
\put(4771,-3796){\makebox(0,0)[lb]{\smash{\SetFigFont{12}{14.4}{\rmdefault}{\mddefault}{\updefault}$\psi_3$}}}
\put(2926,-4246){\makebox(0,0)[lb]{\smash{\SetFigFont{12}{14.4}{\rmdefault}{\mddefault}{\updefault}$\psi_2$}}}
\put(1486,-3166){\makebox(0,0)[lb]{\smash{\SetFigFont{12}{14.4}{\rmdefault}{\mddefault}{\updefault}$\psi_6$}}}
\put(2026,-3571){\makebox(0,0)[lb]{\smash{\SetFigFont{12}{14.4}{\rmdefault}{\mddefault}{\updefault}$\phi_6$}}}
\put(2521,-3931){\makebox(0,0)[lb]{\smash{\SetFigFont{12}{14.4}{\rmdefault}{\mddefault}{\updefault}$\phi_2$}}}
\put(2341,-3661){\makebox(0,0)[lb]{\smash{\SetFigFont{12}{14.4}{\rmdefault}{\mddefault}{\updefault}$\phi_5$}}}
\put(2881,-3661){\makebox(0,0)[lb]{\smash{\SetFigFont{12}{14.4}{\rmdefault}{\mddefault}{\updefault}$\psi_5$}}}
\put(3286,-4336){\makebox(0,0)[lb]{\smash{\SetFigFont{12}{14.4}{\rmdefault}{\mddefault}{\updefault}$\psi_1$}}}
\put(4501,-4336){\makebox(0,0)[lb]{\smash{\SetFigFont{12}{14.4}{\rmdefault}{\mddefault}{\updefault}$\phi_3$}}}
\put(5131,-4336){\makebox(0,0)[lb]{\smash{\SetFigFont{12}{14.4}{\rmdefault}{\mddefault}{\updefault}$\psi_3$}}}
\put(2071,-2941){\makebox(0,0)[lb]{\smash{\SetFigFont{12}{14.4}{\rmdefault}{\mddefault}{\updefault}$\psi_6$}}}
\put(2431,-2986){\makebox(0,0)[lb]{\smash{\SetFigFont{12}{14.4}{\rmdefault}{\mddefault}{\updefault}$\psi_5$}}}
\put(3961,-4336){\makebox(0,0)[lb]{\smash{\SetFigFont{12}{14.4}{\rmdefault}{\mddefault}{\updefault}$\phi_1$}}}
\put(4501,-3661){\makebox(0,0)[lb]{\smash{\SetFigFont{12}{14.4}{\rmdefault}{\mddefault}{\updefault}$\psi_4$}}}
\put(4186,-4201){\makebox(0,0)[lb]{\smash{\SetFigFont{12}{14.4}{\rmdefault}{\mddefault}{\updefault}$\phi_4$}}}
\end{picture}

%% file: a_proof15v.tex
\begin{picture}(0,0)%
\epsfig{file=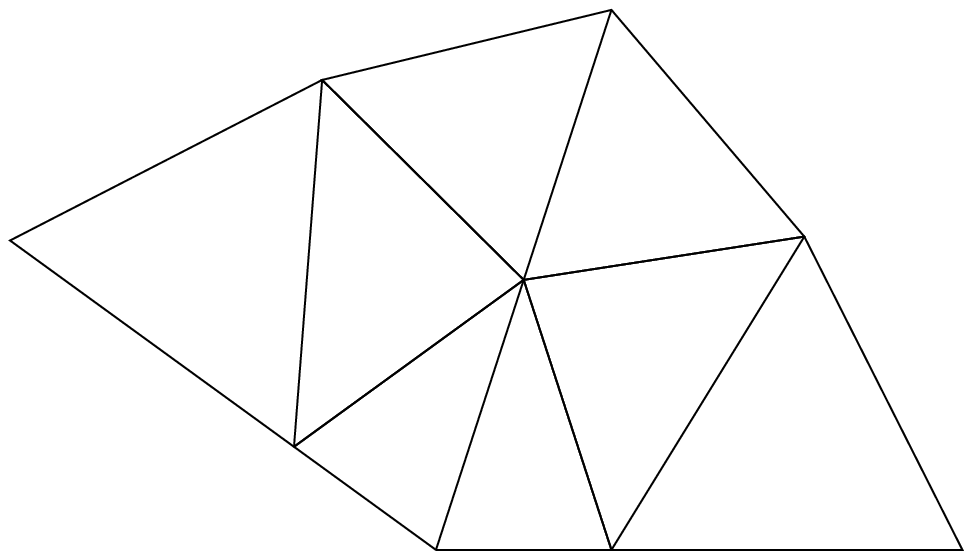}%
\end{picture}%
\setlength{\unitlength}{4144sp}%
\begingroup\makeatletter\ifx\SetFigFont\undefined%
\gdef\SetFigFont#1#2#3#4#5{%
  \reset@font\fontsize{#1}{#2pt}%
  \fontfamily{#3}\fontseries{#4}\fontshape{#5}%
  \selectfont}%
\fi\endgroup%
\begin{picture}(4596,2627)(2104,-5106)
\put(4591,-4921){\makebox(0,0)[lb]{\smash{\SetFigFont{12}{14.4}{\rmdefault}{\mddefault}{\updefault}$\psi_1$}}}
\put(4681,-3976){\makebox(0,0)[lb]{\smash{\SetFigFont{12}{14.4}{\rmdefault}{\mddefault}{\updefault}$\phi_4$}}}
\put(5446,-3841){\makebox(0,0)[lb]{\smash{\SetFigFont{12}{14.4}{\rmdefault}{\mddefault}{\updefault}$\psi_4$}}}
\put(5761,-4021){\makebox(0,0)[lb]{\smash{\SetFigFont{12}{14.4}{\rmdefault}{\mddefault}{\updefault}$\psi_3$}}}
\put(4276,-3796){\makebox(0,0)[lb]{\smash{\SetFigFont{12}{14.4}{\rmdefault}{\mddefault}{\updefault}$\phi_5$}}}
\put(3331,-3211){\makebox(0,0)[lb]{\smash{\SetFigFont{12}{14.4}{\rmdefault}{\mddefault}{\updefault}$\psi_6$}}}
\put(3691,-3301){\makebox(0,0)[lb]{\smash{\SetFigFont{12}{14.4}{\rmdefault}{\mddefault}{\updefault}$\psi_5$}}}
\put(4231,-4156){\makebox(0,0)[lb]{\smash{\SetFigFont{12}{14.4}{\rmdefault}{\mddefault}{\updefault}$\phi_2$}}}
\put(6121,-4921){\makebox(0,0)[lb]{\smash{\SetFigFont{12}{14.4}{\rmdefault}{\mddefault}{\updefault}$\phi_3$}}}
\put(5176,-4921){\makebox(0,0)[lb]{\smash{\SetFigFont{12}{14.4}{\rmdefault}{\mddefault}{\updefault}$\psi_3$}}}
\put(4276,-4921){\makebox(0,0)[lb]{\smash{\SetFigFont{12}{14.4}{\rmdefault}{\mddefault}{\updefault}$\psi_1$}}}
\put(2431,-3706){\makebox(0,0)[lb]{\smash{\SetFigFont{12}{14.4}{\rmdefault}{\mddefault}{\updefault}$\phi_6$}}}
\put(3196,-4246){\makebox(0,0)[lb]{\smash{\SetFigFont{12}{14.4}{\rmdefault}{\mddefault}{\updefault}$\psi_6$}}}
\put(3691,-4606){\makebox(0,0)[lb]{\smash{\SetFigFont{12}{14.4}{\rmdefault}{\mddefault}{\updefault}$\psi_2$}}}
\put(4006,-4831){\makebox(0,0)[lb]{\smash{\SetFigFont{12}{14.4}{\rmdefault}{\mddefault}{\updefault}$\psi_2$}}}
\put(4501,-4246){\makebox(0,0)[lb]{\smash{\SetFigFont{12}{14.4}{\rmdefault}{\mddefault}{\updefault}$\phi_1$}}}
\put(4906,-4741){\makebox(0,0)[lb]{\smash{\SetFigFont{12}{14.4}{\rmdefault}{\mddefault}{\updefault}$\psi_4$}}}
\put(3601,-4291){\makebox(0,0)[lb]{\smash{\SetFigFont{12}{14.4}{\rmdefault}{\mddefault}{\updefault}$\psi_5$}}}
\end{picture}

%% file: a_proof15w.tex
\begin{picture}(0,0)%
\epsfig{file=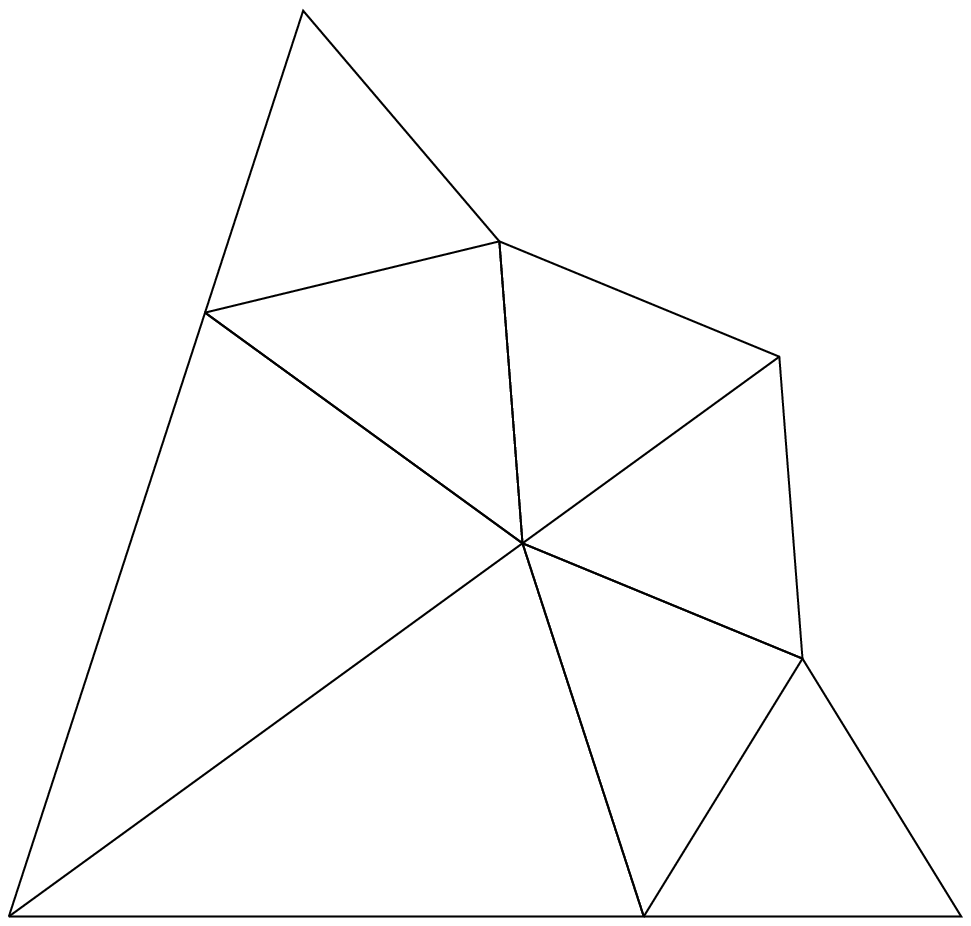}%
\end{picture}%
\setlength{\unitlength}{4144sp}%
\begingroup\makeatletter\ifx\SetFigFont\undefined%
\gdef\SetFigFont#1#2#3#4#5{%
  \reset@font\fontsize{#1}{#2pt}%
  \fontfamily{#3}\fontseries{#4}\fontshape{#5}%
  \selectfont}%
\fi\endgroup%
\begin{picture}(4596,4374)(2779,-7078)
\put(4186,-3211){\makebox(0,0)[lb]{\smash{\SetFigFont{12}{14.4}{\rmdefault}{\mddefault}{\updefault}$\psi_6$}}}
\put(4771,-3841){\makebox(0,0)[lb]{\smash{\SetFigFont{12}{14.4}{\rmdefault}{\mddefault}{\updefault}$\phi_6$}}}
\put(3916,-4066){\makebox(0,0)[lb]{\smash{\SetFigFont{12}{14.4}{\rmdefault}{\mddefault}{\updefault}$\psi_6$}}}
\put(4861,-4111){\makebox(0,0)[lb]{\smash{\SetFigFont{12}{14.4}{\rmdefault}{\mddefault}{\updefault}$\phi_5$}}}
\put(4006,-4336){\makebox(0,0)[lb]{\smash{\SetFigFont{12}{14.4}{\rmdefault}{\mddefault}{\updefault}$\psi_5$}}}
\put(4906,-4966){\makebox(0,0)[lb]{\smash{\SetFigFont{12}{14.4}{\rmdefault}{\mddefault}{\updefault}$\psi_5$}}}
\put(3781,-4561){\makebox(0,0)[lb]{\smash{\SetFigFont{12}{14.4}{\rmdefault}{\mddefault}{\updefault}$\psi_2$}}}
\put(4816,-5326){\makebox(0,0)[lb]{\smash{\SetFigFont{12}{14.4}{\rmdefault}{\mddefault}{\updefault}$\psi_2$}}}
\put(3196,-6451){\makebox(0,0)[lb]{\smash{\SetFigFont{12}{14.4}{\rmdefault}{\mddefault}{\updefault}$\phi_2$}}}
\put(3376,-6856){\makebox(0,0)[lb]{\smash{\SetFigFont{12}{14.4}{\rmdefault}{\mddefault}{\updefault}$\phi_1$}}}
\put(4996,-5641){\makebox(0,0)[lb]{\smash{\SetFigFont{12}{14.4}{\rmdefault}{\mddefault}{\updefault}$\psi_1$}}}
\put(5356,-6856){\makebox(0,0)[lb]{\smash{\SetFigFont{12}{14.4}{\rmdefault}{\mddefault}{\updefault}$\psi_1$}}}
\put(5446,-5641){\makebox(0,0)[lb]{\smash{\SetFigFont{12}{14.4}{\rmdefault}{\mddefault}{\updefault}$\psi_4$}}}
\put(5716,-6631){\makebox(0,0)[lb]{\smash{\SetFigFont{12}{14.4}{\rmdefault}{\mddefault}{\updefault}$\psi_4$}}}
\put(6166,-5911){\makebox(0,0)[lb]{\smash{\SetFigFont{12}{14.4}{\rmdefault}{\mddefault}{\updefault}$\phi_4$}}}
\put(6391,-6136){\makebox(0,0)[lb]{\smash{\SetFigFont{12}{14.4}{\rmdefault}{\mddefault}{\updefault}$\phi_3$}}}
\put(6031,-6856){\makebox(0,0)[lb]{\smash{\SetFigFont{12}{14.4}{\rmdefault}{\mddefault}{\updefault}$\psi_3$}}}
\put(6751,-6856){\makebox(0,0)[lb]{\smash{\SetFigFont{12}{14.4}{\rmdefault}{\mddefault}{\updefault}$\psi_3$}}}
\end{picture}

%% file: a_proof17u.tex
\begin{picture}(0,0)%
\epsfig{file=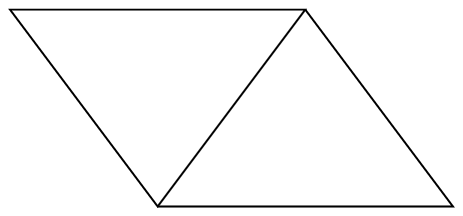}%
\end{picture}%
\setlength{\unitlength}{4144sp}%
\begingroup\makeatletter\ifx\SetFigFont\undefined%
\gdef\SetFigFont#1#2#3#4#5{%
  \reset@font\fontsize{#1}{#2pt}%
  \fontfamily{#3}\fontseries{#4}\fontshape{#5}%
  \selectfont}%
\fi\endgroup%
\begin{picture}(2127,1365)(361,-691)
\put(1036,-331){\makebox(0,0)[lb]{\smash{\SetFigFont{12}{14.4}{\rmdefault}{\mddefault}{\updefault}$\phi_2$}}}
\put(1036,-691){\makebox(0,0)[lb]{\smash{\SetFigFont{12}{14.4}{\rmdefault}{\mddefault}{\updefault}$u_0$}}}
\put(2476,-691){\makebox(0,0)[lb]{\smash{\SetFigFont{12}{14.4}{\rmdefault}{\mddefault}{\updefault}$u_1$}}}
\put(1846,434){\makebox(0,0)[lb]{\smash{\SetFigFont{12}{14.4}{\rmdefault}{\mddefault}{\updefault}$u_3$}}}
\put(361,479){\makebox(0,0)[lb]{\smash{\SetFigFont{12}{14.4}{\rmdefault}{\mddefault}{\updefault}$u_2$}}}
\put(676,209){\makebox(0,0)[lb]{\smash{\SetFigFont{12}{14.4}{\rmdefault}{\mddefault}{\updefault}$\psi_2$}}}
\put(1306,-421){\makebox(0,0)[lb]{\smash{\SetFigFont{12}{14.4}{\rmdefault}{\mddefault}{\updefault}$\phi_1$}}}
\put(2161,-421){\makebox(0,0)[lb]{\smash{\SetFigFont{12}{14.4}{\rmdefault}{\mddefault}{\updefault}$\psi_1$}}}
\put(1711,119){\makebox(0,0)[lb]{\smash{\SetFigFont{12}{14.4}{\rmdefault}{\mddefault}{\updefault}$\chi_1$}}}
\put(1486,254){\makebox(0,0)[lb]{\smash{\SetFigFont{12}{14.4}{\rmdefault}{\mddefault}{\updefault}$\chi_2$}}}
\end{picture}

%% file: a_proof17v.tex
\begin{picture}(0,0)%
\epsfig{file=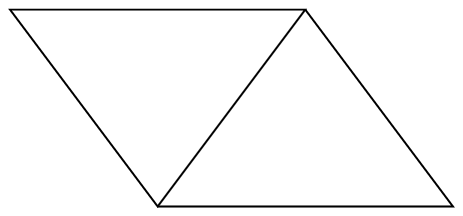}%
\end{picture}%
\setlength{\unitlength}{4144sp}%
\begingroup\makeatletter\ifx\SetFigFont\undefined%
\gdef\SetFigFont#1#2#3#4#5{%
  \reset@font\fontsize{#1}{#2pt}%
  \fontfamily{#3}\fontseries{#4}\fontshape{#5}%
  \selectfont}%
\fi\endgroup%
\begin{picture}(2127,1365)(361,-691)
\put(1036,-331){\makebox(0,0)[lb]{\smash{\SetFigFont{12}{14.4}{\rmdefault}{\mddefault}{\updefault}$\chi_2$}}}
\put(1036,-691){\makebox(0,0)[lb]{\smash{\SetFigFont{12}{14.4}{\rmdefault}{\mddefault}{\updefault}$v_0$}}}
\put(2476,-691){\makebox(0,0)[lb]{\smash{\SetFigFont{12}{14.4}{\rmdefault}{\mddefault}{\updefault}$v_1$}}}
\put(1846,434){\makebox(0,0)[lb]{\smash{\SetFigFont{12}{14.4}{\rmdefault}{\mddefault}{\updefault}$v_3$}}}
\put(361,479){\makebox(0,0)[lb]{\smash{\SetFigFont{12}{14.4}{\rmdefault}{\mddefault}{\updefault}$v_2$}}}
\put(676,209){\makebox(0,0)[lb]{\smash{\SetFigFont{12}{14.4}{\rmdefault}{\mddefault}{\updefault}$\phi_2$}}}
\put(1306,-421){\makebox(0,0)[lb]{\smash{\SetFigFont{12}{14.4}{\rmdefault}{\mddefault}{\updefault}$\chi_1$}}}
\put(2161,-421){\makebox(0,0)[lb]{\smash{\SetFigFont{12}{14.4}{\rmdefault}{\mddefault}{\updefault}$\phi_1$}}}
\put(1711,119){\makebox(0,0)[lb]{\smash{\SetFigFont{12}{14.4}{\rmdefault}{\mddefault}{\updefault}$\psi_1$}}}
\put(1486,254){\makebox(0,0)[lb]{\smash{\SetFigFont{12}{14.4}{\rmdefault}{\mddefault}{\updefault}$\psi_2$}}}
\end{picture}